\documentclass[a4paper,12pt]{article}%
 \usepackage[margin=1in,footskip=0.5in]{geometry}
\usepackage[bottom]{footmisc}
\usepackage{setspace}

\usepackage{amsthm}
\usepackage{amsmath}
\usepackage{amssymb,amsmath,amsthm,amsfonts,mathrsfs,xcolor}
\usepackage{natbib}
\usepackage[colorlinks,citecolor=blue,urlcolor=blue,filecolor=blue,backref=page]{hyperref}
\usepackage{graphicx}

\theoremstyle{plain}
\newtheorem{theorem}{Theorem}
 \newtheorem{remark}{Remark}
 \newtheorem{corollary}{Corollary}
 \newtheorem{lemma}{Lemma}
 \newtheorem{definition}{Definition}
 \newtheorem*{definition*}{Definition}
 
\usepackage{times}

\usepackage{amssymb,amsmath,amsthm,amsfonts,bbm,bm,mathrsfs,xcolor}
 \usepackage{float}  
 \usepackage{subfigure}  
 \usepackage{xr}

\usepackage{bigints}
\DeclareMathOperator*{\argmin}{arg\,min}
\DeclareMathOperator*{\argmax}{arg\,max}

\def\m{\mathcal}
\def\mb{\mathbb}

\def\dd{{\rm d}}

\def\wt{\widetilde}
\def\wh{\widehat}
\def\ov{\overline}

\newcommand{\mnorm}[1]{{\vert\kern-0.25ex\vert\kern-0.25ex\vert #1 
    \vert\kern-0.25ex\vert\kern-0.25ex\vert}}
\newcommand{\bmnorm}[1]{{\big\vert\kern-0.25ex\big\vert\kern-0.25ex\big\vert #1 
    \big\vert\kern-0.25ex\big\vert\kern-0.25ex\big\vert}}
\newcommand{\Bmnorm}[1]{{\Big\vert\kern-0.25ex\Big\vert\kern-0.25ex\Big\vert #1 
    \Big\vert\kern-0.25ex\Big\vert\kern-0.25ex\Big\vert}}
\newcommand{\bbmnorm}[1]{{\bigg\vert\kern-0.25ex\bigg\vert\kern-0.25ex\bigg\vert #1 
    \bigg\vert\kern-0.25ex\bigg\vert\kern-0.25ex\bigg\vert}}
\newcommand{\BBmnorm}[1]{{\Bigg\vert\kern-0.25ex\Bigg\vert\kern-0.25ex\Bigg\vert #1 
    \Bigg\vert\kern-0.25ex\Bigg\vert\kern-0.25ex\Bigg\vert}}

\newcommand{\be}{\begin{equs}}
\newcommand{\ee}{\end{equs}}
 
\newcommand{\keywords}[1]
{
  \small	
  \textbf{\textit{Keywords---}} #1
}

  \title{On the Computational Complexity of Metropolis-Adjusted Langevin Algorithms for Bayesian Posterior Sampling}

\author{Rong Tang\textsuperscript{$\ast$} and Yun Yang\textsuperscript{$\dagger$}}

\date{\textsuperscript{$\ast$}Department of Mathematics, The Hong Kong University of Science and Technology\\
\textsuperscript{$\dagger$}Department of Statistics, University of Illinois Urbana-Champaign}

\begin{document}
\maketitle

\begin{abstract}
In this paper, we examine the computational complexity of sampling from a Bayesian posterior (or pseudo-posterior) using the Metropolis-adjusted Langevin algorithm (MALA). MALA first employs a discrete-time Langevin SDE to propose a new state, and then adjusts the proposed state using Metropolis-Hastings rejection. Most existing theoretical analyses of MALA rely on the smoothness and strong log-concavity properties of the target distribution, which are often lacking in practical Bayesian problems.
Our analysis hinges on statistical large sample theory, which constrains the deviation of the Bayesian posterior from being smooth and log-concave in a very specific way. In particular, we introduce a new technique for bounding the mixing time of a Markov chain with a continuous state space via the $s$-conductance profile, offering improvements over existing techniques in several aspects.
By employing this new technique, we establish the optimal parameter dimension dependence of $d^{1/3}$ and condition number dependence of $\kappa$ in the non-asymptotic mixing time upper bound for MALA after the burn-in period, under a standard Bayesian setting where the target posterior distribution is close to a $d$-dimensional Gaussian distribution with a covariance matrix having a condition number $\kappa$. We also prove a matching mixing time lower bound for sampling from a multivariate Gaussian via MALA to complement the upper bound.

\end{abstract}


\section{Introduction}

 Bayesian inference gains significant popularity during the last two decades due to the advance in modern computing power.
 As a method of statistical analysis based on probabilistic modelling, Bayesian inference allows natural uncertainty quantification on the unknown parameters via a posterior distribution. In the classical Bayesian framework, the data $X^{(n)}=\{X_1,\ldots,X_n\}$ is assumed to consist of i.i.d.~samples generated from a probability distribution $p(X\,|\,\theta)$ depending on an unknown parameter $\theta$ in parameter space $\Theta\subset \mb R^d$. Domain knowledge and prior beliefs can be characterized by a probability distribution $\pi(\theta)$ over $\Theta$ called prior (distribution), which is then updated into a posterior (distribution) $p(\theta\,|\,X^{(n)})$ by multiplying with the likelihood function 
 $$
 \m L_n(\theta;\,X^{(n)}):\,=\prod_{i=1}^n p(X_i\,|\,\theta)
 $$ 
 evaluated on the observed data $X^{(n)}$ using the Bayes theorem.  The classical Bayesian framework relies on the likelihood formulation, which hinders its use in problems where the data generating model is hard to fully specify or is not our primary interest. The pseudo-posterior~\citep{JMLR:v17:15-290,ghosh2020general} idea provides a more general probabilistic inference framework to alleviate this restriction by replacing the negative log-likelihood function in the Bayesian posterior with a criterion function. For example, when applied to risk minimization problems, the so-called Gibbs posteriors~\citep{bhattacharya2020gibbs,syring2020gibbs} use the (scaled) empirical risk function as the criterion function, thus avoiding imposing restrictive assumptions on the statistical model through a fully specified likelihood function.

Despite the conceptual appeal of Bayesian inference,  its practical implementation is a notoriously difficult computational problem. For example, the posterior $p(\theta\,|\,X^{(n)})$ involves a normalisation constant that can be expressed as a multidimensional integral 
$$
\int_\Theta \m L_n(\theta;\,X^{(n)}) \,\pi(\theta)\,\dd\theta. 
$$
This integral is usually analytically intractable and hard to numerically approximate, especially when the parameter dimension $d$ is high. Different from those numerical methods for directly computing the normalisation constant, the Markov chain Monte Carlo (MCMC) algorithm~\citep{hastings1970monte,geman1984stochastic,robert2004monte} constructs a Markov chain, whose simulation only requires evaluations of the likelihood ratio under a pair of parameters, such that its stationary distribution matches the target posterior distribution. Thus, MCMC provides an appealing alternative for Bayesian computation by turning the integration problem into a sampling problem that does not require computing the normalisation constant. Despite its popularity, the theoretical analysis of the computational efficiency of MCMC algorithms is mostly carried out for smooth and log-concave target distributions, and is comparatively rare in the Bayesian literature where a (pseudo-)posterior can be non-smooth and non-log-concave. 
In addition,  precise characterizations of the computational complexity (or mixing time) and its dependence on the parameter dimension $d$ for commonly used MCMC algorithms are important for guiding their practical designs and use.  

{

A widely used MCMC algorithm for sampling from Bayesian posteriors is the Gibbs sampler, which generates samples from a multivariate distribution by iteratively sampling each variable from its conditional distribution, given all other variables. The Gibbs sampler is particularly efficient for Bayesian models with closed-form conditional distributions under conjugate priors. A recent theoretical study by~\cite{ascolani2023dimensionfree} provides a dimension-free mixing time bound for the Gibbs sampler when applied to certain high-dimensional Bayesian hierarchical models. However, it is important to note that each iteration of their algorithm involves the sequential sampling of each dimension of the parameter from its corresponding full conditional distribution. This means that the total number of sampling steps required for the Gibbs sampler to converge is at least linear in the parameter dimension, which is larger than our sub-linear $d^{1/3}$ scaling of the needed sampling steps for MALA. On the other hand, the per-step cost of MALA can be linear in $d$ because of gradient computation, while that of Gibbs sampling can be much lower, especially under weak dependence (although in the worst case, computing each conditional distribution may also require $O(d)$ complexity). On a separate note, we would like to mention that although MALA has a per-iteration cost linear in $d$ to compute the gradient, the computation across different dimensions can be parallelized. In contrast, Gibbs sampling must sequentially scan over all its components and cannot be made parallel in order to maintain the detailed balance property. Additionally, the high efficiency of the Gibbs sampler often relies on the use of conjugate priors that facilitate closed-form conditional distributions. However, for complex Bayesian models, such a conjugate prior may not exist, as is the case in Bayesian quantile regression, discussed in~\cite{YU2001437}, or linear regression with heavy-tailed noise (like Student's t-distributions).   Moreover, there are situations where people tend to use specific non-conjugate priors for particular reasons. For example,  sparsity-induced priors such as the spike and slab priors (with heavy-tailed slabs) are widely used in regression analysis for facilitating variable selection. In these complicated scenarios, one might have to resort to using MALA or, more broadly, the Metropolis-Hastings (MH) algorithm, to draw samples from the Bayesian posterior.

On the other hand, the Metropolis-Hastings (MH) algorithm provides a more flexible alternative. An MH algorithm produces samples by proposing and then accepting or rejecting these proposals based on a specified acceptance criterion. A key advantage of the MH algorithm is its ability to handle Bayesian (pseudo-)posterior distributions without requiring explicit knowledge of the normalization constant or the full conditional distributions.}
One of the most popular MH algorithms is the Metropolis random walk (MRW), a zeroth-order method that queries the value of the target density ratio under two points per iteration. \cite{JMLR:v20:19-306} shows that for a log-concave and smooth target density, the $\varepsilon$-mixing time in total variation distance (the number of iterations required to converge
to an $\varepsilon$-neighborhood of stationary distribution in the total variation distance) for MRW is at most $\m O\big(d\log ({1}/{\varepsilon})\big)$. On the other hand, the $\m O(d)$ scaling limit of \cite{gelman1997weak} suggests that their linear dependence on dimension $d$ is optimal. For a class of Bayesian pseudo-posteriors that can be non-smooth and non-log-concave, it has been shown in~\cite{10.2307/30243694} that as the sample size $n$ grows to infinity while the parameter dimension $d$ does not grow too quickly relative to $n$ so that the pseudo-posterior satisfies a Bernstein-von Mises (asymptotic normality) result, then MRW for sampling from the target pseudo-posterior constrained on an approximate compact set with a warm start has an asymptotic total variation $\varepsilon$-mixing time upper bound as $\m O_p\big(d^2 \log (1/{\varepsilon})\big)$.

{
Another prominent class of MH algorithms is the Metropolis-adjusted Langevin algorithm (MALA), which utilizes additional gradient information about the target density. Although this approach requires computing the gradient and can be costlier than zeroth-order methods that only use function evaluations, the development of automatic differentiation tools~\citep{paszke2017automatic, margossian2019review} has simplified this task for many explicit and smooth densities. These tools make the computational demands for gradient computation comparable to those for evaluating the density itself. Furthermore, it has been demonstrated that MALA tends to have a lower mixing time in comparison to the MRW.} For example,~\cite{pmlr-v134-chewi21a} show that if the negative log-density (will be referred to as potential) of the target distribution is twice continuously differentiable and strongly convex, then the $\varepsilon$-mixing time in $\chi^2$ divergence for MALA with a warm start scales as $\Theta(d^{1/2})$ modulo polylogarithmic factors in $\varepsilon$.  Additionally,~\cite{roberts1998optimal} and~\cite{pmlr-v134-chewi21a} show that the optimal dimension dependence for MALA is $d^{1/3}$ for some product measures satisfying stringent conditions like the standard Gaussian. However, for Bayesian (pseudo-)posteriors, it is common that the smoothness and strong convexity properties of the log-density assumed in literature are not satisfied. {
For instance, consider Bayesian quantile regression with a quantile level $\tau$. Given a dataset $X^{(n)}=\{X_i=(\widetilde X_i, Y_i)\}_{i=1}^n$ consisting of covariates and response variables, the posterior distribution then takes the form of $\pi_{n}(\theta|X^{(n)})\propto \exp\big(-\sum_{i=1}^n (Y_i-\widetilde X_i^T\theta)(\tau-\mathbf{1}(Y_i<\widetilde X_i^T\theta))\big)\,\pi(\theta)$, where $\mathbf{1}(\cdot)$ denotes the indicator function. An important feature of this example is that the resulting Bayesian posterior is neither differentiable owing to the discontinuity introduced by the indicator function, nor strongly log-concave.} For such non-differentiable densities, we slightly extend the MALA by using any subgradient to replace the gradient in its algorithm formulation. Theoretically, it is natural to investigate: 
 {
\begin{quote}\normalsize
  What is the optimal dimension (and condition number) dependence when using MALA to sample from a possibly non-smooth and non-log-concave (pseudo-)posterior density, in light of the asymptotic Gaussian nature of the posterior as predicted by statistical large sample theory?
\end{quote}
Moreover, it would be insightful to determine to what extent we can diverge from a Gaussian distribution while preserving the dimension dependence as sampling from a Gaussian distribution, and how various factors, such as the dimensionality, sample size and density smoothness, affect the deviance of the posterior from the Gaussian distribution.}

\smallskip

\noindent {\bf Our contributions.}
 In this work, we show an upper bound on the $\varepsilon$-mixing time of MALA for sampling from a class of possibly non-smooth and non-log-concave distributions with non-product forms (c.f.~Condition A for a precise definition) with an $M_0$-warm start (defined in~Section~\ref{sec:problem}) as $\m O\big(\max\big\{d^{1/3}\log (\varepsilon^{-1}\log M_0),\,\log M_0\big\}\big)$, which matches (up to logarithmic terms in $(M_0,\varepsilon)$) the lower bound result proved in~\cite{pmlr-v134-chewi21a} that the mixing time of MALA for the standard Gaussian is at least $\m O(d^{1/3})$. 
  Specially, our condition requires the target distribution (after proper rescaling by the sample size $n$) to be close to a multivariate Gaussian subject to small perturbations. 
  {
  We verify that a wide class of Gibbs posteriors~\citep{bhattacharya2020gibbs,syring2020gibbs}, including conventional Bayesian posteriors defined through likelihood functions, meets our condition under a minimal set of assumptions. In particular, our theory provides an explicit upper bound condition on the growth of parameter dimension $d$ relative to sample size $n$, stated in a non-asymptotic manner, that is, $d\leq c\frac{n^{\kappa_1}}{\log n}$, where $\kappa_1$ depends on the regularity of the density function (c.f.~Theorem~\ref{th:Gibbsmixing}). Specifically, for less smooth density functions, a smaller dimension $d$ is necessary to maintain the $d^{1/3}$ scaling of the mixing time guarantee, which is also supported by our numerical results in Section~\ref{sec:num}.}

{
In addition, our result illustrates that the mixing time of MALA exhibits a linear dependence on the condition number $\kappa$ of the covariance matrix (which may have a polynomial dependence on the dimension in some ill-conditioned cases) of the approximating multivariate Gaussian. Our bound matches the mixing time scaling of Gaussian targets with condition number $\kappa$, and is therefore optimal. For the sake of completeness, we derive a matching lower bound in Appendix~\ref{app:lower}. In our lower bound analysis, we extend the proof of Theorem 1 in~\cite{pmlr-v134-chewi21a}, which primarily focuses on a standard Gaussian target distribution. In addition, we also carefully keep track of the dependence on the condition number in our derivation, which allows us to establish a lower bound that explicitly demonstrates a linear dependence on the condition number and also matches with our upper bound. }

  It is worthwhile mentioning that our Condition A does not require the distance between the target posterior and the multivariate  Gaussian distribution to vanish as $n$ tends to infinity; while in the context of Bayesian posteriors, these distances indeed decay to zero under minimal assumptions on the statistical model. Therefore, our mixing time result is more generally applicable to problems beyond Bayesian posterior sampling, for example, to optimization of approximately convex functions via simulated annealing~\citep{belloni2015escaping}, where the target distribution can deviate from being smooth and strongly log-concave by a finite amount. In such settings, the computational complexity of sampling
  algorithms scales as $\m O(d^{1/3})$ with the variable dimension $d$ under reasonably good initialization while that of a wide class of gradient-based optimization algorithms may scale exponentially~\citep{ma2019sampling}.

  Our result on the $\m O(d^{1/3})$ dimension dependence for the mixing time of MALA after the burn-in period for the perturbed Gaussian class strengthens our understanding of sampling from non-smooth and non-log-concave distributions. It also partly fills the gap between the optimal $d^{1/3}$ mixing time for a class of sufficiently regular product distributions derived from the scaling limit approach in~\cite{roberts1998optimal} and the $d^{1/2}$ lower bound on the class of all log-smooth and strongly log-concave distributions obtained in~\cite{pmlr-v134-chewi21a}, by identifying a much larger class of distributions of practical interest that attain the optimal $d^{1/3}$ dimension dependence. Moreover, we introduce a somewhat more general average conductance argument based on the $s$-conductance profile in Section~\ref{mixingboundscp} to improve the warming parameter dependence without deteriorating the dimension dependence. More specifically, our mixing time upper bound improves upon existing results~\citep[e.g.][]{pmlr-v134-chewi21a} in the dependence on the warming parameter $M_0$ from logarithmic to doubly logarithmic (the $\log\log(M_0)$ term in Theorem~\ref{thmala}) when $\log M_0\leq d^{\frac{1}{3}}$, by adapting the $s$-conductance profile and the log-isoperimetric inequality device~\citep{JMLR:v21:19-441}, or more generally, the log-Sobolev inequality device~\citep{lovasz1999faster,kannan2006blocking}, to our target distribution class. {
  Our constraint of $d^{\frac{1}{3}}$ on $\log M_0$ can be overly strong for general target distributions in practice.  For instance, in the case of distributions possessing product forms, such as a pair of isotropic Gaussians with varying means, $\log M_0$ tends to increase linearly with the dimension $d$.  However, for Bayesian posterior with smooth density,  we may leverage its asymptotic distribution to construct more effective warm starts (c.f. Lemma~\ref{lemma:Warm} and Corollary~\ref{cor:smoothloss}). }In addition, we study a variant of MALA where the (sub-)gradient vector in the Langevin SDE is preconditioned by a matrix for capturing the local geometry, for example, the Fisher information matrix in the context of Bayesian posterior sampling, and  we illustrate in our Corollaries~\ref{cor:smoothloss} and~\ref{co:quantile} that MALA with suitable preconditioning may improve the convergence of the sampling algorithm even though the target density is non-differentiable.

  Our analysis is motivated by the statistical large sample theory suggesting the Bayesian posterior to be close to a multivariate Gaussian. We develop mixing time bounds of MALA for sampling from
  general Gibbs posteriors (possibly with increasing parameter dimension and non-smooth criterion function) by establishing  non-asymptotic Bernstein-von Mises results, applying techniques from empirical process theory, including chaining, peeling, and localization. Due to the delicate analysis in our mixing time upper bound proof that utilizes the explicit form of Gaussian distributions for bounding the acceptance probability in each step of MALA, we obtain a better dimension dependence of $d^{1/3}$ than the $d^{1/2}$ dependence derived for general smooth and log-concave densities. In addition, by utilizing  our $s$-conductance profile technique,  we can obtain a mixing time upper bound for sampling from the original Bayesian posterior instead of a truncated version considered in~\cite{10.2307/30243694}.


\noindent {\bf Organization.} The rest of the paper is organized as follows. In Section~\ref{sec:prelim}, we describe the background and formally formulate the theoretical problem of analyzing the computational complexity of MALA for Bayesian posterior sampling that is addressed in this work. In Section~\ref{mixingboundscp}, we briefly review some common concepts and existing techniques for analyzing the computational complexity (in terms of mixing time) of a Markov chain, and introduce our improved technique based on $s$-conductance profile.
In Section~\ref{sec:MALA_mixtime}, we apply the generic technique developed in Section~\ref{mixingboundscp} to analyze MALA for Bayesian posterior sampling.
In Section~\ref{sec:application}, we specialize the general mixing bound of MALA to the class of Gibbs posteriors, and apply it to both Gibbs posteriors with smooth and non-smooth loss functions.
Section~\ref{sec:proof_sketch} sketches the main ideas in proving the MALA mixing time bound and discuss some main differences with existing proofs.
Some numerical studies are provided in Section~\ref{sec:num}, where we empirically compare the convergence of MALA and MRW. All proofs and technical details are deferred to the appendices in the supplementary material.

\noindent {\bf Notation.}  
 For two real numbers, we use $a\wedge b$ and $a\vee b$ to denote the maximum and minimum between $a$ and $b$. For two distributions $p$ and $q$, we use $\|p-q\|_{  \rm TV}=\frac{1}{2}\int|p(x)-q(x)|\, \dd x$ to denote their the total variation distance and $\chi^2(p,\,q)$ to denote their $\chi^2$ divergence.  We use $\|\cdot\|_p$ to denote the usual vector $\ell_p$ norm, and suppress the subscript when $p=2$. We use $\mathbf{0}_d$ to denote the $d$-dimensional all zero vector, and $B_r(x)$ to denote the closed ball centered at $x$ with radius $r$ (under the $\ell_2$ distance) in the Euclidean space; in particular, we use $B_r^d$ to denote $B_r(\mathbf{0}_d)$ when no ambiguity may arise. We use $\mb S^{d}=\big\{x\in \mb R^{d+1}:\, \|x\|=1\big\}$ to denote the $d$-dimensional sphere.  We use $N_d(\mu, \Sigma)$ to denote the $d$-dimensional multivariate Gaussian distribution with mean vector $\mu\in \mathbb{R}^d$ and covariance matrix $\Sigma\in \mathbb{R}^{d\times d}$, and d suppress the subscript when $d=1$. We use $\m P(K)$ to denote the set of probability measures on a set $K$. For a function $f:\mb R^d\to \mb R$, we use $\nabla f(x)$ to denote the $d$-dimensional gradient vector of $f$ at $x$ and ${\rm Hess}(f(x))$ to denote the Hessian matrix of $f$ at $x$. For a matrix $J$, we use $\mnorm{J}_{ \rm  op}$ and $\mnorm{J}_{ \rm F}$ to denote its operator norm and Frobenius norm respectively, and use $\lambda_{\max}(J)$ and $\lambda_{\min}(J)$ to denote the maximal and minimal eigenvalues of $J$. Throughout, $C$, $c$, $C_0$, $c_0$, $C_1$, $c_1$, \ldots are generically used to denote positive constants independent of $n,d$ whose values might change from one line to another.

 \section{Background and Problem Setup}\label{sec:prelim}
 We first review the Bayesian (pseudo-)posterior framework and the Metropolis-adjusted Langevin algorithm (MALA).
 After that, we discuss an extension of MALA to handle the case where the target density is non-smooth by using the subgradient to replace the gradient and formulate the theoretical problem to be addressed in this work.

 \subsection{Bayesian pseudo-posterior}\label{pseudo-posterior}
 A standard Bayesian model consists of a prior distribution (density) $\pi(\theta)$ over parameter space $\Theta\subset\mb R^d$ as the marginal distribution of the parameter $\theta$ and a sampling distribution (density) $p(X\,|\,\theta)$ as the conditional distribution of the observation random variable $X$ given $\theta$. 
 After obtaining a collection of $n$ observations $X^{(n)}=\{X_1,X_2,\cdots,X_n\}$ modelled as $n$ independent copies of $X$ given $\theta$, we update our beliefs about $\theta$ from the prior by calculating the posterior distribution (density)
 \begin{align}\label{eqn:Bayes_post}
     p(\theta\,|\,X^{(n)}) = \frac{\exp\big\{\log\pi(\theta)+\log \m L_n(\theta;\,X^{(n)})\big\}}{\int_\Theta \exp\big\{\log\pi(\theta)+\log \m L_n(\theta;\,X^{(n)})\big\}\,\dd\theta}, \quad\theta\in\Theta,
 \end{align}
 where recall that $\m L_n(\theta;\,X^{(n)})=\prod_{i=1}^n p(X_i|\theta)$ is the likelihood function.
Despite the Bayesian formulation, in our theoretical analysis, we will adopt the frequentist perspective by assuming the data $X^{(n)}$ to be i.i.d.~samples from an unknown data generating distribution $\m P^\ast:\,=p(X\,|\,\theta^\ast)$, where $\theta^\ast$ will be referred to as the true parameter, or simply truth, throughout the rest of the paper.

In many real situations, practitioners may not be interested in learning the entire data generating distribution $\m P^\ast$, but want to draw inference on some characteristic as a functional $\theta=\theta(\m P^\ast)$ of $\m P^\ast$, which alone does not fully specify $\m P^\ast$.
An illustrative example is the quantile regression where the goal is to learn the conditional quantile of the response given the covariates; however, the conventional Bayesian framework requires a full specification of the condition distribution by imposing extra restrictive assumptions on the model, which may lead to model misspecification and sacrifice estimation robustness.
A natural idea to alleviate the limitation of requiring a well-specified likelihood function is to replace the log-likelihood function $\log \m L_n(\theta;\,X^{(n)})$ in the usual Bayesian posterior~\eqref{eqn:Bayes_post} by a criterion function $\m C_n(\theta;\, X^{(n)})$. The resulting distribution,
\begin{align}\label{eqn:Bayes_q_post}
     \pi_n(\theta\,|\,X^{(n)}) = \frac{\exp\big\{\log\pi(\theta)+\m C_n(\theta;\,X^{(n)})\big\}}{\int_\Theta \exp\big\{\log\pi(\theta)+ \m C_n(\theta;\,X^{(n)})\big\}\,\dd\theta}, \quad\theta\in\Theta,
\end{align}
is called the Bayesian pseudo-posterior with criterion function $\m C_n:\, \Theta\times \m X^n\to\mb R$, and we may use the shorthand $\pi_n(\cdot)$ to denote the  pseudo-posterior $\pi_n(\cdot|X^{(n)})$ when no ambiguity may arise. A popular choice of a criterion function is $\m C_n(\theta;\, X^{(n)})=-\alpha\, n\, \m R_n(\theta)$, where
$$\m R_n(\theta):\,=n^{-1} \, \sum_{i=1}^n \ell(X_i,\,\theta)$$ is the empirical risk function induced from a loss function $\ell:\,\m X\times\Theta \to\mb R$, and $\alpha\in(0,\infty)$ is the learning rate parameter. The corresponding Bayesian pseudo-posterior is called the Gibbs posterior associated with loss function $\ell$ in the literature~\citep[e.g.\!][]{bhattacharya2020gibbs,syring2020gibbs}.
In particular, the usual Bayesian posterior~\eqref{eqn:Bayes_post} is a special case when the loss function is $\ell(X,\theta)=-\log p(X\,|\,\theta)$ and $\alpha=1$. For Bayesian quantile regression, we may take the check loss function $\ell(x,q) = (q-x)\cdot \big(\tau-\mathbf{1}(q<x)\big)$ for a given quantile level $\tau\in (0,1)$, since the $\tau$-th quantile of any one-dimensional random variable $X$ corresponds to the population risk function minimizer $\argmin_{q\in\mb R}\mb E[\ell(X,q)]$.

A direct computation of either the posterior $p(\theta\,|\,X^{(n)})$ or the pseudo-posterior~\eqref{eqn:Bayes_q_post} involves the normalisation constant (the denominator) as a $d$-dimensional integral, which is often analytically intractable unless the prior distributions form a conjugate family to the likelihood (criterion) function.  In practice, Markov chain Monte Carlo (MCMC) algorithm~\citep{hastings1970monte,geman1984stochastic,robert2004monte} is instead employed as an automatic machinery for sampling from the (pseudo-)posterior, whose implementation is free of the unknown normalisation constant. The aim of this paper is to provide a rigorous theoretical analysis on the computational complexity of a popular and widely used class of MCMC algorithms described below. In particular, we are interested in characterizing a sharp dependence of their mixing times on the parameter dimension in the context of Bayesian posterior sampling.

\subsection{Metropolis-adjusted Langevin algorithm}\label{intro:MALA}
Consider a generic (possibly unnormalized) density function $f(\theta)=\exp\{-U(\theta)\}$ defined on a set $\Theta\subset\mb R^d$, where $U:\, \Theta\to\mb R$ is called the potential (function) associated with $f$. For example, in the Bayesian setting with target posterior~\eqref{eqn:Bayes_q_post}, we can take $U(\theta) = -\log\pi(\theta) - \m C_n(\theta;\,X^{(n)})$.
Suppose our goal is to sample from the probability distribution $\mu$ induced by $f$, where $\mu(A)=\frac{\int_A f(\theta)\, \dd\theta}{\int_\Theta f(\theta)\, \dd\theta}$ for any measurable set $A\subset \Theta$. 
Metropolis-adjusted Langevin algorithm (MALA), as an instance of MCMC with a special design of the proposal distribution, aims at 
producing a sequence of random points $\{\theta_k\}_{k\geq 0}$ in $\Theta$ such that the distribution of $\theta_k$ approaches $\mu$ as $k$ tends to infinity, so that for sufficiently large $k_0$, the $k_0$-th iterate $\theta_{k_0}$ can be viewed as a random variable approximately sampled from the target distribution $\mu$. In practice, every $k_0$ iterates from the chain can be collected (called thinning), which together form approximately independent draws from $\mu$.

Specifically, given step size $\wt h>0$ and initial distribution $\mu_0$ on $\Theta$, MALA produces $\{\theta_k\}_{k\geq 0}$ sequentially as follows: for $k=0,1,2,\ldots$,
\begin{enumerate} 
\item ({\bf Initialization}) If $k=0$, sample $\theta_0$ from $\mu_0$;
\item ({\bf Proposal}) If $k\geq 1$, given previous state $\theta_{k-1}$, generate a candidate point $y_k$ from proposal distribution $N_d\big(\theta_{k-1}-\wt h\nabla U(\theta_{k-1}), \,2\wt h\, I_d\big)$ whose density function is denoted as  $Q(\theta_{k-1},\,\cdot)$,
or equivalently, 
$$
y_k = \theta_{k-1}-\wt h\,\nabla U(\theta_{k-1}) +\sqrt{2\wt h} \,z_k,\quad\mbox{with } z_k\sim N_d(0,\, I_d).
$$
\item ({\bf Metropolis-Hasting rejection/correction}) Set acceptance probability $A(\theta_{k-1},\,y_{k}):\,= 1\wedge \alpha(\theta_{k-1},\,y_{k})$ with acceptance ratio statistic $$
\alpha(\theta_{k-1},y_{k}):\,=\frac{f(y_k)\,\cdot \,Q(y_k,\,\theta_{k-1})}{f(\theta_{k-1})\cdot Q(\theta_{k-1},\,y_k)}.
$$
Flip a coin and accept $y_k$ with probability $A(\theta_{k-1},\,y_{k})$ and set $\theta_k=y_k$; otherwise, set $\theta_k=\theta_{k-1}$.
\end{enumerate}
It is straightforward to verify that MALA described above produces a Markov chain whose transition kernel is
\begin{equation}\label{eqn:MALA_tran_1}
    T(\theta,\,\dd y)=\Big( \underbrace{1-\int_{\Theta} A(\theta,\,y)\, Q(\theta,\,y)\, \dd y}_{\mbox{\small rejection probability}} \Big) \cdot \delta_{\theta}(\dd y)+A(\theta,\,y)\,Q(\theta,\,y)\,\dd y,
\end{equation}
 where $\delta_\theta$ denotes the point mass measure at $\theta$.   In practice, the target density $f$ can be non-smooth at certain point $\theta\in \Theta$, and we address this issue by replacing the gradient $\nabla  U(\theta)$ with any of its subgradient $\widetilde\nabla  U(\theta)$\footnote{A subgradient of a function $f:\mb R^d\to \mb R$ at point $x\in\mb R^d$ is a vector $g\in\mb R^d$ such that
$f(y) \geq f(x) + \langle g,\, y-x \rangle + \m O(\|y-x\|)$ as $y\to x$} in MALA. That means, the proposal distribution $Q$ is being chosen as  $N_d(\theta_{k-1}-\wt h\,\wt\nabla U(\theta_{k-1}),2\wt h)$ and other aspects of the MALA algorithm remain unchanged. Furthermore, MALA can be generalized by introducing a symmetric positive-definite preconditioning matrix $\wt I\in \mb {R}^{d\times d}$, so that the proposal $Q$ in MALA  is modified as $N_d(\theta_{k-1}-\wt h\widetilde{I}\,\wt\nabla U(\theta_{k-1}),2\wt h\,\widetilde{I})$. It has been shown that~\citep{https://doi.org/10.1111/j.1467-9868.2010.00765.x,5947220} for a suitable preconditioning matrix, the resulting preconditioned MALA can help to alleviate the issue caused by the anisotropicity of the target measure. We illustrate both empirically  (c.f. Appendix~\ref{addex}) and theoretically (c.f.~Corollary~\ref{cor:smoothloss}) that a suitable preconditioning matrix may improve the convergence of the sampling algorithm for Bayesian posteriors.  {
As a common practice~\citep{JMLR:v21:19-441,https://doi.org/10.1002/rsa.3240040402} to simplify the analysis of MALA, in this paper, we consider the $\zeta$-lazy version of MALA, where at each iteration,  the chain is forced to remain unchanged with probability $\zeta$. The corresponding Markov transition kernel of the $\zeta$-lazy version of MALA is given by 
\begin{equation}\label{eqn:MALA_tran}
    \begin{aligned}
        T^{\zeta}(\theta,\dd y)=\big(1-(1-\zeta)\cdot\int_{\Theta}A(\theta,y)\,Q(\theta,y)\, \dd y\big)\cdot\delta_{\theta}(\dd y)+(1-\zeta)\cdot A(\theta,y)Q(\theta,y)\dd y.
    \end{aligned}
\end{equation}  }

 A closely related algorithm is the unadjusted Langevin algorithm~\citep[ULA,][]{durmus2017nonasymptotic,cheng2018sharp,roberts1996exponential,pmlr-v65-dalalyan17a}, which corresponds to discretization of the following Langevin stochastic differential equation (SDE), 
 $$
 \dd X_t= -\nabla U(X_t)\,\dd t+ \sqrt{2}\,\dd B_t, \quad t>0,
 $$
 and does not have
 the Metropolis-Hasting correction step 3. As a consequence, the stationary distribution of ULA is of order $\m O(\sqrt{dh})$ away from $\mu$ under several commonly used metrics~\citep{durmus2019analysis}. Due to this error, even in the strongly log-concave scenario, unlike MALA which requires 
 at most poly-$\log(1/\varepsilon)$ iterations with a constant step size $h$ to get one sample distributed close from $\mu$ with accuracy $\varepsilon$, ULA requires poly-$(1/\varepsilon)$ iterations and an $\varepsilon$-dependent choice of $h$~\citep{durmus2019analysis}.
 
 Another closely related algorithm is the classical Metropolis random walk (MRW), which instead uses $N_d\big(\theta_{k-1}, \,2\wt h\, I_d\big)$ without the gradient term in the proposal distribution $Q$. As we will see,  by using the extra gradient information,   the dimension dependence of the mixing time can be improved from $\m O(d)$~\citep{gelman1997weak,JMLR:v20:19-306} to $\m O(d^{1/3})$ for sampling from Bayesian posteriors.

 \subsection{Problem setup}\label{sec:problem}
 The goal of this paper is to characterize the mixing time of MALA for sampling from the Bayesian pseudo-posterior $\pi_n$ defined in~\eqref{eqn:Bayes_q_post}.  Assume we have access to a \emph{warm start} defined as follows.
 \begin{definition}
 We say $\mu_0$ is an $M_0$-warm start with respect to the stationary distribution $\mu$, if $\mu_0(E)\leq M_0\, \mu(E)$ holds for all Borel set $E\subset \mb R^d$, and we call $M_0$ the warming parameter.
 \end{definition}
 We state our problem as \emph{characterizing the $\varepsilon$-mixing time in $\chi^2$ divergence of the Markov chain produced by (preconditioned) MALA starting from an arbitrary $M_0$-warm start $\mu_0$ for obtaining draws from $\pi_{n}(\theta)$}, which is mathematically defined as the maximum of the minimal number of steps required for the chain to be within $\varepsilon^2$-$\chi^2$ divergence from its stationary distribution, over $M_0$-warm starts, or
 \begin{equation*}
 \begin{aligned}
      &\qquad\tau_{\rm mix}(\varepsilon,M_0)= {\max}\,\{   \tau_{\rm mix}(\varepsilon, \mu_0)\,:\,\mu_0 \text{ is an } M_0\text{-warm start with respect to }\pi_n \}\\
      &\qquad\qquad\qquad\qquad\text{with }\,\tau_{\rm mix}(\varepsilon, \mu_0)={\inf} \big\{k\in\mb N:\,\sqrt{\chi^2\big(\mu_k,\,\pi_{n}\big)}\leq \varepsilon\big\},
 \end{aligned}
 \end{equation*}
 where $\mu_k$ denotes the probability distribution obtained after $k$ steps of the Markov chain. Note that a mixing time upper bound in $\chi^2$ divergence implies that in total variation distance since $\|p-q\|_{\rm TV} \leq \sqrt{\chi^2(p,\,q)}$.

 \section{Mixing Time Bounds via  $s$-Conductance Profile}\label{mixingboundscp}
In this section, we introduce a general technique of using $s$-conductance profile to bound the mixing time of a Markov chain.  We first review some common concepts and previous results in Markov chain convergence analysis, and then provide an improved analysis for obtaining a sharp mixing time upper bound of MALA in this work.

\smallskip
\noindent {\bf Ergodic Markov chains:} {
Given a Markov transition kernel $T(\cdot,\,\cdot)$  with stationary distribution $\mu\in \m P(\mb R^d)$, the ergodic flow of a set $S$ is defined as 
$$
\phi(S)=\int_{S} \bigg\{\int_{S^c} T(\xi,\,\dd y)\bigg\}\,\mu(\dd\xi).
$$
The ergodic flow captures the mass of points leaving $S$ (i.e., $T(\xi, S^c)=\int_{S^c} T(\xi,\,\dd y)$) on average under stationary distribution $\mu$ in one step of the Markov chain.} A Markov chain is said to be ergodic if $\phi(S)>0$ for all measurable set $S\subset \mb R^d$ with $0<\mu(S)<1$. Let $\mu_k$ denote the probability distribution obtained after $k$ steps of a Markov chain. If the Markov chain is ergodic, then $\mu_k\to \mu$ as $k\to\infty$ in total variation distance; see, for example, Corollary 1.6 of~\cite{https://doi.org/10.1002/rsa.3240040402}.

\smallskip
\noindent {\bf Conductance of Markov chain and rapid mixing:}  The (global) conductance of an ergodic Markov chain characterizes
the least relative ratio between $\phi(S)$ and the measure $\mu(S)$ of $S$, and is formally defined as 
$$
\Phi=\inf\bigg\{\frac{\phi(S)}{\mu(S)}:\,  0<\mu(S)\leq \frac{1}{2}\bigg\}.
$$
{
A Markov chain with low conductance tends to become trapped in a subset of its states, whereas one with high conductance has more freedom to explore and transition across its entire state space.}
The conductance is related to the spectral gap\footnote{The spectral gap is define as $\Lambda=\inf\{ \m E(f,f)/{\rm Var}_{\mu}(f)\,:\, f\in L^2(\mu), {\rm Var}_{\mu}(f)>0\}$, where $\m E(f,g)=\int (f(x)-g(y))^2 T(x,\dd y)\,\dd \mu(x)$ is the \emph{Dirichlet form}.} of the Markov chain via Cheeger’s inequality~\citep{cheeger2015lower}, and thus can be used to characterize the convergence of the Markov chain. For example, Corollary 1.5 in~\cite{https://doi.org/10.1002/rsa.3240040402} shows that if $\mu_0$ is an $M_0$-warm start with respect to the stationary distribution $\mu$, then 
$$
\| \mu_k - \mu\|_{  \rm TV} \leq \sqrt{M_0}\, \Big(1-\frac{\Phi^2}{2}\Big)^k,\quad k\geq 0.
$$
{
Furthermore, some people consider the more flexible notion of $s$-conductance, defined as 
$$
\Phi_s:=\inf\bigg\{ \frac{\phi(S)}{\mu(S)-s}:\, s<\mu(S)\leq \frac{1}{2}\bigg\},\quad\mbox{for }s\in(0,1/2),
$$
which restricts the infimum over all sets with a probability greater than $s$. This restriction avoids including sets in the conductance bound that have poor conductance but receive negligible probability, which should be less significant to the overall mixing of the Markov chain. Specifically for sampling from Bayesian posteriors, this refined analysis allows us to focus our calculations on these ``highest posterior regions" while avoiding some unwieldy tail probability regions (e.g., the region defined in Condition A.3).}
Using the $s$-conductance, Corollary 1.6 in~\cite{https://doi.org/10.1002/rsa.3240040402} proves a similar bound implying the exponential convergence of the algorithm up to accuracy level $s$ as 
$$
\| \mu_k - \mu\|_{  \rm TV} \leq M_0\,s+M_0\,\Big(1-\frac{\Phi_s^2}{2}\Big)^k,\quad k\geq 0.
$$
Consequently, the $\varepsilon$-mixing time with respect to the total variation distance of the Markov chain starting from  an $M_0$-warm start can be upper bounded by $\frac{2}{\Phi_s^2}\log \frac{2M_0}{\varepsilon}$ if we choose $s=\frac{\varepsilon}{2M_0}$.

  \smallskip
 \noindent {\bf Conductance profile of Markov chain:}  Instead of controlling mixing times  via  a worst-case conductance bound,  some recent works have introduced more refined methods based on the conductance profile.  The conductance profile is defined as the following collection of conductance,
 $$
 \Phi(v):=\inf\bigg\{\frac{\phi(S)}{\mu(S)}\,:\,0<\mu(S)\leq v\bigg\}, \quad \mbox{indexed by }v\in\Big(0,\,\frac{1}{2}\,\Big].
 $$  
 Note that the classic conductance constant $\Phi$ is a special case that can be expressed as $\Phi=\Phi(\frac{1}{2})$. Based on the conductance profile,~\cite{JMLR:v21:19-441} consider the concept of $\Omega$-restricted conductance profile for a convex set $\Omega$, given by 
 \begin{equation*}
    \Phi_{\Omega}(v):=\inf\bigg\{\frac{\phi(S)}{\mu(S\cap \Omega)}\,:\,0<\mu(S\cap \Omega)\leq v\bigg\},\quad v\in\Big(0,\,\frac{\mu(\Omega)}{2}\,\Big].
 \end{equation*}
  It has been shown in~\cite{JMLR:v21:19-441} that given an $M_0$-warm start $\mu_0$,  if $$
  \mu(\Omega)\geq 1-\frac{\varepsilon^2}{3M_0^2} \quad\mbox{and}\quad \Phi_{\Omega}(v)\geq \sqrt{B\log \frac{1}{v}} \ \ \mbox{for all } \ v\in \Big[\,\frac{4}{M_0},\,\frac{1}{2}\,\Big],
  $$
  then the $\varepsilon$-mixing time in $\chi^2$ divergence of the chain is bounded from above by $\m O\big(\frac{1}{B}\log (\frac{\log M_0}{\varepsilon})\big)$.  Therefore, compared with the (global) conductance, employing the technique of conductance profile may improve the warming parameter dependence in the mixing time bound from $\log M_0$ to $\log\log M_0$. This improvement from a logarithmic dependence to the double logarithmic dependence may dramatically sharpen the mixing time upper bound, since in a typical Bayesian setting $M_0$ may grow exponentially in the dimension $d$.  However, one drawback of the conductance profile technique from~\cite{JMLR:v21:19-441} is that the high probability set $\Omega$ should be constrained to be convex (Lemma 4 of~\cite{JMLR:v21:19-441}) to bound the $\Omega$-restricted  conductance profile $\Phi_{\Omega}(v)$. This convexity constraint may cause $\Phi_{\Omega}(v)$ to have a worse dimension dependence compared with the complexity analysis using the $s$-conductance $\Phi_s$.

 In order to address the above issues of previous analysis,  we introduce the following notion of \emph{$s$-conductance profile }, which combines ideas from the $s$-conductance and conductance profile,
\begin{equation*}
    \Phi_s(v):=\inf \left\{\frac{\phi(S)}{\mu(S)-s}\,\bigg|\, s< \mu(S)\leq v\right\} \quad \mbox{indexed by } \  s\in \Big(0,\,\frac{1}{2}\,\Big) \ \ \mbox{and} \ \ v\in \Big(s,\,\frac{1}{2}\,\Big].
\end{equation*}
The $s$-conductance profile evaluated at $v=\frac{1}{2}$ corresponds to the $s$-conductance that is commonly-used in previous study for analyzing the mixing time of Markov chain~\citep{pmlr-v134-chewi21a,JMLR:v20:19-306}. We show in the following lemmas that a lower bound on the $s$-conductance profile can be translated into an upper bound on the mixing time in $\chi^2$-squared divergence. {
We formulate here an informal result and postpone a more detailed statement to Appendix~\ref{app:profile}.

 \begin{lemma}[\bf Mixing time bound via $s$-conductance profile (informal)]\label{lemma:conductance_informal} For any error tolerance $\varepsilon\in(0,1)$, the  mixing time in $\chi^2$ divergence of the $\zeta$-lazy version of MALA over $M_0$-warm starts can be bounded as
\begin{equation*}
    \tau_{\rm mix}(\varepsilon,M_0)\lesssim \zeta^{-1}\cdot\Big(  \int_{\frac{4}{M_0}}^{\frac{1}{2}} \frac{\dd v}{v\,\Phi_s^2(v)}+ \frac{1}{\Phi^2_s(\frac{1}{2})}\log(\frac{1}{\varepsilon})\Big),\quad  s=\frac{\varepsilon^2}{16M_0^2}.
\end{equation*}

\end{lemma}

It is worth noting that since $\Phi_s(v)$ is a decreasing function of $v$, by replacing $\Phi_s(v)$ with its lower bound $\Phi_s=\Phi_s(\frac{1}{2})$, one can obtain a mixing time bound via $s$-conductance. However, instead of simply considering the worst case, the integral $\int_{\frac{4}{M_0}}^{\frac{1}{2}} \frac{\dd v}{v\,\Phi_s^2(v)}$ averages over $\Phi_s(v)$, offering a possible improvement in the dependence on warming parameter $M_0$.  To establish a lower bound for the $s$-conductance profile, we can employ the ``overlap argument" frequently used in the literature~\citep{pmlr-v134-chewi21a,JMLR:v21:19-441,10.2307/30243694,wu2022minimax}, that is, 1.~prove a log-isoperimetric inequality for $\mu$; 2.~bound the total variation distance between $T(x,\cdot)$ and $T(z,\cdot)$ for any two sufficiently close points $x, z$ in a high probability set (not necessarily convex) of $\mu$. We leave a detailed description of this argument to Appendix~\ref{app:profile}.}

Among previous works of mixing time analysis of MALA,~\cite{JMLR:v21:19-441} study the problem of sampling from general smooth and strongly log-concave densities, using the technique of $\Omega$-restricted conductance profile. Their bound has a double logarithmic $\log\log M_0$ dependence on the warmth parameter $M_0$ under certain regime (of step size $h$), and a sub-optimal $\m O(d)$-dependence on the dimension.
On the other hand,~\cite{pmlr-v134-chewi21a} study the same problem as~\cite{JMLR:v21:19-441} and obtain a mixing time bound with an optimal $\m O(d^{\frac{1}{2}})$-dependence, based on the $s$-conductance technique. However, the bound in~\cite{pmlr-v134-chewi21a} has a quadratic dependence on $\log M_0$. By utilizing our $s$-conductance profile argument, when $\log M_0$ and $h^{-1}$ are not of constant order, we can improve their bounds from $h^{-1}\log(\frac{M_0}{\epsilon})$ to $\max\{h^{-1}\log(\frac{\log M_0}{\epsilon}),\, \log M_0\}$, where $h$ is the step size used in Theorem 3 of~\cite{pmlr-v134-chewi21a}.

 \section{Mixing Time of MALA}\label{sec:MALA_mixtime}
In this section, we describe our main result by providing an upper bound to the mixing time of (preconditioned) MALA for sampling from the Bayesian pseudo-posterior $\pi_n$. We consider the $\zeta$-lazy version of MALA and assume that a warm start is accessible, which is a common assumption~\citep[e.g.][]{JMLR:v20:19-306,mangoubi2019nonconvex}. For example, Corollary~\ref{cor:smoothloss} in Section~\ref{sec:smoothloss} provides a construction of $M_0$-warm start for general Gibbs posterior with smooth criterion function, where $M_0$ is bounded above by an $(n,d)$-independent constant.

Note that the Bayesian pseudo-posterior with criterion function $\m C_n$ can be rewritten as 
\begin{align}
    &\qquad\qquad \pi_{n}(\theta\,|\,X^{(n)}) = \frac{\exp\big\{-V_n\big(\sqrt{n}(\theta-\wh\theta)\big)\big\}}{\int_{\Theta} \exp\big\{-V_n\big(\sqrt{n}(\theta-\wh\theta)\big)\big\}\,\dd\theta} \quad \forall \theta\in\Theta, \label{Eqn:loc_posterior}\\
    &\qquad\mbox{where}\quad \hat\theta=\underset{\theta\in \Theta}{\argmax}\,\,\m C_n(\theta) \quad\mbox{and}\\ 
    & V_n(\xi)=-\m C_n\Big(\widehat{\theta}+\frac{\xi}{\sqrt{n}};X^{(n)}\Big)+\m C_n\big(\widehat{\theta}\,; X^{(n)}\big)-\log \pi\Big(\widehat{\theta}+\frac{\xi}{\sqrt{n}}\Big)+\log \pi(\widehat\theta\,) \notag
\end{align}
is the corresponding rescaled potential (function). In the expression of $V_n$, we deliberately added two terms independent of $\xi$ so that $V_n(0)=0$ for simplifying the analysis. Motivated by the classical Bernstein-von Mises (BvM) theorem\footnote{When sample size $n$ is large, the Bayesian posterior is close to the Gaussian distribution $N_d(\widehat{\theta}_{\rm MLE},\,n^{-1}\m J^{-1})$, where $\widehat{\theta}_{\rm MLE}$ is the maximum likelihood estimator and $\m J$ the Fisher information matrix.}~\citep{van2000asymptotic, Ghosh2003} for Bayesian posteriors, we impose following conditions on $V_n$, stating that $V_n(\xi)$ is close to a quadratic form and the subgradient of $V_n(\xi)$ employed in MALA is close to a linear form, uniformly over a high probability set of the rescaled target measure $\pi_{\rm loc}=(\sqrt{n}(\cdot-\wh\theta))_{\#}\pi_n$.\footnote{We use $\mu=G_{\#}\nu$ to denote the push forward measure so that for any measurable set $A$, $\mu(A)=\nu(G^{-1}(A))$.} Here $\pi_{\rm loc}$ corresponds to the  measure of the localized random variable $\xi=\sqrt{n}(\theta-\wh\theta)$ for $\theta\sim \pi_n(\theta|X^{(n)})$, and the transformation $\sqrt{n}(\cdot-\wh\theta)$ makes the limiting distribution of $\xi$ zero-centered and has constant-order variances.

\vspace{0.5em}
\noindent\textbf{Condition A:}  {\it  There exists 
a tolerance $\varepsilon\in (0,1)$, preconditioning matrix $\wt I$, step size parameter $h$ (rescaled by $n$), warming parameter $M_0$, numbers $R,\wt\varepsilon_0,\wt\varepsilon_1\geq 0$, $\rho_1,\rho_2>0$ and  a symmetric positive definite matrix $J\in \mb R^{d\times d}$  so that 
\begin{enumerate}
    \item for any $\xi\in K=\{x:\|\wt I^{-1/2}x\|\leq R\}$\footnote{Here the notation $A^{-1/2}$ of a symmetric positive definite matrix $A$ means the inverse of its matrix square root $A^{1/2}$.}
 \begin{equation*}
 \begin{aligned}
  \big|V_n(\xi)-\frac{1}{2}{\xi^TJ\xi} \big|\leq  \wt\varepsilon_0 \quad\mbox{and}\quad 
   \big\|\widetilde{\nabla}V_n(\xi)-J\xi\|\leq\widetilde\varepsilon_1,
       \end{aligned}
 \end{equation*}
where $\widetilde{\nabla}V_n(\xi)$ is a subgradient of $V_n(\xi)$;
\item $\rho_1 I_d \preceq \wt J=\widetilde{I}^{1/2} J \widetilde{I}^{1/2} \preceq \rho_2 I_d$;
\item $\pi_n\big(\sqrt{n}\,\|\wt {I}^{-1/2}(\theta-\wh\theta)\|\leq R/2\big)\geq 1-\exp(-4\wt\varepsilon_0)\cdot\frac{h\rho_1\varepsilon^2}{M_0^2}$ and $R\geq 8\sqrt{d/\lambda_{\min}(\wt J)}\,$.
\end{enumerate}
}

\smallskip

{
The first inequality in Condition A.1 requires that $V_n(\xi)$ can be uniformly approximated by the quadratic term $\frac{1}{2}\xi^TJ\xi$ with an approximation error $\wt\varepsilon_0$. This requirement is implied by the classical BvM result, which is commonly utilized in MCMC mixing time analysis for Bayesian posterior sampling~\citep{10.2307/30243694,ascolani2023dimensionfree}. It is noteworthy that we do not impose any smoothness or convexity constraints on $V_n(\xi)$, and the deviation characteristic $\wt\varepsilon_0$ can take any value. We also keep track of the impact of this deviation in the final mixing time bound, as reflected in Theorem 3, where we explicitly show the dependency of the mixing time on this approximation error, $\wt\varepsilon_0$. The result reveals that the mixing time exhibits an exponential dependence on $\wt\varepsilon_0$.
The second inequality in Condition A.1 assumes that the subgradient of $V_n(\xi)$ can be approximated by the linear term $J\xi$ with an approximation error $\wt\varepsilon_1$. Although less standard, this condition is crucial since $\wt\varepsilon_1$ governs the efficacy of the subgradient used in MALA to adjust the proposal distribution and facilitate faster exploration of the parameter space. As we will see in Theorem 3, a small $\wt\varepsilon_1$ enables MALA, leveraging (sub)gradient information, to improve upon MRW in terms of mixing time.
Condition A.2 requires the asymptotic covariance matrix $J$, after rescaling by the preconditioning matrix, to have its maximum eigenvalue upper-bounded by $\rho_2$ and its minimum eigenvalue lower-bounded by $\rho_1$. The condition number $\kappa=\frac{\rho_2}{\rho_1}$ serves as an indicator of how well the preconditioning matrix $\wt I$ is chosen to alleviate issues arising from the anisotropy of the target distribution. As we will see from Theorem 3, a small $\kappa$ will lead to a lower mixing time.
The last condition (Condition A.3) assumes that the radius $R$ of the compact set $K$, considered in Condition A.1, is sufficiently large. This ensures that $K$ is a high probability set under $\pi_{\rm loc}$. This assumption guarantees that the region where the density $\pi_{\rm loc}$ (or $\pi_n$) deviates significantly from a Gaussian form, and is possibly non-smooth and non-log-concave, is negligible, thereby reducing the chances of the Markov chain becoming trapped in such regions.

In summary, Condition A requires the localized (rescaled) posterior $\pi_{\rm loc}=(\sqrt{n}(\cdot-\wh\theta))_{\#}\pi_n$ to be close to a Gaussian distribution $N_d(0,J^{-1})$, so that we can analyze the mixing time of MALA for sampling $\pi_n$ or $\pi_{\rm loc}$ (note that the complexity for sampling from $\pi_n$ with step size $\wt h=h/n$ is equivalent to that from $\pi_{\rm loc}$ with rescaled step size $h$) by comparing its transition kernel $T$ expressed in~\eqref{eqn:MALA_tran} with the transition kernel $T^{\Delta}$ induced from the MALA for sampling the Gaussian distribution. Interestingly, we find that as long as  the deviance of $\pi_{\rm loc}$ to Gaussian is sufficiently small but not necessarily diminishing as $n,d\to \infty$, some key properties (more precisely, conductance lower bound) of $T^\Delta$ guarantee that the fast mixing of MALA will be inherited by $T$, so that the mixing time associated with $T$ can be controlled.  Using this argument, we prove a mixing time upper bound without imposing the smoothness and strongly convexity assumptions on $V_n(\xi)$ that are restrictive and commonly assumed in the literature for analyzing the convergence of MALA~\citep{pmlr-v134-chewi21a,JMLR:v21:19-441}.  As a concrete example, under mild assumptions, Condition A holds for a broad class of Gibbs posteriors~\citep{bhattacharya2020gibbs} mentioned in Section~\ref{pseudo-posterior} where the criterion function $\m C_n$ is proportional to the negative empirical risk function $\m R_n$, as long as $d$ is relatively small compared to $n$ (see Lemma~\ref{th1} and Lemma~\ref{lemmatail} in Appendix~\ref{sec:proofmixing} for details). Now we are ready to state the following theorem.}

\begin{theorem}[\bf MALA mixing time upper bound]\label{thmala}
Let  $\pi_{n}$ defined in~\eqref{Eqn:loc_posterior}  be the target distribution and $\zeta\in (0,\frac{1}{2}]$ be a lazy parameter.   Assume Condition A holds for a tolerance $\varepsilon$, warming parameter $M_0$, sample size $n$, preconditioning matrix $\wt I$,  rescaled step size $h$, and some $R>0$, $\wt\varepsilon_1\geq 0$, $\rho_2\geq\rho_1>0$, and that there exists a small enough absolute $(n,d)$-independent constant $c_0$  so that the step size can be expressed as $\wt h=h/n$ with  
\begin{equation*}
    h=c_0\cdot\bigg[\rho_2\Big(d^{\frac{1}{3}}+d^{\frac{1}{4}}\Big(\wt \varepsilon_0+\log \frac{M_0d\kappa}{\varepsilon}\Big)^{\frac{1}{4}}+\Big(\wt \varepsilon_0+\log \frac{M_0d\kappa}{\varepsilon}\Big)^{\frac{1}{2}}+\mnorm{\wt I}_{\rm op}R^2\wt\varepsilon_1^2\Big)\bigg]^{-1}, \text{ where } \kappa=\frac{\rho_2}{\rho_1},
     \end{equation*}
  then the $\zeta$-lazy version of MALA with  proposal distribution $N_d(\theta_{k-1}-\wt h\widetilde{I}\,\wt\nabla U(\theta_{k-1}),2\wt h\,\widetilde{I})$ and step size $\wt h$ has a maximal $\varepsilon$-mixing time in $\chi^2$ divergence over $M_0$-warm starts being bounded as
 \begin{equation}\label{mixingbound}
     \tau_{\rm mix}(\varepsilon,M_0)\leq \frac{C_1\exp(4\wt\varepsilon_0)}{\zeta} \cdot\bigg\{\bigg[\rho_1^{-1}\exp(8\wt\varepsilon_0)\cdot h^{-1}\log\big(\frac{\log M_0}{\varepsilon}\big)\bigg]\vee \log M_0\bigg\},
 \end{equation}
 where $C_1$ is an $(n,d)$-independent constant.
 \end{theorem}

  The mixing time bound~\eqref{mixingbound} is proved using the technique of $s$-conductance profile introduced in Section~\ref{mixingboundscp}. A similar mixing time bound can be obtained if when consider the sampling of $\pi_{\rm loc}$ constrained on the high probability set $K=\{x:\|\wt I^{-1/2}x\|\leq R\}$, which is adopted by~\cite{10.2307/30243694} for analyzing the mixing time of MRW; however, our result does not require such a constraining step. According to Theorem~\ref{thmala}, for a fixed tolerance (accuracy level) $\varepsilon$, the $\varepsilon$-mixing time is determined by the parameter dimension $d$, warming parameter $M_0$, preconditioning matrix $\wt I$, approximation errors $\wt\varepsilon_0$, $\wt\varepsilon_1$ of the potential and the gradient, radius $R$ of the high probability set of $\pi_{\rm loc}$ and the precision matrix $J$ of the Gaussian approximation to $\pi_{\rm loc}$. {
  The derived mixing time bound is exponentially dependent on $\wt\varepsilon_0$, implying that a bound that is polynomial in $d$ can only be attained if $\wt\varepsilon_0$ is either constant-order or logarithmic in $d$.} The fourth term $\mnorm{\wt I}_{\rm op}R^2\wt\varepsilon_1^2$ in the expression of $h$ will be dominated by others once $\wt\varepsilon_1$ is sufficiently small. For example, suppose $\wt I =I_d$, $\log \frac{M_0\kappa}{\varepsilon}=\m O(d)$ and $\pi_{\rm loc}$ has a sub-Gaussian type tail behavior, or 
  $$
  \pi_{\rm loc}\big(\|\xi\|\geq c_1(\sqrt{d}+t)\big)\leq \exp(-c_2 \,t^2), \quad t>0,
  $$
  then we can choose the radius as $R=\m O(\sqrt{d})$, and the term $\mnorm{\wt I}_{\rm op}R^2\wt\varepsilon_1^2$ will be dominated by the $\m O(d^{\frac{1}{3}})$ term once $\wt\varepsilon_1=\m O(d^{-\frac{1}{3}})$. This suggests that a $d^{\frac{1}{3}}$-mixing time upper bound is achievable as long as the (sub)gradient used in MALA deviates from a linear form with approximation error at most $d^{-\frac{1}{3}}$, which is independent of the sample size.  Therefore,  when $d\ll n$, it is safe to fix a mini-batch dataset for computing the (sub)gradient in MALA instead of using the full batch. As another remark, our theorem also gives a tight mixing time upper bound $\m O(d)$ of MRW by taking $\wt{\varepsilon}_1=O(1)$, corresponding to the case where the gradient estimate is completely uninformative.

Our mixing time bound has a linear dependence (modulo logarithmic term) on the condition number $\kappa={\rho_2}/{\rho_1}$, which matches the best condition number dependence for MALA under strong convexity~\citep{wu2022minimax} and we show the tightness of the condition number dependence in Theorem~\ref{thmala:lower} of Appendix~\ref{app:lower}.
Moreover, by introducing preconditioning matrix $\wt I$, a small condition number can be obtained once $\wt I$ acts as a reasonable estimator to $J^{-1}$, which will lead to a faster mixing time when $J$ is ill-conditioned. On the other hand, assume $\kappa$ is bounded above by an $(n,d)$-independent constant and 
$$
\big(\, \mnorm{\wt I}_{\rm op}\, R^2\, \wt\varepsilon_1^2\, \big)\vee \log \Big(\frac{M_0}{\varepsilon}\Big)\leq d^{\frac{1}{3}}, 
$$
we have $\tau_{\rm mix}(\varepsilon,\mu_0)\leq C_1\, d^{\frac{1}{3}}\log (\frac{\log M_0}{\varepsilon})$. This upper bound matches the lower bound proved in~\cite{pmlr-v134-chewi21a} that the mixing time of MALA for sampling from the standard Gaussian target is at least $\m O(d^{\frac{1}{3}})$, and it improves the warming parameter dependence from $\log M_0$ to $\log(\log M_0)$ compared with the upper bound proved in~\cite{pmlr-v134-chewi21a}. Therefore, in order to attain the best achievable mixing time $\m O(d^{\frac{1}{3}})$, we need to  find a initial distribution $\mu_0$  that is close to $\pi_{n}$, so that the warming parameter $M_0$ can be controlled.
{
For a generic log-concave distribution, it has been shown that a warm start with warming parameter $M_0$ polynomial in $d$ can be obtained with $d^{\frac{1}{2}}$ complexity, as demonstrated by~\cite{altschuler2023faster}. However, efficiently obtaining a poly($d$) warm start for general non-log-concave sampling problems is infeasible. Fortunately, in our Bayesian posterior sampling context, although $\pi_{n}$ may not be log-concave, large sample asymptotic theory (refer to Section~\ref{sec:application}, for instance) ensures that $\pi_{n}$ is approximately Gaussian. Therefore, using the Gaussian distribution $N_d(\wh \theta,n^{-1}\wt I)$, constrained on a compact set, as the initialization $\mu_0$, is a natural choice. To support this initialization scheme, the following lemma provides an upper bound for the corresponding warming parameter $M_0$.}
 \begin{lemma}[\bf Warming parameter control]\label{lemma:Warm}
Suppose Condition A is satisfied. For any compact set $K\subset \mb R^d$, the initial distribution as
\begin{equation*}
    \mu_0=N_d(\wh \theta,n^{-1}\widetilde{I})|_{\{\theta\,:\,\sqrt{n}(\theta-\wh\theta)\in K\}}
\end{equation*}
is $M_0$-warm with respect to $\pi_{n}$, where
\begin{equation*}
    \log M_0\leq -\log \pi_n\big(\{\theta\,:\,\sqrt{n}(\theta-\wh\theta)\in K\}\big)\,+\,\underset{\xi\in K}{\sup}{\big|\,\xi^T(\wt I^{-1}-J)\xi\big|}\,+\,2\cdot\underset{\xi \in K}{\sup}|V_n(\xi)-\frac{1}{2}x^TJx|.
\end{equation*}
 \end{lemma}

\noindent 
{
In order to construct a feasible warm start using Lemma~\ref{lemma:Warm}, it is necessary to compute the maximizer $\wh\theta$ of the criterion function $\m C_n(\theta)$. An inaccurate approximation of $\wh\theta$ may cause the warming parameter $M_0$ to grow linearly with the sample size $n$, a similar observation also noted in studies by~\cite{ascolani2023dimensionfree,10.2307/30243694}. While it is generally challenging to obtain solutions for non-convex optimization problems, there are cases where optimizing a nearly quadratic function can be much easier compared to sampling from a nearly Gaussian distribution. A specific example is Bayesian quantile regression, where the estimation of $\hat{\theta}$ can be efficiently achieved using linear programming techniques. Our theoretical results also suggest that under Condition A, we can control the warming parameter $M_0$ in MALA by choosing a reasonable estimator $\wt I$ for the inverse asymptotic covariance matrix $J^{-1}$ of $\pi_{\rm loc}$. For instance, if $\wt I$ is chosen to be the identity matrix and $J$ has a bounded operator norm, then $\log M_0$ should be of order $\m O(d)$. Furthermore, in Bayesian Gibbs posterior sampling, where the loss function $\ell$ is continuously twice differentiable, a viable option for approximating $J^{-1}$ could be the plug-in estimator:
$$
{\widetilde I}=\bigg\{\frac{1}{|S|}\sum_{i\in S} {\rm Hess}_{\theta}(\ell(X_i,
 \wh\theta))\bigg\}^{-1},
 $$ 
 where $S$ is a subset of ${1,2,\cdots,n}$, and ${\rm Hess}_{\theta}(\ell(x,\theta))$ denotes the Hessian matrix of $\ell(x,\cdot)$ evaluated at $\theta$. Notably, since the warming parameter $M_0$ can be of order $\m O(d^{1/3})$ for achieving the best possible mixing time, it is feasible to compute the plug-in estimator using only a mini-batch of data, the size of which depends solely on the dimension, rather than the full dataset. Further details can be found in Corollary~\ref{cor:smoothloss}.}
 
\smallskip

 \noindent {
 According to Lemma~\ref{lemma:Warm} and Theorem~\ref{thmala},  a reasonably good approximation $\wt I$ to matrix $J$ in Condition A will improve both the mixing time of MALA after burn-in period and the initialization affecting the burn-in.  For completeness, we also provide an experiment in Appendix~\ref{addex} for investigating the impact of the preconditioning matrix and initial distribution on the performance of MALA. }However, in some complicated problems, especially when $\log \pi_{\rm loc}$ is not differentiable, a good estimator for the matrix $J$ may not be easy to construct. One possible strategy is to use adaptive MALA~\citep{Atchade_2006}, where the preconditioner $\widetilde I$ and step size $h$ are updated in each iteration by using the history draws. It has been empirically shown in~\cite{Atchade_2006} that adaptive MALA outperforms non-adaptive counterparts in many interesting applications.  
 We leave a rigorous theoretical analysis of adaptive MALA as a future direction.  


 \section{Sampling from Gibbs Posteriors}\label{sec:application}
 
 Recall from Section~\ref{pseudo-posterior} that a Gibbs posterior is a Bayesian pseudo-posterior defined in~\eqref{eqn:Bayes_q_post} with the criterion function $\m C_n(\theta;\, X^{(n)})=-\alpha\, n\, \m R_n(\theta)$, where $\alpha$ is an $(n,d)$-independent positive learning rate and $\m R_n(\theta):\,=n^{-1} \, \sum_{i=1}^n \ell(X_i,\,\theta)$ is the empirical risk function induced from a loss function $\ell:\,\m X\times\Theta \to\mb R$. In this section, we first provide generic conditions under which Condition A for Theorem~\ref{thmala} can be verified for the the Gibbs posterior so that the mixing time bound of the corresponding MALA can be applied. After that, we specialize the result to two representative cases: Gibbs posterior with a generic smooth loss function, and Gibbs posterior in Bayesian quantile regression where the check loss function is non-smooth.
 
 Firstly, we make the following conditions on the population level risk function $\m R(\theta)=\mb E[\ell(X,\theta)]$. Recall that $\theta^\ast={\arg\min}_{\theta\in \Theta} \m R(\theta)$ denotes the true parameter. The key idea is that although the sample level risk function (i.e.~empirical risk function) $\m R_n$ is allowed to be non-smooth, but as the sample size $n$ grows, it becomes closer and closer to the population level risk function $\m R(\theta)$, which can be properly analyzed if smooth.

  \vspace{0.5em}
  
  \noindent \textbf{Condition B.1 (Risk function):} {\it For $(n,d)$-independent constants $(C'
 ,C,r)>0$ and $( \gamma_0,\gamma_1,\gamma_2 )\geq 0$:
   \begin{enumerate}
     \item  $\m R(\theta)$ is twice differentiable with mixed partial derivatives  of order two being uniformly bounded by $C$ on $B_r(\theta^\ast)$; for any $\theta\in \Theta$, $\m R(\theta)-\m R(\theta^\ast)\geq C'\,d^{-\gamma_0}\,(d^{-\gamma_1}\wedge\|\theta-\theta^\ast\|^2)$.
     \item Let $\m H_{\theta}$ denote the Hessian of $\m R$ at $\theta$. For any $\theta\in B_r(\theta^\ast)$,   $\mnorm{\m H_{\theta}-\m H_{\theta^\ast}}_{  \rm  op}\leq C\,d^{\gamma_2}\|\theta-\theta^\ast\|$.
  \end{enumerate}
}
{
\noindent Condition B.1 imposes two requirements. Firstly, the population level risk function $\mathcal{R}(\cdot)$ must possess a unique global minimizer $\theta^*$. This condition ensures that when the empirical risk $\mathcal{R}_n$ in the Gibbs posterior is substituted with $\mathcal{R}$, the resulting distribution  $\pi^*(\theta)\propto\exp(-\alpha \,n\, \m R(\theta))\pi(\theta)$ will be unimodal, thereby preventing the Markov chain from getting stuck in any local mode. Note that this condition is equivalent to the identifiability of the parameter in the model, and therefore is natural to assume. Secondly, the risk function should exhibit sufficient smoothness and local strong convexity in the vicinity of $\theta^*$. This property enables a reliable Gaussian approximation for the local shape of $\pi^*(\theta)$ around $\theta^*$, which is again a standard assumption and holds when the Fisher information matrix is not singular. Next, we introduce the following assumption of Lipschitz continuity for the loss function $\ell$.}
   \vspace{0.5em}
   
\noindent \textbf{Condition B.2 (Loss function):} {\it There exist $(n,d)$-independent constants $C>0$ and  $\gamma\geq 0$ such that for any $x\in \m X$ and  $(\theta,\,\theta')\in \Theta^2$, it holds that $|\ell (x,\theta)-\ell(x,\theta')|\leq C\, d^{\gamma}\,\|\theta-\theta'\|$. }
  
\vspace{0.5em}
\noindent If the loss function has uniformly bounded derivatives with respect to $\theta$, that is, $\big|\frac{\partial \ell (X,\theta)}{\partial\theta_j}\big|\leq C$ holds for any $j\in [d]$, $x\in \mathcal{X}$, and $\theta\in \Theta$, where $C$ is a constant independent of $n$ and $d$, then Condition B.2 holds with $\gamma=\frac{1}{2}$. Next, we introduce a function $g: \mathcal{X}\times \Theta \to \mathbb{R}^d$ that satisfies the following conditions.

 \vspace{0.5em}
  
  \noindent\textbf{Condition B.3 (Subgradient of loss function):} {\it There exist some $(n,d)$-independent  constants $(C,r,\beta_1)>0$ and  $(\gamma_3,\gamma_4)\geq 0$ so that:
  \begin{enumerate} 
      \item For any $\theta\in B_r(\theta^\ast)$, it holds $\mb E [g(X,\theta)]=\nabla \m R(\theta)$ and $\sup_{x\in \m X} \|g(x,\theta)\|\leq C\,d^{\gamma}$, where $\gamma$ is the same as that defined in Condition B.2.
      
      \item Let $d_n^{g}(\theta,\theta')=\sqrt{n^{-1}\sum_{i=1}^n \|g(X_i,\theta)-g(X_i,\theta')\|^2}\,$ be a pseudo-metric in $\Theta$.\footnote{$d_n^g(\theta,\theta)=0$ and $d_n^g$ satisfies the symmetric property and triangle inequality, but can be zero for two distinct points.} The logarithm of the $\varepsilon$-covering number of $B_r(\theta^\ast)$ with respect to $d_n^{g}$ is upper bounded by $C\,d\log(\frac{nd}{\varepsilon})$.
     
     \item For any $v\in \mb S^{d-1}$ and $\theta,\, \theta'\in B_r(\theta^\ast)$, it holds that $\mathbb{E} \big[\big(v^Tg(X,\theta)-v^Tg(X,\theta')\big)^2 \big]\leq C\,d^{\gamma_3}\,\|\theta-\theta'\|^{2\beta_1}$ and  $\mathbb{E}\big[\big(\ell(X,\theta)-\ell(X,\theta')-g(X,\theta')(\theta-\theta')\big)^2\big]\leq Cd^{\gamma_3}\, \|\theta-\theta'\|^{2+2\beta_1}$.
  
     \item  Let $\Delta_{\theta^\ast}=\mb E[\,g(X,\theta^\ast)\,g(X,\theta^\ast)^T]$ be the covariance matrix of the ``score vector'' $g(X,\theta^\ast)$. It holds that $\m H_{\theta^\ast}^{-1}\Delta_{\theta^\ast} \m H_{\theta^\ast}^{-1}\preceq Cd^{\gamma_4}\,I_d$.
  \end{enumerate}
  }

{
\noindent Conditions B.3.1 relaxes the pointwise differentiability requirement for the loss function $\ell(x,\theta)$ with respect to $\theta$. In fact,  in many statistical applications, the expectation in the population-level risk function $\m R(\theta)=\mb E[\ell(X,\theta)]$ has the smoothing effect of rendering $\m R$ to be twice differentiable. For instance, we can choose $g(x,\cdot)$ as the gradient (or any subgradient) of $\ell(x,\cdot)$ for $x\in\m X$ when $\ell$ is (or not) differentiable. Moreover,  the boundedness assumption on the covering number  in  Condition B.3.2 allows us to uniformly control the random fluctuation of the empirical mean $\frac{1}{n}\sum_{i=1}^n g(X_i,\theta)$ away from the gradient of $\m R(\theta)$. Condition B.3.3 can be interpreted as ``smooth'' assumptions on the loss function at the population level, quantified by $\beta_1$:  by taking expectations with respect to the data $X$, the first term controls the Lipschitz constant of $g(X,\cdot)$, while the second term controls the remainder term of the first-order Taylor expansion of $\ell(X,\cdot)$, where the gradient is replaced with $g(X,\cdot)$.   Condition B.3.4 assumes the boundedness of the operator norm of the matrix $\m H_{\theta^\ast}^{-1}\Delta_{\theta^\ast} \m H_{\theta^\ast}^{-1}$.  This matrix represents the limiting covariance matrix for the sampling distribution of the empirical risk minimizer $\wh\theta$, scaled by the sample size, i.e., $\sqrt{n}(\wh\theta-\theta)$ converges in distribution to $N_d(0,\m H_{\theta^\ast}^{-1}\Delta_{\theta^\ast} \m H_{\theta^\ast}^{-1})$.  This assumption allows us to provide an explicit bound on the deviance of $\wh\theta$, which represents the asymptotic mean of the Gibbs posterior, from $\theta^*$.  It is important to highlight that Conditions B.1-B.3 can cover the common scenario where the loss function is continuously twice differentiable (see Corollary~\ref{cor:smoothloss}). Furthermore, these conditions also apply to more general cases with non-smooth loss functions, such as quantile regression (see Corollary~\ref{co:quantile}).}

Additionally, we assume the following smoothness condition for the prior distribution and compactness of the parameter space.

   \vspace{0.5em}

  \noindent\textbf{Condition B.4 (Prior and parameter space):} {\it There exist positive $(n,d)$-independent constants $(C,r)$ so that the parameter space $\Theta$ satisfies  $B_r(\theta^\ast)\subset \Theta\subset [-C,C]^d$, and for any $\theta \in \Theta$, $\|\nabla (\log \pi) (\theta)\|\leq C\sqrt{d}$.}

\vspace{0.5em}
\noindent
{
The posterior density is defined to be zero for values of $\theta$ outside the parameter space $\Theta$, ensuring that MALA rejects any proposed states that go beyond the boundaries of $\Theta$. The assumption of compactness for the parameter space is primarily for technical convenience and is commonly made in Bayesian literature~\citep{10.1214/12-EJS675, yang2012bayesian}.  However, it is possible to relax this requirement by assuming the exponential tail behavior of the prior distribution, which will only incur extra logarithmic terms in the final result.   Finally, we made the following conditions to the preconditioning matrix $\wt I$. }

 \vspace{0.5em}
  
 \noindent\textbf{Condition C (Preconditioning matrix):} {\it There exist  some $(n,d)$-independent constants $C$ so that the preconditioning matrix  $\wt I$ satisfies that
 \begin{enumerate}
     \item  $\mnorm{\wt I^{-1}}_{\rm op}\mnorm{\wt I}_{\rm op}\leq C \mnorm{\m H_{\theta^*}}_{\rm op}\mnorm{\m H_{\theta^*}^{-1}}_{\rm op}$; 
     \item  ${\mnorm{\wt I}_{\rm op}\mnorm{(\wt I^{\frac{1}{2}}H_{\theta^*}\wt I^{\frac{1}{2}})^{-1}}_{\rm op}}\leq C\, \mnorm{\m H_{\theta^*}^{-1}}_{\rm op}$.
 \end{enumerate}
 }
 
The requirement for the preconditioning matrix $\wt I$ holds when $\wt I$ and its inverse has constant-order eigenvalues, such as the identity matrix that is conventionally used in MALA. On the other hand, it can also cover the case when $\wt I$ acts as a reasonable estimator to $\m H_{\theta^*}^{-1}$ (i.e, $\widetilde{I}^{1/2} \m H_{\theta^*} \widetilde{I}^{1/2}$ and its inverse  has constant-order eigenvalues). \\

\noindent We now state the following theorem that provides a mixing time bound for sampling from a Gibbs posterior using MALA. Note that the (sub)gradient $g$ is used for constructing the proposal in each step MALA.
  
  \begin{theorem}[\bf Complexity of MALA for Bayesian sampling]\label{th:Gibbsmixing}
   Consider sampling from the Bayesian Gibbs posteriors where $\m C_n(\theta;\, X^{(n)})=-n\,\alpha\,\m R_n(\theta)$.   Under Conditions B.1-B.4 and Condition C,  consider positive numbers $\rho_1,\rho_2$,  warming parameter $M_0$ and tolerance $\varepsilon$ satisfying (1) $\rho_1 I_d \preceq \widetilde{I}^{1/2} \m H_{\theta^*} \widetilde{I}^{1/2} \preceq \rho_2 I_d$; (2) $\log (\frac{M_0}{\varepsilon})\leq C_1\,(d^{\gamma_5}+\log n)$ for $(n,d)$-independent constants $C_1$ and $\gamma_5\geq 1$. There exists a constant $\kappa_1$ depends only on $(\beta_1,\gamma,\gamma_0,\gamma_1,\cdots,\gamma_5)$ so that
 if $d\leq c \frac{n^{\kappa_1}}{\log n}$ for a small enough constant $c$,  then with probability at least $1-n^{-1}$, the mixing time bound~\eqref{mixingbound}  in Theorem~\ref{thmala} holds for $\wt \varepsilon_0=1$ and
\begin{equation*}
    h=c_0\cdot\bigg[\rho_2\Big(d^{\frac{1}{3}}+d^{\frac{1}{4}}\big(\log \frac{M_0d\kappa}{\varepsilon}\big)^{\frac{1}{4}}+\big(\log \frac{M_0d\kappa}{\varepsilon}\big)^{\frac{1}{2}}\Big)\bigg]^{-1}, \text{ where } \kappa=\frac{\rho_2}{\rho_1},
     \end{equation*}
     where $c_0$ is an $(n,d)$-independent constant.
\end{theorem}

\begin{remark}
  {
  Theorem~\ref{th:Gibbsmixing} is proved by verifying Condition A for Bayesian Gibbs posteriors.  The parameter $\kappa_1$  sets an upper bound on how the dimensionality of the parameter space $d$ can grow in relation to the sample size $n$. A smaller $\kappa_1$ value implies that a larger dataset is necessary for the target posterior to be well-approximated by a Gaussian distribution. The expression for $\kappa_1$ is given by:
    \begin{equation}\label{eqn:kappa1}
 \begin{aligned}
      &\kappa_1
      = {\frac{1}{1+2\gamma+6\gamma_0+4\gamma_2+\gamma_4}}\wedge  {\frac{\beta_1}{1+\gamma_3+[(2\gamma_0)\vee ((\gamma_5+\gamma_0)(1+\beta_1))]}}\wedge \frac{1}{\gamma_0+\gamma_1+\gamma_5}\\
      &\wedge {\frac{1}{2\gamma+2\gamma_0+2\gamma_1+[2\vee (1+\gamma_4)]}}\wedge {\frac{1}{3\gamma_5+\gamma_0+[(2\gamma)\vee (\gamma_4+2\gamma_2+\gamma_0)\vee(2\gamma_2+2\gamma_0)]}}.
       \end{aligned}
 \end{equation}
    From the expression, $\kappa_1$ tends to be smaller if the loss function exhibits low smoothness, that means, $\beta_1$ is small.} The classical proof of the Gaussian approximation of Bayesian posteriors with smooth likelihoods is based on the Taylor expansion of the likelihood function around $\widehat\theta$~\citep[e.g.~see][]{Ghosh2003}. For the general non-smooth cases, we instead apply the Taylor expansion to the population level risk function $\m R$ and use chaining and localization techniques in the empirical process theory to relate it to the sample version. Moreover, we keep track of the parameter dimension dependence, making Theorem~\ref{th:Gibbsmixing} adaptable to more general cases under increasing dimension.
\end{remark}

\subsection{Gibbs posterior with smooth loss function}\label{sec:smoothloss}
 One representative example of Gibbs posterior satisfying Conditions B.1-B.4  is the one equipped with a smooth loss function.   More specifically, we  need Condition B.1 for the local convexity of the risk function, Condition B.4 for the smoothness of the prior  and the following smoothness condition to the loss function.
 
\vspace{0.5em}
 \noindent\textbf{Condition B.3' (Smoothness of loss function):}{\it There exist some $(n,d)$-independent constants $C>0$ and $(\gamma,\gamma_2,\gamma_3,\gamma_4)\geq 0$ so that (1) the loss function is twice differentiable so that for any $x\in \m X$ and $\theta\in \Theta$,  $\|\nabla_{\theta}\ell(x,\theta)\|\leq Cd^{\gamma}$; $\mnorm{{\rm Hess}_{\theta}(\ell(x,\theta))}_{\rm op}^2\leq  Cd^{\gamma_3}$;\footnote{We use $\nabla_{\theta}\ell(x,\theta)$ and ${\rm Hess}_{\theta}(\ell(x,\theta))$ to denote the gradient and Hessian matrix of $\ell_x(\cdot)=\ell(x,\cdot)$ evaluated at $\theta$, respectively. } and for any $\theta,\theta'\in \Theta$, $\vert\kern-0.25ex\vert\kern-0.25ex\vert{\rm Hess}_{\theta}(\ell(x,\theta))-{\rm Hess}_{\theta}(\ell(x,\theta'))\vert\kern-0.25ex\vert\kern-0.25ex\vert_{\rm op}\leq C\, d^{\gamma_2}\|\theta-\theta'\|$; (2) let $\Delta_{\theta^*}=\mb{E}[\nabla_{\theta} \ell(X,\theta^*)\nabla_{\theta} \ell(X,\theta^*)^T]$, then $\m H_{\theta^*}^{-1}\Delta_{\theta^*}\m H_{\theta^*}^{-1}\preceq C\, d^{\gamma_4}I_d $.}

\begin{corollary}[\bf Sampling from smooth posteriors]\label{cor:smoothloss}
 Consider the Bayesian Gibbs posterior with loss function $\ell$. Suppose (1) Conditions B.1, B.3' and B.4 hold; (2) the warming parameter $M_0$ and tolerance $\varepsilon$ satisfying  $\log (\frac{M_0}{\varepsilon})\leq C_1\,(d^{\gamma_5}+\log n)$ for $(n,d)$-independent constants $C_1$ and $\gamma_5\geq 1$; (3) $d\leq c\frac{n^{\kappa_1}}{\log n}$ for a small enough constant $c$, where $\kappa_1$ is defined in~\eqref{eqn:kappa1} with $\beta_1=1$. Then there exists an $(n,d)$-independent constant $c_0$ so that it holds with probability at least $1-n^{-1}$ that 
 \begin{enumerate}
     \item  consider the identity preconditioning matrix $\wt I=I_d$. the mixing time upper bound~\eqref{mixingbound} holds for any $\rho_1\leq \rho_2$ so that $\rho_1I_d\preceq\m H_{\theta^*}\preceq\rho_2I_d$, $\log (\frac{\rho_1}{\rho_2})\leq C_1 d^{\gamma_5}$ and
 \begin{equation*}
    h=c_0\cdot\bigg[\rho_2\cdot\Big(d^{\frac{1}{3}}+d^{\frac{1}{4}}\big(\log \frac{M_0d}{\varepsilon}\big)^{\frac{1}{4}}+\big(\log \frac{M_0d}{\varepsilon}\big)^{\frac{1}{2}}\Big)\bigg]^{-1};
     \end{equation*}
     \item  consider the inverse empirical Hessian matrix $\wt I=\big(|S|^{-1}\sum_{i\in S} {\rm Hess}_{\theta}(\ell(X_i,\wh\theta))\big)^{-1}$, where $S\subset \{1,2,\cdots,n\}$ with $|S|\geq C_2\,d^{\gamma_3+2\gamma_0+7/3}$ for a large enough $(n,d)$-independent constant $C_2$, then the mixing time upper bound~\eqref{mixingbound} holds with $\rho_1=\frac{1}{2}$ and
 \begin{equation*}
    h=c_0\cdot\bigg[\Big(d^{\frac{1}{3}}+d^{\frac{1}{4}}\big(\log \frac{M_0d}{\varepsilon}\big)^{\frac{1}{4}}+\big(\log \frac{M_0d}{\varepsilon}\big)^{\frac{1}{2}}\Big)\bigg]^{-1};
     \end{equation*}
moreover,  let $\mu_0=N_d(\wh\theta,n^{-1}\wt I)\big|_{\{\theta:\sqrt{n}\wt I^{-\frac{1}{2}}(\theta-\wh\theta)\|\leq 3c_1\sqrt{d}\}}$, where $c_1$  is a constant so that
   $c_1\geq 3\vee \underset{i\in [d],j\in [d]}{\sup} \frac{\partial^2\m R(\theta^*)}{\partial \theta_i\partial\theta_j} $, then $\mu_0$ is $M_0$-warm with respect to $\pi_n$ with $\log M_0\leq d^{\frac{1}{3}}$.
 \end{enumerate}
 
 \end{corollary}
 
 When the Hessian matrix $\m H_{\theta^*}$ is ill-conditioned, introducing the preconditioning matrix $\wt I=\big(|S|^{-1}\\\sum_{i\in S} {\rm Hess}_{\theta}(\ell(X_i,\wh\theta))\big)^{-1}$ may lead to a faster mixing. Furthermore, if the tolerance satisfying $ \log (\frac{1}{\varepsilon})=\m O(d^{\frac{1}{3}})$,  then the second statement of Corollary~\ref{cor:smoothloss} can lead to an optimal mixing time bound $\m O\big(d^{\frac{1}{3}}\log(\frac{1}{\varepsilon})\big)$.

\subsection{Bayesian quantile regression}\label{Sec:Quantile_reg}
We consider Bayesian quantile regression as a representative example where the loss function is non-smooth.
 Specifically, in quantile regression~\citep{10.2307/1913643}, for a fixed $\tau \in (0,1)$, the $\tau^{th}$ quantile $q_{\tau}(Y|\widetilde X)$ of the response $Y\in \mathbb{R}$ given the covariates $\widetilde{X}\in \mathbb{R}^d$ is modelled as $q_{\tau}(Y|\widetilde{X})=\widetilde{X}^T\theta^\ast$. Here we consider the homogeneous case where the error $e=Y-\wt X^T\theta^*$  is independent of the covariates $\wt X$. Given a set of $n$ i.i.d. samples $X^{(n)}=\{X_i=(\widetilde X_i,Y_i)\}_{i\in[n]}$, the quantile regression solves the following convex optimization problem: 
 $$
 \widehat{\theta}=     {\arg\min}_{\theta\in \Theta} \sum_{i=1}^n \Big[(Y_i-\widetilde{X}_i^T\theta) \cdot \big(\tau-\mathbf{1} (Y_i<\widetilde{X}_i^T\theta) \big)\Big],
 $$
 where the loss function $\ell_{q}\big((\widetilde X, Y),\theta\big)=(Y-\widetilde{X}^T\theta)\cdot\big(\tau-\mathbf{1} (Y<\widetilde{X}^T\theta)\big)$ is referred to as the check loss. The minimization of the check loss function is equivalent to the maximization of a likelihood function formed by combining independently distributed asymmetric Laplace densities~\citep{YU2001437}. The posterior for Bayesian quantile regression can thus be formed by assuming a (possibly misspecified) asymmetric Laplace distribution (ALD) for the response, which is  
 $$
 \pi_{n}(\theta)\propto \exp\big(-n\,\m R_n(\theta)\big)\,\pi(\theta),\quad \theta\in\mb R^d,
 $$
 with $\pi(\theta)$ being a prior on $\Theta$ and $\m R_n(\theta)=n^{-1}\sum_{i=1}^n \ell_q(X_i,\theta)$ being the empirical risk function. Furthermore,  by adding a multiplier $\alpha>0$ to the likelihood, we can obtain the Gibbs (or tempered) posterior.
 
Since the loss function $\ell_q(X,\theta)$ for quantile regression is not differentiable when $Y=\widetilde X^T\theta$, in order to sampling from the Gibbs posterior associated with Bayesian quantile regression  using the (preconditioned) MALA, we need to consider the subgradient of $\ell_q$ with respect to $\theta$, given by 
$$
g(X,\theta)=\big(\mathbf{1}(Y<\widetilde X^T\theta)-\tau\big)\, \widetilde X, \quad X=(\widetilde X, Y), \ \ \theta\in \mb R^d.
$$
The following corollary quantifies the computational complexity for sampling from $\pi_{n}$ using MALA. We first state the required conditions.
\vspace{0.5em}

\noindent\textbf{Condition D.1:}  {\it There exist $(n,d)$-independent constants $(C,C') >0$ and $(\alpha_0,\alpha_1)\geq 0$ such that  (1) the support $\m X$ of the covariates $\widetilde X$ is included in $[-C,C]^d$; (2) for any $v\in \mb S^{d-1}$, $  \mb E|\widetilde{X}^Tv|^2\geq C'd^{-\alpha_0}$ and $\mb E|\widetilde{X}^Tv|^3\leq Cd^{\alpha_1}$.}

\vspace{0.5em}
\noindent\textbf{Condition D.2:} {\it Let $f_{e}(\cdot)$ denote the probability density function of the homogeneous error $e=Y-\widetilde X^T\theta^\ast$, then there exist $(n,d)$-independent constants $(C,C') >0$ such that  (1) $\int_{-\infty}^{0} f_{e}(z)dz=\tau$; (2) $f_{e}(0)>C'$  and $\sup_{e\in \mb R^d}f_{e}(e)\leq C$; (3) for any $e_1,e_2\in \mb R$, $|f_{e}(e_1)- f_{e}(e_2)|\leq C|e_1-e_2| $.}

\vspace{0.5em}

\noindent Condition D.1 assumes the compactness of the covariate space and the positive definiteness of the gram matrix $\mathbb{E}[\widetilde{X}\widetilde{X}^T]$.
Condition D.2 introduces several regularity conditions on the distribution of the error $e = Y - \widetilde{X}^T\theta^\ast$: (1) The error term $e$ is independent of the covariates. (2) The model is correctly specified, meaning that $\widetilde{X}^T\theta^*$ corresponds to the $\tau$-th quantile of the response variable $Y$ given $\wt X$. (3) The density function $f_e(\cdot)$ of the error term is positive at the origin and Lipschitz continuous. Under the assumption of homogeneous errors, the limiting covariance matrix of the posterior distribution of interest is given by $n^{-1}(f_e(0)\cdot \mb{E}[\wt X\wt X^T])^{-1}$. In this case, a natural choice for the preconditioning matrix is the inverse of the empirical Gram matrix, denoted as  $\widetilde I=\big(|S|^{-1}\sum_{i\in S} \widetilde X_i\widetilde X_i^T\big)^{-1}$ where $S\subset\{1,2,\cdots,n\}$. It is worth noting that similar analyses can be carried out for the case of heterogeneous errors, but the limiting covariance matrix will be more complex.

 \begin{corollary}[\bf Sampling from non-smooth posteriors]\label{co:quantile}
 Suppose Conditions D.1, D.2, and B.4 are satisfied, and the warming parameter $M_0$ and tolerance $\varepsilon$ satisfying  $\log (\frac{M_0}{\varepsilon})\leq C_1\,(d^{\alpha_2}+\log n)$ for $(n,d)$-independent constants $C_1$ and $\alpha_2\geq 1$. Assume $d\leq c(\frac{n^{\widetilde \alpha}}{\log n})$ with $\widetilde \alpha=\frac{1}{2+4\alpha_1+7\alpha_0}\wedge \frac{1}{2+3\alpha_0+2\alpha_1+3\alpha_2}$ and a small enough constant $c$, and let the inverse empirical Gram matrix $\widetilde I=\big(|S|^{-1}\sum_{i\in S} \widetilde X_i\widetilde X_i^T\big)^{-1}$ be the preconditioning matrix,  where $S\subset \{1,2,\cdots,n\}$ with $|S|\geq C_2\, d^{\alpha_1+2\alpha_0+3/2}\log n$ for a large enough $(n,d)$-independent constant $C_2$, then it holds with probability larger than $1-\frac{1}{n}$ that that the mixing time upper bound~\eqref{mixingbound} is true with $\rho_1=\frac{1}{2}f_{e}(0)$ and
 \begin{equation*}
    h=c_0\cdot\bigg[{f_{e}(0)}\cdot\Big(d^{\frac{1}{3}}+d^{\frac{1}{4}}\big(\log \frac{M_0d}{\varepsilon}\big)^{\frac{1}{4}}+\big(\log \frac{M_0d}{\varepsilon}\big)^{\frac{1}{2}}\Big)\bigg]^{-1}
     \end{equation*}
 with $c_0$ being an $(n,d)$-independent constant.
 \end{corollary}

 {
 Corollary~\ref{co:quantile} illustrates the implications of applying our theory to non-smooth posteriors. A key observation is that in the large-sample regime, although the potential function associated with the Bayesian posterior may be non-smooth, its population-level counterpart is smooth (as per Condition B.3). This allows MALA, using sub-gradients, to effectively sample from non-smooth posteriors.
 Moreover, while our theory is applicable to non-smooth posteriors, the smoothness of the posterior density function influences its convergence to a Gaussian limit as $n$ grows, as captured by the parameter $\beta_1$ in Condition B.3. A posterior with higher smoothness (or larger $\beta_1$) will converge more rapidly to a Gaussian distribution, as demonstrated in Lemma~\ref{th1}, which in turn leads to an improved (higher) acceptance rate of MALA. For example, in Bayesian quantile regression, the smoothness parameter $\beta_1$ is at most $\frac{1}{2}$. In contrast, for posterior densities with smooth loss functions, $\beta_1$ can be taken as $1$. Interestingly, our theoretical result also leads a practical guideline: when applying MALA to sample from a less smooth Bayesian posterior densities, a relatively larger sample size $n$ is needed to maintain the sampling efficiency. Otherwise, if $n$ is not sufficiently large relative to the dimension, the non-smoothness of the Bayesian posterior can result in a lower acceptance rate and slower mixing times for MALA; see our simulation results in Section 7 for some empirical evidence. }

\section{Proof Sketch of Theorem~\ref{thmala}}\label{sec:proof_sketch}
  In this section, we provide a sketched proof about how to utilize the general machinery of $s$-conductance profile developed in Section~\ref{mixingboundscp} to analyze the mixing time of MALA under Condition A. We consider the identity preconditioning matrix (i.e. $\wt I=I_d$) in this sketch for simplicity, and the case for general preconditioning matrix can be proved by considering the transformation $G(\theta)=\sqrt{n}\widetilde{I}^{-\frac{1}{2}}(\theta-\wh\theta)$, see Appendix~\ref{Proof1} for further details. 
  
  Let $T^{\zeta}_x(\dd y)=T^{\zeta}(x,\dd y)$ denote the Markov transition kernel of the $\zeta$-lazy version of MALA for sampling from $\pi_{\rm loc}$ as described in Section~\ref{sec:MALA_mixtime} with rescaled step size $h$. To apply Lemma~\ref{lemma:conductance_informal}, we first need to establish a log-isoperimetric inequality, which is a property of $\pi_{\rm loc}$ alone and is not specific to MALA. This step can be done by adapting existing proofs of a log-isoperimetric inequality for Gaussians~(e.g.~Lemma 16 of~\cite{JMLR:v21:19-441}) to $\pi_{\rm loc}$ via a perturbation analysis (see Lemma~\ref{lemmalogiso} and its proof in the appendix for details). Second, we need to apply an overlap argument for bounding the total variation distance between $T^{\zeta}_x(\cdot)$ and $T^{\zeta}_z(\cdot)$ for $x$ and $z$ satisfying $\|x-z\|\leq C \sqrt{h}$ and belonging to a high probability set $E$ under $\pi_{\rm loc}$. This step utilizes the structure and properties of MALA algorithm, and we briefly sketch its proof below (details can be found in Lemma~\ref{boundTV} in the appendix) and discuss its difference from existing proofs. 
  
{
We construct a high probability set as $E=\{\xi\in B_{R/2}^d: \big|\xi^T\wt J^3\xi-{\rm tr}(\wt J^2)\big|\leq r_d\}\cap \{\xi\in B_{R/2}^d: \big|\xi^T\wt J^2\xi-{\rm tr}(\wt J)\big|\leq r_d/\rho_2\}$, where the value of $r_d$ makes $\pi_{\rm loc}(E)\geq 1-2\frac{h\rho_1\varepsilon^2}{M_0^2}$ based on the last property of Condition A (details can be found in Lemma~\ref{lemmaprobE}). Recall the acceptance probability $A(x,y)=1\wedge \frac{\pi_{\rm loc}(y)\,Q(y,x)}{\pi_{\rm loc}(x)\,Q(x,y)}$ and denotes  
 $\ov A(x,y)=1\wedge \frac{\ov \pi(y)\, Q(y,x)}{\ov \pi(x)\, Q(x,y)}$ with $\ov\pi$ being the density of the Gaussian $N_d(0,J^{-1})$. By comparing $\pi_{\rm loc}$ and $\ov \pi$ using Condition A, we can get the following inequality: 
\begin{equation}\label{decompositionTV}
    \begin{aligned}
      &\ \|T^{\zeta}_{x}- T^{\zeta}_{z}\|_{TV}\leq  1-(1-\zeta)\int_{ B_{R}^d}\min\Big( A(x,y)  Q(x,y),   A(z,y)  Q(z,y)\Big)\,\dd y \\
 &\leq 1-\frac{1}{2}(1-\zeta) \exp(-2\wt\varepsilon_0) \cdot \bigg(\int_{ B_{R}^d}\ov A(x,y)Q(x,y)\,\dd y+\int_{ B_{R}^d}\ov A(z,y)  Q(z,y)\,\dd y\\
      &\qquad-\int_{ B_{R}^d}\big|\ov A(x,y)  Q(x,y)-  \ov A(z,y)  Q(z,y)\big|\,\dd y\bigg)\\
      & \leq 1-(1-\zeta) \exp(-2\wt\varepsilon_0) \cdot \Big(1-\int_{B_R^d}  Q(x,y)(1-\ov A(x,y))\,\dd y-\int_{B_R^d}   Q(z,y)(1-\ov A(z,y))\,\dd y\\
       &\qquad -\| Q_x-  Q_z\|_{  \rm TV}-\frac{1}{2}\int_{ (B_{R}^d)^c}  Q(x,y)\,\dd y-\frac{1}{2}\int_{ (B_{R}^d)^c}   Q(z,y)\,\dd y\Big).
    \end{aligned}
     \end{equation}
We will separately bound the terms on the right hand side of~\eqref{decompositionTV} as follows. The last term $\frac{1}{2}\int_{ (B_{R}^d)^c}  Q(x,y)\,\dd y+\frac{1}{2}\int_{ (B_{R}^d)^c}   Q(z,y)\,\dd y$ can be upper bounded by $\frac{1}{6}$ using the condition of $R$ in Condition A. For the remaining terms, let $Q_x$ denote the probability measure with density function $Q(x,\cdot)$, now we use Condition A by comparing $Q_x$ with the proposal distribution 
 $$
 Q^{\Delta}_x:=N_d(x-hJx,\,2hI_d)
 $$ 
 of MALA for sampling from the Gaussian $N_d(0,J^{-1})$, leading to
 \begin{equation}\label{decompositionTV3.1}
    \begin{aligned}
     \int_{B_R^d}Q(x,y)\,\big(1-\ov A(x,y)\big)\, \dd y&\leq 2\,\|Q_x-Q^{\Delta}_x\|_{  \rm TV}+ {\int_{\mb R^d}\Big|Q^{\Delta}(x,y)-\frac{\ov\pi(y)Q^{\Delta}(y,x)}{\ov\pi(x)}\Big|\,\dd y} \\
     &+ \int_{B_R^d}\Big|\frac{\ov\pi(y)Q^{\Delta}(y,x)}{\ov\pi(x)}-\frac{\ov \pi (y)Q(y,x)}{\ov\pi(x)(x)}\Big|\,\dd y,
    \end{aligned}
\end{equation}
where we use $Q^{\Delta}(x,\cdot)$ to denote the density function of $Q^{\Delta}_x$.
It  then can be proved using Condition A and Pinsker's inequality after some careful calculations (see Lemmas~\ref{lemma3} and~\ref{boundCD} in the appendix) that 
$$
\int_{B_R^d}  Q(x,y)(1-\ov A(x,y))\,\dd y+\int_{B_R^d}   Q(z,y)(1-\ov A(z,y))\,\dd y+\| Q_x-  Q_z\|_{  \rm TV}\leq 1/3.
$$ 
Our proof of Lemma~\ref{lemma3} for bounding $\int_{\mb R^d} \big|Q^{\Delta}(x,y)- {\ov\pi(y)Q^{\Delta}(y,x)}/{\ov\pi(x)}\big|\,\dd y$ is technically similar to that of Proposition 38 in~\cite{pmlr-v134-chewi21a} for bounding the mixing time of MALA with a standard Gaussian target (i.e.~$\overline{\pi}=N_d(0,I_d)$). The non-trivial part in our analysis lies in keeping track of the dependence on the maximal and minimal eigenvalues of $J$. 
 Finally, we can obtain 
 \begin{equation*}
       \| T^{\zeta}_{x}- T^{\zeta}_{z}\|_{TV}\leq  1-\frac{1-\zeta}{2} \exp(-2\wt\varepsilon_0).
   \end{equation*}
With the lower bound on $\pi_{\rm loc}(E)$ and the upper bound on $\| T^{\zeta}_{x}- T^{\zeta}_{z}\|_{TV}$, we are then able to apply the $s$-conductance profile argument to control the mixing time. 
 }

 \begin{remark}
      It is worth mentioning that the analysis in~\cite{JMLR:v21:19-441} requires the high probability set, which is set $E$ in our case, to be convex. This requirement will deteriorate the $d$ dependence of the mixing time bound since $\|T^\zeta_x-T^\zeta_z\|_{\rm TV}$ for $x,z\in E$ can no longer be controlled under a large step size $h$ as ours. This motivates us to introduce the more flexible notion of \emph{$s$-conductance profile} that extends the commonly used conductance profile~\citep{10.1214/EJP.v11-300,JMLR:v21:19-441} and  $s$-conductance~\citep{https://doi.org/10.1002/rsa.3240040402}. Analysis based on the $s$-conductance profile leads to a better warming parameter dependence than that obtained in~\cite{pmlr-v134-chewi21a,10.2307/30243694} without affecting our obtained dimension dependence (based on $s$-conductance). A complete proof of this theorem is included in Appendix~\ref{Proof1}. Similar analysis can also be carried over for analyzing general smooth and strictly log-concave densities to improve the warming parameter dependence~\citep[e.g.][]{pmlr-v134-chewi21a,10.2307/30243694}.
 
 \end{remark}

{
\section{Numerical Study}\label{sec:num}
In this section, we conduct an empirical study  to explore how the performance of MALA varies across different dimensions and sample sizes when targeting different Bayesian posteriors.

 \subsection{Set up}
 
We carry out the experiment using two examples: Bayesian linear regression and Bayesian median regression. For Bayesian linear regression, the corresponding Bayesian posterior is given by: 
$$
\pi_n^ {\rm mean}(\theta\,|\,X^{(n)})\propto \exp\Big(-\frac{1}{2}\sum_{i=1}^n \big\|Y_i-\widetilde{X}_i^T\theta\big\|^2\Big)\,\pi(\theta), \ \ \theta\in\mb R^d.
$$
For Bayesian median regression, the Bayesian posterior is given by 
\begin{equation*}
    \pi_n^ {\rm med}(\theta\,|\,X^{(n)})\propto \exp\Big(-\frac{1}{2}\sum_{i=1}^n \big|Y_i-\wt X_i^T\theta\big|\Big)\,\pi(\theta), \ \ \theta\in\mb R^d.
\end{equation*}
We choose the parameter dimension $d$ from the set $\{15,20,30,40,\cdots,100\}$ and sample size $n$ from $\{500,1000,2000,5000,500(d/15),500(d/15)^{3/2},500(d/15)^2\}$. The covariates $\widetilde X$ are generated from a multivariate Gaussian distribution  with zero mean and identity covariance matrix. For Bayesian linear regression, we generate a random error variable $e$ follows a standard normal distribution, and for Bayesian median regression, $e$ follows a Laplace distribution with location parameter $\mu=0$ and scale parameter $b=2$. The response variable $Y$ is  given by  $Y=\widetilde X^T\theta^\ast+e$ with $\theta^\ast=(1,1,\cdots,1)$. We consider the parameter space $\Theta=[-100,100]^d$ and the prior is chosen to be a uniform distribution over $\Theta$. We then use MALA to sample from the Bayesian posterior $\pi_n^ {\rm mean}$ and $\pi_n^ {\rm med}$.

 \subsection{Results}

 \begin{figure}[tp]
  
\centering  
\vspace{1em}
\subfigure[Acceptance probability ($\pi^{\rm med}_n$)]{
 \includegraphics[trim={0.1cm 0.1cm 0.1cm 0.25cm}, clip,
 width=0.46\textwidth]{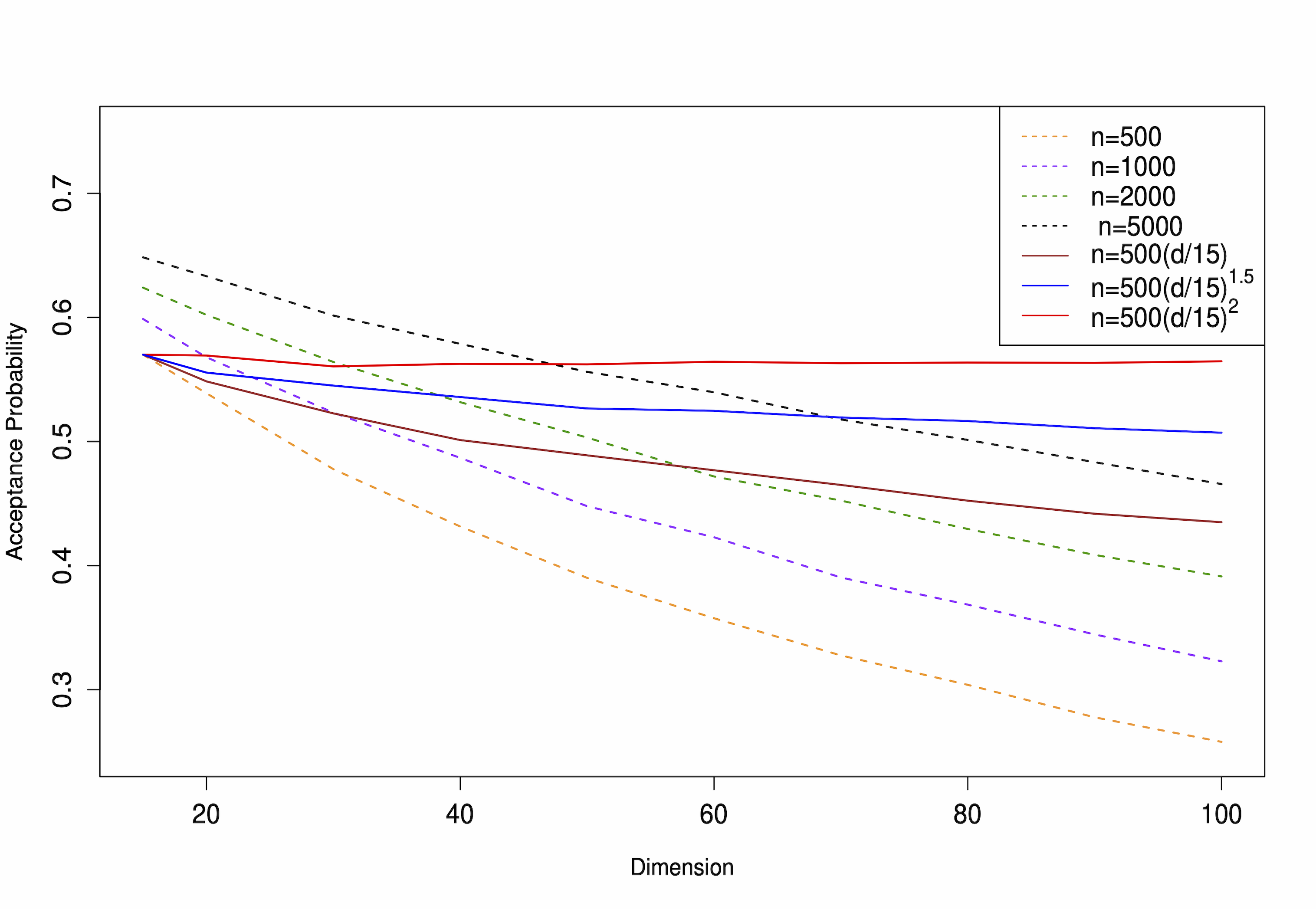}
 \label{Fig:ACCP_q}}
\subfigure[Log effective sample size ($\pi^{\rm med}_n$)]{
 \includegraphics[trim={0.1cm 0.1cm 0.1cm 0.25cm}, clip, width=0.46\textwidth]{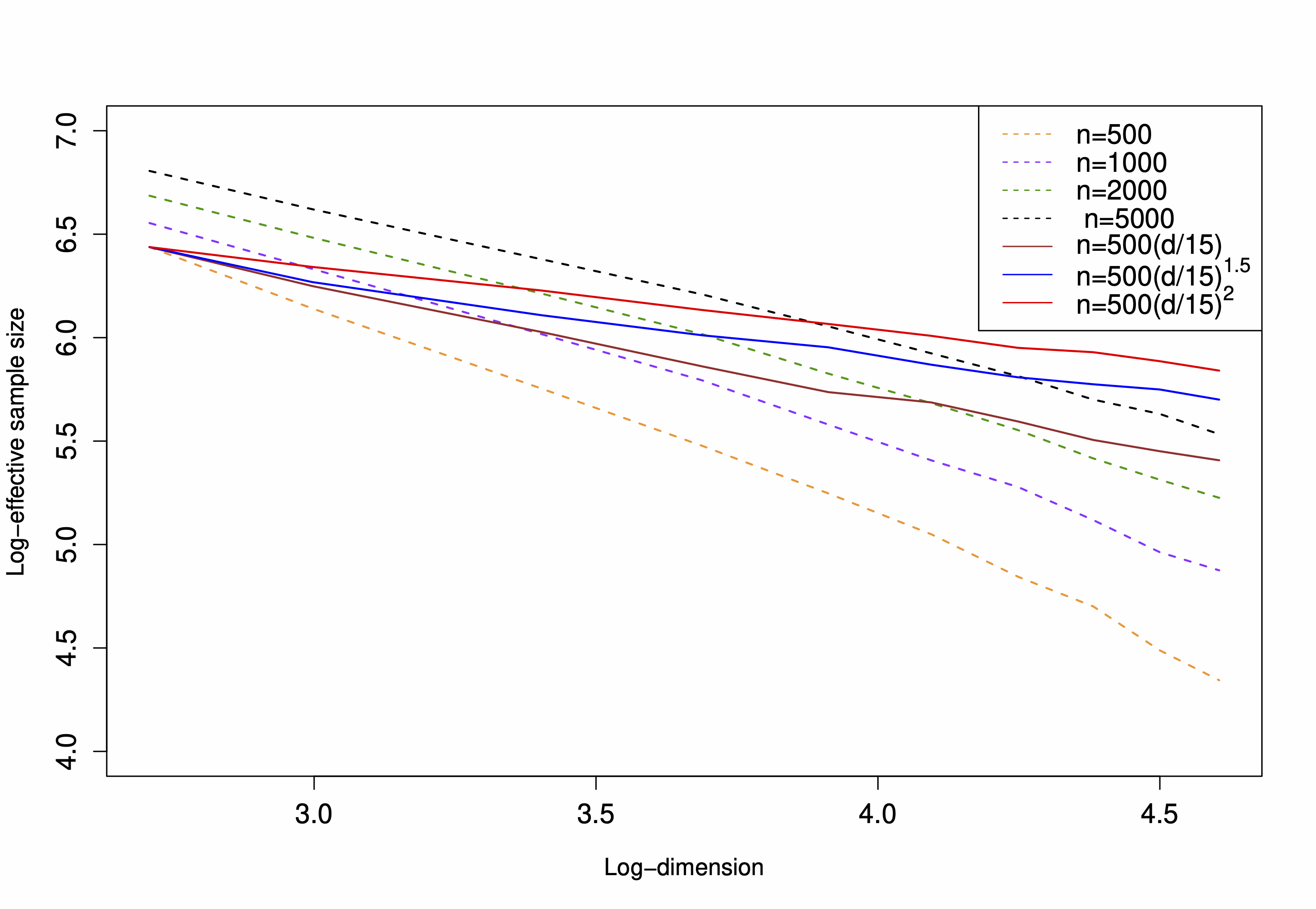}
 \label{Fig:ESS_q}
 \vspace{2em}}
\vspace{1em}
\subfigure[Acceptance probability ($\pi^{\rm mean}_n$)]{
 \includegraphics[trim={0.1cm 0.1cm 0.1cm 0.25cm}, clip,
 width=0.46\textwidth]{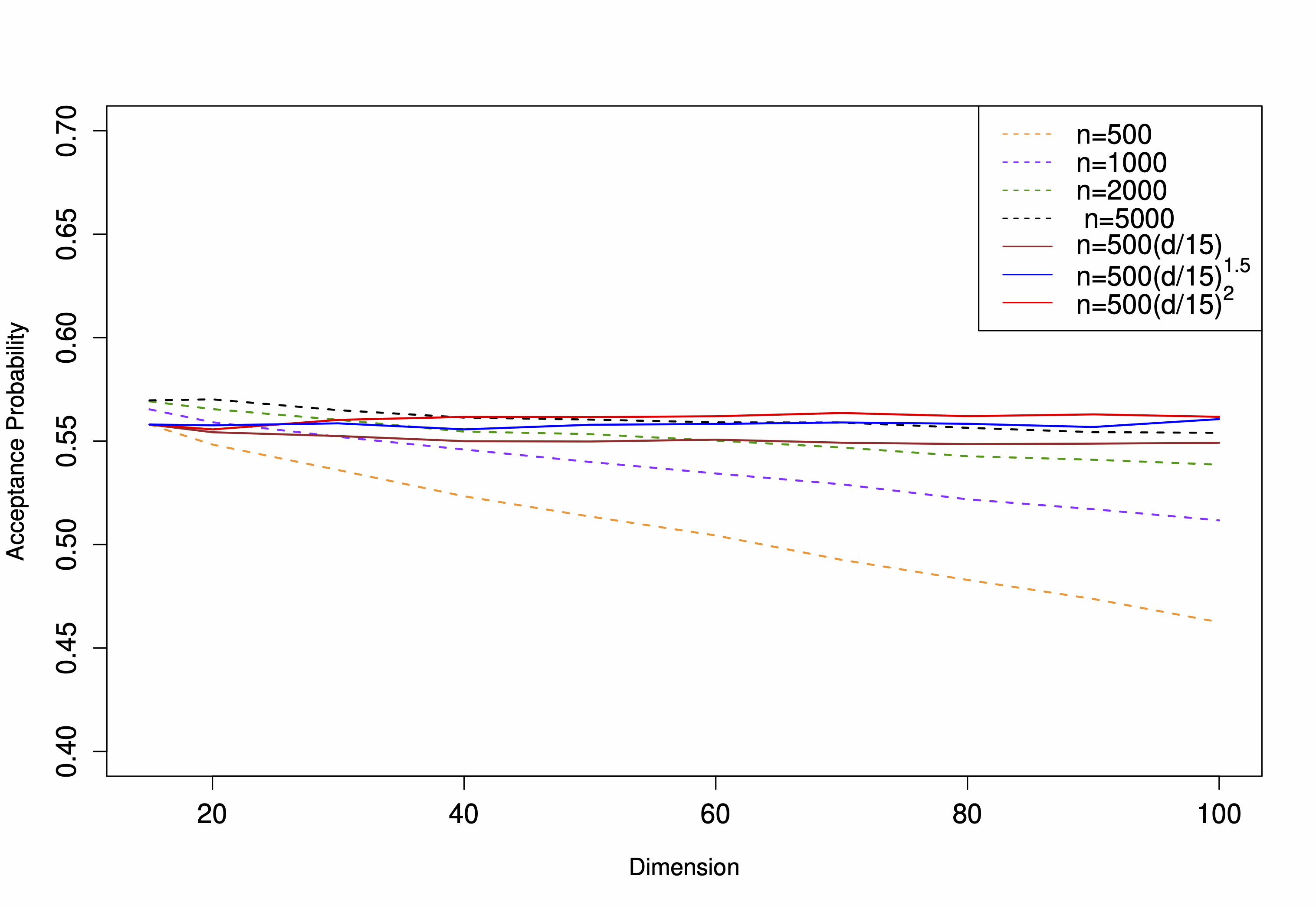}
 \label{Fig:ACCP_r}}
\subfigure[Log effective sample size ($\pi^{\rm mean}_n$)]{
 \includegraphics[trim={0.1cm 0.1cm 0.1cm 0.25cm}, clip, width=0.46\textwidth]{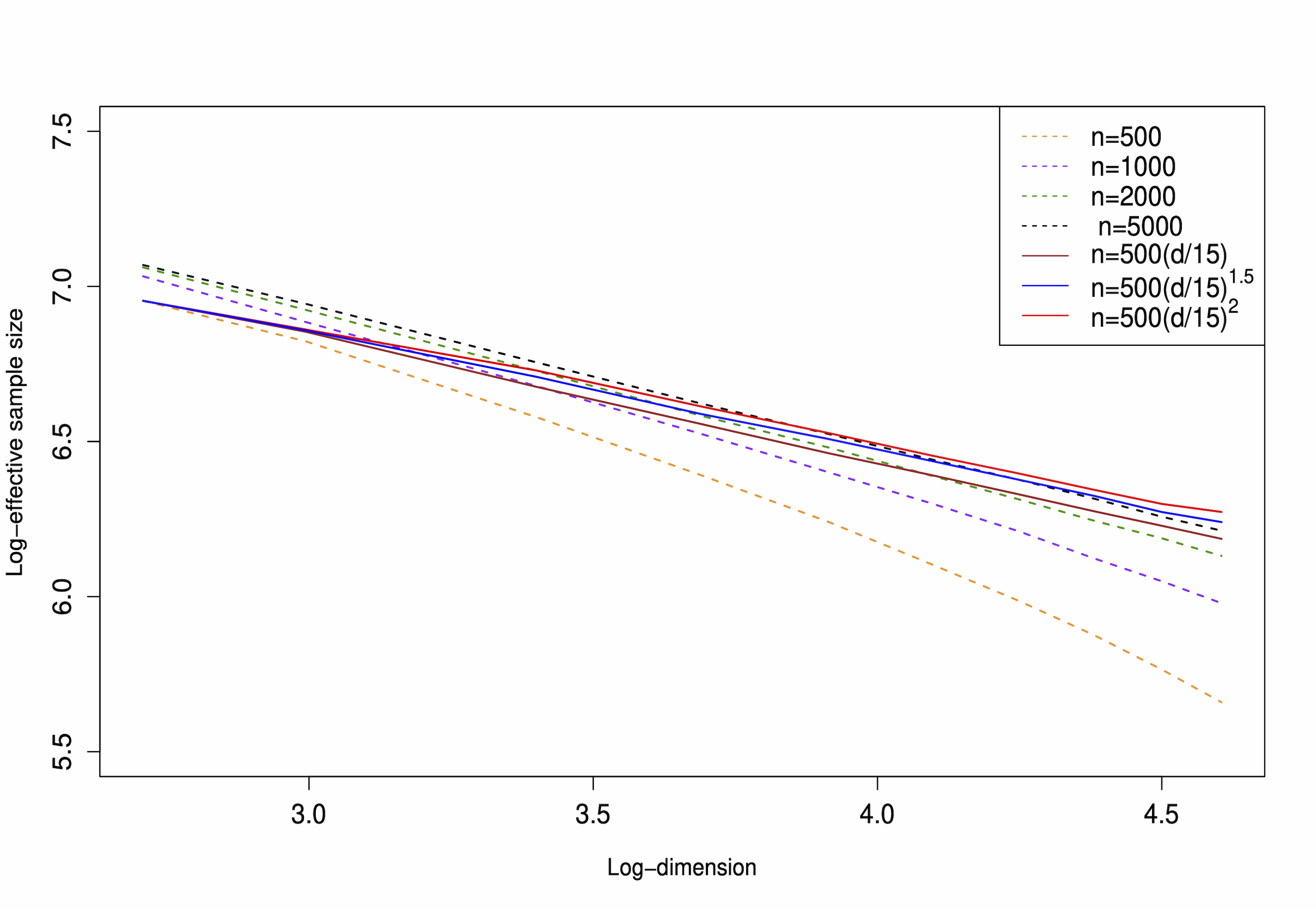}
 \label{Fig:ESS_r}
 \vspace{2em}}
 
 \caption{ Plots (a) and (c) report the average acceptance probabilities of MALA when sampling from the posterior in Bayesian quantile regression (denoted as $\pi^{\rm med}_n$) and  Bayesian linear regression (denoted as $\pi^{\rm mean}_n$) respectively, across various sample sizes ($n$) and dimensions ($d$), with the step size $\wt h=c d^{-\frac{1}{3}}n^{-1}$. Plots (b) and (d) present the relationship between the logarithm of the effective sample size and the logarithm of the dimension for sampling from $\pi^{\rm med}_n$ and $\pi^{\rm mean}_n$ respectively.  As we can see, when the sample size increases with the dimension at a rate of $d^{2}$, by choosing steps sizes with scaling $d^{-\frac{1}{3}}n^{-1}$, the acceptance probabilities roughly remain constant and the change in the logarithmic effective sample sizes exhibit  slopes close to $-\frac{1}{3}$ for both examples of Bayesian linear regression and median regression. On the other hand, when $n$ remains a constant, in both cases, the acceptance probabilities will decrease as $d$ becomes larger, and  the changes in the logarithmic effective sample size exhibit slopes smaller than $-\frac{1}{3}$. However, compared to $\pi^{\rm med}_n$, the decreases in acceptance probability and effective sample size are  much slower for sampling from $\pi^{\rm mean}_n$. In particular, when $n$ increases with $d$ at a linear rate, the acceptance probabilities for $\pi^{\rm mean}_n$ roughly remain constant, while there is obvious decrease in the acceptance probabilities for $\pi^{\rm med}_n$.}
 
\label{Fig_mean}
\end{figure}
In general, estimating the mixing time of a Markov chain is a challenging task. Instead, we utilize the effective sample size~\citep{gelman1995bayesian} as a metric to assess the mixing of MALA.  The effective sample size of $N$ Markov samples, denoted as $N_{\rm eff}$, quantifies the amount of information lost due to correlations in the chain, and plays a role similar to the number of independent draws in the standard central limit theorem~\citep{brooks2011handbook}.   The effective sample size of a sequence is formally defined in terms of the autocorrelations within the sequence at different lags,  i.e.,  $N_{\text {eff }}=\frac{N}{1+2 \sum_{t=1}^{\infty} \rho_t}$, with $\rho_t$ being the autocorrelation  at lag $t$.  Details of the estimation of $N_{\rm eff}$ can be found in Section 11.5 of~\cite{gelman2013bayesian}. It is worth noting that, theoretically, the ratio $\frac{N}{N_{\rm eff}}$ can be controlled by the inverse of the spectral gap~\citep{kloeckner2019effective}, which governs the convergence of the Markov chain.

Taking into account our theoretical findings regarding the convergence of  MALA with an appropriate warm start, we compute the effective sample sizes after a burn-in period of $1000$ iterations, totaling $5000$ iterations. We choose the step size $\wt h=c_1 d^{-\frac{1}{3}}n^{-1}$, where $c_1=4.28$ for Bayesian median regression and $c_1=1.39$ for Bayesian linear regression. These choices of $c_1$ ensure that the overall acceptance probability in each example closely approximates $0.574$ as suggested by~\cite{roberts1998optimal}. The preconditioning matrix $\wt I$ is chosen to be the identity matrix. 

Figure~\ref{Fig_mean}  present the trends of the average acceptance probability and the logarithm of the effective sample size when sampling from $\pi^{\rm med}_n$ and $\pi^{\rm mean}_n$,  considering varying sample sizes and dimensions. When $n$ remains unchanged for varying $d$, we observe a decrease in the acceptance probability as $d$ grows larger in both cases. Additionally, the trends of the logarithmic effective sample size exhibit slopes smaller than $-\frac{1}{3}$. The reason for this phenomenon is that, the deviance $\wt{\varepsilon}_0$ of the target posterior from the Gaussian distribution, stated in Theorem~\ref{thmala}, will increase with $d$ when the sample size remains unchanged. Consequently, when $d$ is sufficiently large, the mixing time will deviate significantly from $\m O(d^{\frac{1}{3}})$ and the acceptance probability will decrease rapidly when employing a step size of order $d^{-\frac{1}{3}}n^{-1}$. Another interesting observation is that the  decreases in acceptance probability and effective sample size are  much slower when sampling from $\pi^{\rm mean}_n$ compared to sampling from $\pi^{\rm med}_n$. One factor results in this phenomenon can be the smoothness of the loss function used in $\pi^{\rm mean}_n$, which aids the convergence of the Gibbs posterior to the Gaussian distribution. Specifically, Lemma~\ref{th1} in Appendix~\ref{sec:proofmixing} demonstrates that a Gibbs posterior with  a smooth loss function will converge to a Gaussian distribution with a rate of $\m O(n^{-1/2})$ for a fixed $d$, while the Gibbs posterior used in Bayesian quantile regression approaches a Gaussian distribution at a rate of $\m O(n^{-1/4})$. Therefore, under the same $n$ and $d$, the approximation error $\wt{\varepsilon}_0$ for $\pi^{\rm mean}_n$ is much smaller than $\pi^{\rm med}_n$. Additionally, we can see from Figure~\ref{Fig_mean} that, for achieving a constant acceptance probability and effective sample size at an order of $d^{-\frac{1}{3}}$ when $d$ ranges from $15$ to $100$,  the condition $d=\m O(\sqrt{n})$ is required for sampling from $\pi^{\rm med}_n$, while the condition $d=\m O(n)$ suffices for sampling from $\pi^{\rm mean}_n$.}

  \section{Conclusion and Discussion} 
 In this paper, we studied the sampling complexity of Bayesian (pseudo-)posteriors using MALA under large sample size, covering cases where the posterior density is non-smooth and/or non-log-concave. A variant of MALA that includes a preconditioning matrix was also considered. While our analysis for the preconditioned MALA suggests an adaptive MALA with a data-driven preconditioning matrix may be preferable, its rigorous theoretical analysis may leave as our future work.  When applying our main result to Bayesian inference, we mainly considered the Gibbs posterior, while similar analysis may carry over to other types of Bayesian pseudo-posterior, such as Bayesian empirical likelihood~\citep{10.1093/biomet/90.2.319}, and we leave this for future research. {
 Another challenge lies in constructing a suitable warm start that satisfies $\log M_0 \leq d^{\frac{1}{3}}$. Obtaining a warm start efficiently for general non-log-concave sampling can be challenging. However, the asymptotic Gaussian nature of the Bayesian posterior may aid in the construction of such a warm start, and it is possible to develop specific algorithms tailored to particular problems that leverage the Gaussian asymptotics. For instance, in Bayesian quantile regression, one can determine the point estimator $\hat{\theta}$ using linear programming and utilize the Gaussian asymptotic properties of the posterior to construct initializations. A more detailed exploration of this topic is left for future research.}

 \bibliographystyle{plainnat}
 
 \bibliography{ref} 
  
 \newpage
\appendix
\begin{center}
{\bf\Large Appendix}
\end{center}

 
 
We summarize some necessary notation and definitions in the appendix. We use $\mathbf{1}_{A}$ to denote the indicator function of a set $A$ so that $\mathbf{1}_{A}(x) = 1$ if $x\in A$ and zero otherwise. For two sequences $\{a_n\}$ and $\{b_n\}$, we use the notation $a_n \lesssim b_n$ and $a_n \gtrsim b_n$ to mean $a_n \leq Cb_n$ and $a_n \geq C b_n$, respectively, for some constant $C>0$ independent of $n,d$. In addition, $a_n \asymp b_n$ means that both $a_n \lesssim b_n$ and $a_n\gtrsim b_n$ hold, and $a_n=\m O(b_n)$ if $a_n\lesssim b_n$; $a_n=\Theta(b_n)$ if $a_n\asymp b_n$. We use $\mathbf N(\mathcal{F},d_n,\varepsilon)$ to denote the $\varepsilon$-covering number  of $\m F$ with respect to pseudo-metric $d_n$. Throughout, $C$, $c$, $C_0$, $c_0$, $C_1$, $c_1$, \ldots are generically used to denote positive constants independent of $n,d$ whose values might change from one line to another. We denote $\m L^2(\pi)$ to be the space of square integrable functions under measure $\pi$. For a transition kernel  $T:  \Theta\times \m B(\Theta)\to \mb R $ of a reversible Markov chain with invariant distribution $\pi$, where $\m B(\Theta)$ is  the Borel-sigma algebra on $\Theta$,  the Dirichlet form $\m{E}:\m L_2(\pi)\times \m L_2(\pi)\to \mb R$ associated with the transition kernel $T$ is given by $\m{E}(g,h)=\frac{1}{2}\int_{x,y\in \Theta^2} (g(x)-h(y))^2T(x,\,\dd y)\pi(\dd x)$.

 
 
   
\section{Additional Results}
\subsection{Additional Simulation}\label{addex}
\begin{figure}[h]
    \centering
    \includegraphics[width=0.46\textwidth]{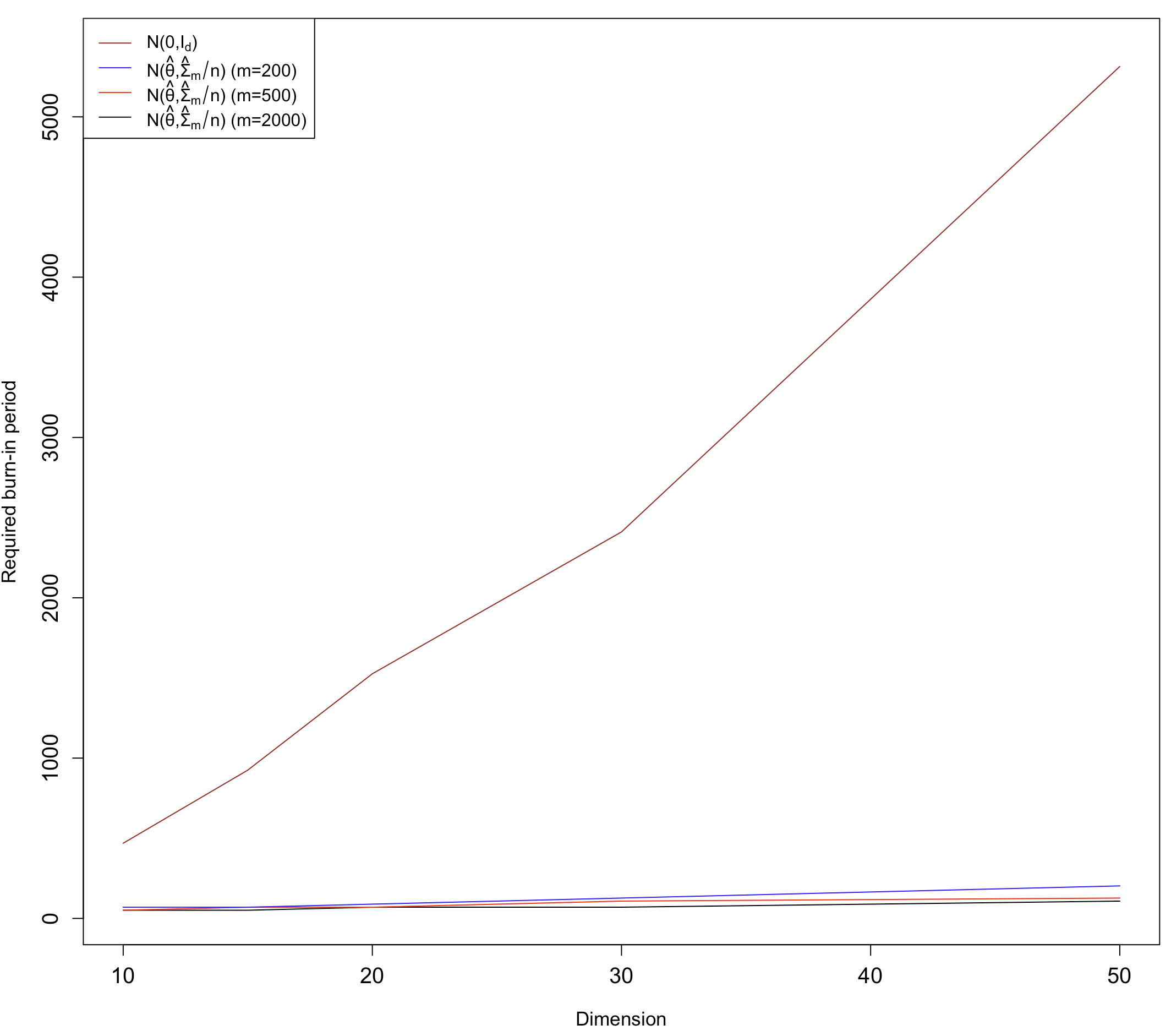}
    \caption{The figure shows the minimal burn-in period required to attain a Gelman-Rubin statistic below $1.1$. It compares two scenarios: MALA with an initial distribution of $\mathcal{N}(0, I_d)$ and a preconditioning matrix of $I_d$, and MALA with an initial distribution of $\mathcal{N}(\wh{\theta}, \wh{\Sigma}_m/n)$ and a preconditioning matrix of $\hat{\Sigma}_m$. We can see the utilization of $\hat{\Sigma}_m$ for constructing the initial distribution and preconditioning matrix can significantly accelerates the convergence of MALA. }
    \label{Fig_burnin}
\end{figure}

In this section, we carry out experiment using Bayesian linear regression with the following posterior
\begin{equation*}
    \pi_n^{\rm mean}(\theta|X^{(n)})\propto \exp\Big(-\frac{1}{2}\sum_{i=1}^n \|Y_i-\wt X_i^T\theta\|^2\Big)\pi(\theta),\quad\theta\in \mb R^d.
\end{equation*}
for exploring the impact of the preconditioning matrix and initial distribution in MALA. We set the sample size $n=2000$ and choose the parameter dimension $d$ from set $\{10,15,20,30,50\}$. The covariates $\wt X$ are generated from a multivariate Gaussian distribution with zero mean and the covariance matrix $\Sigma$ given by a diagonal matrix with elements
\begin{equation*}
\Big(\underbrace{\sqrt{d},\sqrt{d},\cdots,\sqrt{d}}_{\big[\frac{d}{2}\big]},\underbrace{\frac{1}{\sqrt{d}},\frac{1}{\sqrt{d}},\cdots,\frac{1}{\sqrt{d}}}_{d-\big[\frac{d}{2}\big]}\Big).
\end{equation*}
We consider two choices for the preconditioning matrix: one is the inverse (mini-batch) empirical gram matrix $\wh \Sigma_m=(m^{-1}\sum_{j=1}^m \wt X_j\wt X_j^T)^{-1}$, which is an estimator to the covariance matrix  $n^{-1}\Sigma^{-1}$ of the posterior rescaled by $n$. Here, we consider values of $m=\{200,500,2000\}$. The other choice is the standard identity matrix. For the initial distribution, we also consider two options: one is $\m N(\wh\theta,n^{-1}\wh\Sigma_m)$ with $\wh\theta$ being the regression point estimator, as suggested in Corollary~\ref{cor:smoothloss}; and another choice is the standard normal distribution $\m N(0,I_d)$. Figure~\ref{Fig_burnin} displays the minimum number of iteration required for achieving a Gelman-Rubin statistic smaller than $1.1$, which is a common-used rule for determining the burn-in period~\citep{flegal2008markov,roy2020convergence,gelman1992inference}. We observe that choosing the initial distribution as $\m N(\wh\theta,n^{-1}\wh \Sigma_m)$ allows the chain to converge in a very short period, whereas using $\m N(0,I_d)$ requires a much longer time for convergence. Furthermore, we note that the mini-batch size does not significantly affect the required burn-in period, as choosing $m=200$ is sufficient for fast convergence.

\begin{figure}[h]
  
\centering  
\vspace{1em}
\subfigure[Log step size]{
 \includegraphics[trim={0.1cm 0.1cm 0.1cm 0.25cm}, clip,
 width=0.46\textwidth]{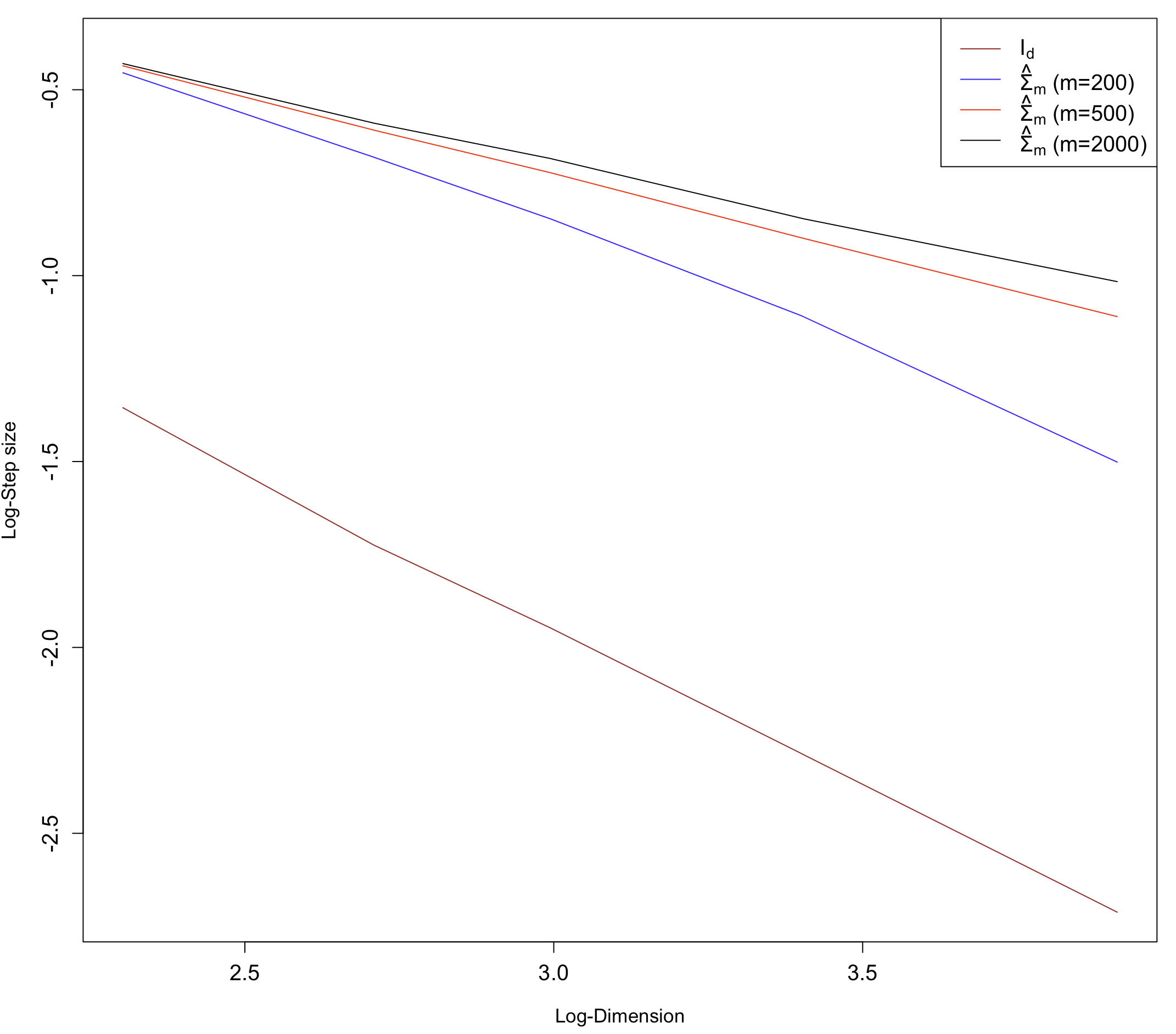}}
\subfigure[Log effective sample size]{
 \includegraphics[trim={0.1cm 0.1cm 0.1cm 0.25cm}, clip, width=0.46\textwidth]{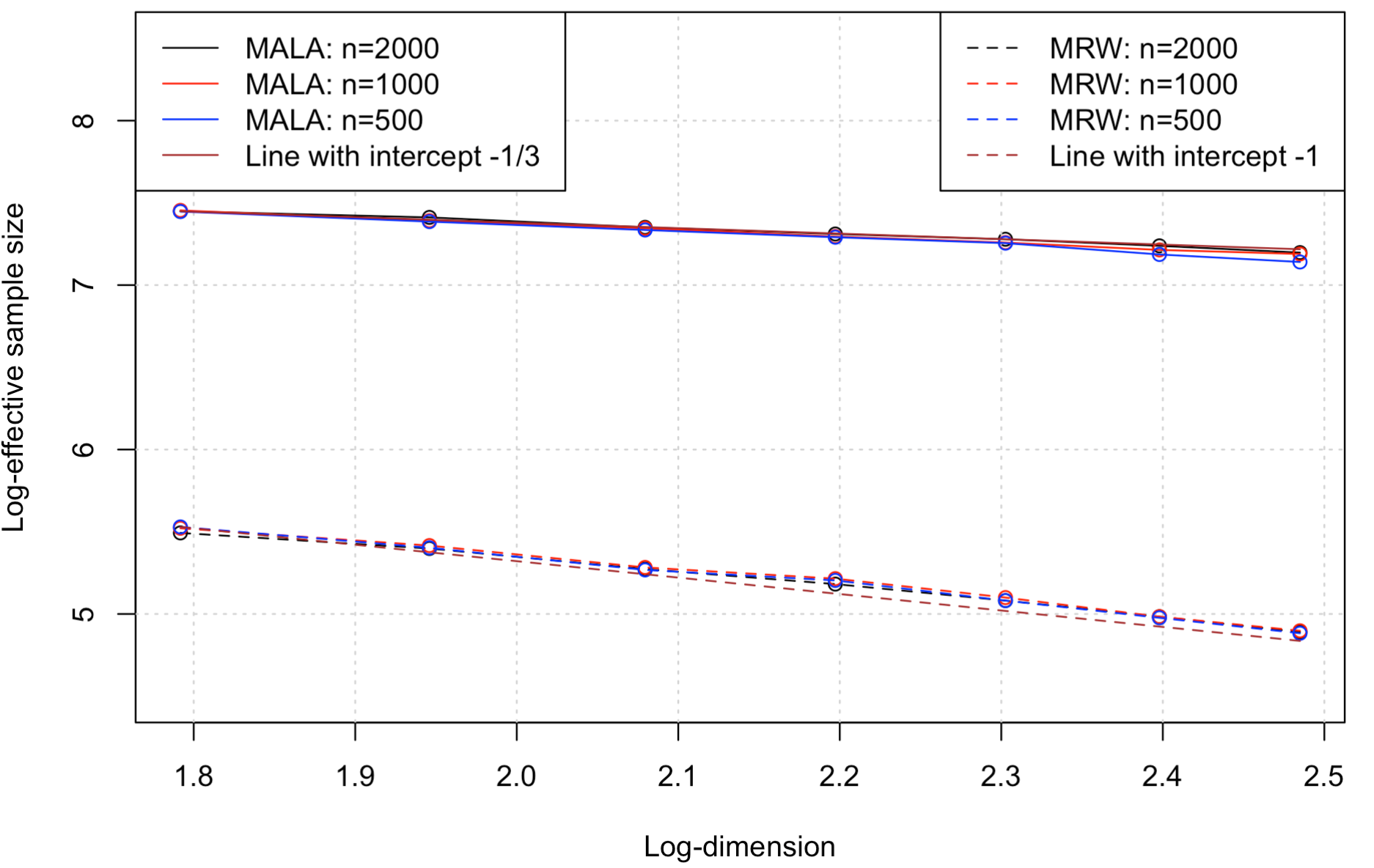}
 \vspace{2em}}

 \caption{Plot (a) illustrates the logarithm of the maximum step size allowed to achieve an average acceptance probability close to 0.57 for various preconditioning matrices and dimensions. Plot (b) illustrates the relationship between the logarithm of the effective sample size and the logarithm of the dimension. The results demonstrate that choosing the preconditioning matrix based on the inverse of the empirical gram matrix enables the use of larger step sizes and leads to a higher effective sample size. Additionally, the disparity between different cases for various values of $m$ becomes more pronounced as the dimension $d$ increases. This is because the approximation error of $\wh{\Sigma}_m$ increases with higher dimensions, necessitating a larger batch size for accurate estimation.}
 
\label{Fig_precondition}
\end{figure}

 Figure~\ref{Fig_precondition} illustrates the largest step size allowed for achieving an average acceptance probability close to $0.57$, as well as the effective sample size, after a total number of $5000$ iterations with a burn-in period of $1000$. 
We observe that utilizing the inverse empirical gram matrix enables a larger step size and leads to a larger effective sample size. Additionally, we find that the best performance is achieved when the batch size $m$ is chosen to be equal to the sample size $n$. This is because a larger batch size provides a better estimator for $\Sigma^{-1}=(\mathbb{E}[\widetilde{X}\widetilde{X}])^{-1}$, resulting in a rescaled covariance matrix $\widehat{\Sigma}_m^{\frac{1}{2}} (\mathbb{E}[\widetilde{X}\widetilde{X}])^{-1}\widehat{\Sigma}_m^{\frac{1}{2}}$ with a smaller condition number. However, when $d\leq 20$, choosing $m=500$ instead of using the full batch does not result in significant loss in performance.

\subsection{Lemmas Related to $s$-conductance Profile}\label{app:profile}

 \begin{lemma}[\bf Mixing time bound via $s$-conductance profile]\label{lemma:mixingtime}
 Consider a reversible,\footnote{A Markov chain with transition kernel $T:\m X\times \m B(\m X)\to \mb R$  and stationary distribution $\mu$ is called reversible if $\mu(\dd x)T(x,\dd y)=\mu(\dd y)T(y,\dd x)$ holds for any $x,y\in \m X$.} irreducible,\footnote{A Markov chain with transition kernel $T:\m X\times \m B(\m X)\to \mb R$ is irreducible if for all $x,y\in \m X$, there is a natural number $k>0$ so that $T^k(x,\dd y)>0$, where $T^k$ is the $k$-step transition  kernel.} $\zeta$-lazy\footnote{A Markov chain is said to be $\zeta$-lazy if at each iteration, the chain is forced to stay at previous iterate with probability $\zeta$. The laziness of Markov chain is also assumed in previous analysis based on $s$-conductance~\citep{https://doi.org/10.1002/rsa.3240040402} and conductance profile~\citep{JMLR:v21:19-441}.}  and smooth Markov chain\footnote{We say that the Markov chain satisfies the smooth chain assumption
if its transition probability function $T$ can be expressed in the form
$T(x,\,\dd y) =\theta(x,y)\,\dd y+\alpha_x\delta_x(\dd y)$ where $\theta$ is a transition density function and $\delta_x$ is the Dirac measure at $x$.} with stationary distribution $\mu$. For any error tolerance $\varepsilon\in(0,1)$, the maximal mixing time in $\chi^2$ divergence of the chain over $M_0$-warm starts can be bounded as
\begin{equation*}
    \tau_{\rm mix}(\varepsilon,M_0)\leq \frac{16}{\zeta}\,\int_{\frac{4}{M_0}}^{\frac{1}{2}} \frac{\dd v}{v\,\Phi_s^2(v)}+ \frac{64}{\zeta}\, \int_{\frac{1}{2}}^{\frac{4\sqrt{2}}{\varepsilon}} \frac{\dd v}{v\,\Phi_s^2(\frac{1}{2})},
\end{equation*}
where $s=\frac{\varepsilon^2}{16M_0^2}$.
\end{lemma}

We can calculate the second term of the upper bound above explicitly as $\frac{64}{\zeta\Phi^2_s(\frac{1}{2})}\log\big(\frac{8\sqrt{2}}{\xi}\big)$.
The next lemma shows that the $s$-conductance profile can be lower bounded given one can: 1.~prove a log-isoperimetric inequality for $\mu$; 2.~bound the total variation distance between $T(x,\cdot)$ and $T(z,\cdot)$ for any two sufficiently close points $x, z$ in a high probability set (not necessarily convex) of $\mu$, which will be referred to as the overlap argument.

\begin{lemma}[\bf $s$-conductance profile lower bound]\label{lemma:conductance}
Consider a Markov chain with Markov transition kernel $T$ and stationary distribution $\mu$.   Given a tolerance $\varepsilon\in (0,1)$ and warming parameter $M_0$, if there are two sets $K$, $E$, and positive numbers $\lambda$, $\psi$, $\omega$ so that
\begin{enumerate}
    \item the probability measure of $\mu$ constrained on $K$, denoted as $\mu|_{K}(\cdot)=\frac{\mu(\cdot\, \cap K)}{\mu(K)}$, satisfies the following log-isoperimetric inequality:
    \begin{equation*}
        \mu|_{K}(S_3)\geq \lambda\cdot t\cdot \min\big\{ \mu|_{K}(S_1),\, \mu|_{K}(S_2)\big\}\cdot \sqrt{\log\Big(1+\frac{1}{ \min\big\{ \mu|_{K}(S_1), \,\mu|_{K}(S_2)\big\}}\Big)}\ ,
    \end{equation*}
    for any partition\footnote{$\bigcup_{j=1}^J A_j$ forms a partition of set $\Omega$ means $\Omega = \bigcup_{j=1}^J A_j$ and $\{A_j\}_{j=1}^J$ are mutually disjoint.}  $K=S_1\cup S_2\cup S_3$ satisfying $\inf_{x\in S_1, z\in S_2}\|x-z\|\geq t$;
    \item for any $x,z\in E$, if $\|x-z\|\leq \psi$,  then $\|T(x,\cdot)-T(z,\cdot)\|_{\rm TV}\leq 1-\omega$;
    \item it holds that $\mu(E)\geq 1-(\lambda\psi\wedge 1)\,\frac{\varepsilon^2}{256M_0^2}$ and  $\mu(K)\geq 1-(\lambda\psi\wedge 1)\,\frac{\varepsilon^2}{256M_0^2}$;
\end{enumerate}
then the $s$-conductance profile $\Phi_s(v)$ with $s=\frac{\varepsilon^2}{16M_0^2}$ can be bounded from below by
\begin{equation*}
  \Phi_s(v)\geq \frac{\omega}{4}\,\min\bigg\{1,\,\frac{\lambda\,\psi}{9} \sqrt{\log\big(1+\frac{1}{v}\big)}\,\bigg\}.
\end{equation*}
\end{lemma}

\smallskip
\noindent By combining this lemma with Lemma~\ref{lemma:mixingtime}, we obtain that if the assumptions in Lemma~\ref{lemma:conductance} hold, then the mixing time of the chain can be bounded as
\begin{equation} 
    \begin{aligned}
         \tau_{\rm mix}(\varepsilon,M_0)&\leq \frac{C_1}{\zeta\omega^2}\, \log M_0+ \frac{C_1}{\zeta\omega^2}\, \lambda^{-2}\psi^{-2} \log(\log M_0) + \frac{C_1}{\zeta\omega^2}\,  \lambda^{-2}\psi^{-2}\log \frac{1}{\varepsilon},
    \end{aligned}
\end{equation}
for some universal constant $C_1$. Therefore, the problem of bounding the mixing time can be converted to verify the assumptions in Lemma~\ref{lemma:conductance}. 
\subsection{Lower Bound of Mixing Time}\label{app:lower}

\begin{theorem}[\bf MALA mixing time lower bound]\label{thmala:lower}
Consider a positive definite preconditioning matrix $\wt I\in \mb R^{d\times d}$, and the target distribution defined as a multivariate normal  $\ov \pi=\m N(0,J^{-1})$, where $J\in \mb R^{d\times d}$ is a covariance matrix with $\wt I^{\frac{1}{2}}J\wt I^{\frac{1}{2}}={\rm diag}(\rho_2,\rho_2,\cdots,\rho_2,\rho_1)$. Assume $1\leq \kappa=\frac{\rho_2}{\rho_1}\leq c_1\cdot d^{c_2}$ for some $c_1,c_2>0$. Then there exists an integer $N$ that depends only on $c_1,c_2$ and universal constants $c_3,c_4$ such that for any $d>N$, $M_0\geq 2$, step size $h>0$ and tolerance $\varepsilon\in (0,1)$,  the $\frac{1}{2}$-lazy version of preconditioned MALA for sampling from $\ov \pi$ has the following  mixing time lower bound in  $\chi^2$ divergence
\begin{equation*}
 \tau_{\rm mix}(\varepsilon,M_0)\geq {c_3\, \kappa \, \Big(\frac{d}{\log (d\kappa)}\Big)^{\frac{1}{3}}} \log \left(\frac{c_4}{\varepsilon}\right).
\end{equation*}
 \end{theorem}

\noindent A proof of Theorem~\ref{thmala:lower} is provided in Appendix~\ref{sec:proofmixing}, part of which is adapted from~\cite{pmlr-v134-chewi21a,wu2022minimax}. Note that the worst-case construction used in~\cite{wu2022minimax} does not satisfy our condition A. As a result, our lower bound has a different dimension dependence of $d^{1/3}$ than that in~\cite{wu2022minimax} of $d^{1/2}$. Additionally, unlike Theorem 1 of \cite{pmlr-v134-chewi21a}, which considers a standard Gaussian target distribution (i.e., $\kappa=1$), our lower bound has an explicit linear dependence on the condition number. From  Theorem~\ref{thmala} and Theorem~\ref{thmala:lower}, we can see that when $\log \big(\frac{M_0\kappa}{\varepsilon}\big)=\m O(d^{\frac{1}{3}})$, our mixing time upper bound and lower bound match up to some logarithmic terms of $(d,\kappa)$ and a double logarithmic term of $M_0$. 
  
 \section{Proof of Main Results}\label{app:main}

\subsection{Proof of Theorem~\ref{thmala} (MALA mixing time upper bound)}\label{Proof1}


\noindent Note that combined with Lemma~\ref{lemma:mixingtime}, if the assumptions in Lemma~\ref{lemma:conductance} holds, we have 
\begin{equation}\label{eqn:mixingfinal}
    \begin{aligned}
         \tau_{\rm mix}(\varepsilon,\mu_0)&\leq \frac{C}{\zeta\omega}\,\int_{\frac{4}{M_0}}^{\frac{1}{2}} \frac{1}{v} \,\dd v +\frac{C}{\zeta\omega}\, \int_{\frac{4}{M_0}}^{\frac{1}{2}} \lambda^{-2}\psi^{-2} \frac{1}{v\log(1+\frac{1}{v})}\,\dd v+ \frac{C}{\zeta\omega}\, \int_{\frac{1}{2}}^{\frac{4\sqrt{2}}{\varepsilon}} \lambda^{-2}\psi^{-2} \frac{1}{v}\,\dd v\\
         &\leq \frac{C_1}{\zeta\omega}\, \log M_0+ \frac{C_1}{\zeta\omega}\, \lambda^{-2}\psi^{-2} \log(\log M_0) + \frac{C_1}{\zeta\omega}\,  \lambda^{-2}\psi^{-2}\log \frac{1}{\varepsilon},
    \end{aligned}
\end{equation}
where the last inequality follows equation (18) of~\cite{JMLR:v21:19-441}.
Now it remains to verify the assumptions in Lemma~\ref{lemma:conductance}.
Fix a lazy parameter $\zeta\in (0,\frac{1}{2}]$. Consider a linear transformation $G:\mb R^d\to \mb R^d$ defined as $G(\theta)= \sqrt{n}\widetilde{I}^{-\frac{1}{2}}(\theta-\wh \theta)$, and let $\widetilde{\mu}_k=G_{\#}\mu_k$ denote the push forward measure of $G$ by $\mu_k$ for $k\in \mb N$ and $\widetilde\pi_{\rm loc}$ denote the push forward measure of $G$ by $\pi_{n}$.  Then it holds that 
 \begin{equation*}
     M_0=\underset{A: \pi_{n}(A)>0} {\sup}\frac{\mu_0(A)}{\pi_{n}(A)}=\underset{A: \widetilde\pi_{\rm loc}(A)>0} {\sup}\frac{\widetilde\mu_0(A)}{\widetilde\pi_{\rm loc}(A)}.
 \end{equation*}
 Moreover, by the invariability of $\chi^2$ measure to linear transformation, we have $\chi^2(\mu_k,\pi_{n}) =\chi^2(\widetilde\mu_{k},\widetilde\pi_{\rm loc})$.  
Define $ \wt Q(\xi,\cdot)$ be the density function of the multivarite normal $ N_d(\xi-h \widetilde{I}^{\frac{1}{2}}\widetilde{\nabla}V_n( \widetilde{I}^{\frac{1}{2}}\xi),2h\,I_d)$, and the corresponding Markov transition kernel
 \begin{equation*}
    \wt T(\xi,\dd y)=\left[1-(1-\zeta)\cdot\int \wt A(\xi,y) \wt Q(\xi, y)\,\dd y\right]\mathbf{1}_{\xi}(\dd y)+(1-\zeta)\cdot \wt Q(\xi,y) \wt A(\xi,y)\,\dd y
 \end{equation*}
 with
 \begin{equation*}
     \wt A(\xi,y)=1\wedge \frac{\wt\pi_{\rm loc}(y)\wt Q(y,\xi)}{\wt \pi_{\rm loc}(\xi)\wt Q(\xi,y)}.
 \end{equation*}
  We have the following lemma.
  \begin{lemma}\label{transform}
     For any $k\in \mb N$, $\widetilde{\mu}_k=G_{\#}\mu_k$ is the probability distribution obtained after $k$ steps of a Markov chain with transition kernel $\wt T$ and initial distribution $\widetilde\mu_0$.
  \end{lemma}
 
  \noindent It remains to calculate the mixing time of $\wt \mu_k$ converging to $\widetilde\pi_{\rm loc}$, which is equivalent to verify the assumptions in Lemma~\ref{lemma:conductance} for Markov transition kernel $  \wt T(\xi,\cdot)$ with stationary distribution $\wt \pi_{\rm loc}$. Recall $K= \{x\,:\, \|\wt I^{-\frac{1}{2}}x\|\leq R\}$. By Condition A, firstly we have 
 \begin{equation*}
 \begin{aligned}
 &\underset{\widetilde\xi\in B_{R}^d}{\sup}\big|V_n( \widetilde{I}^{\frac{1}{2}}\widetilde\xi)-\frac{1}{2}{\widetilde\xi^T \widetilde{I}^{\frac{1}{2}}J \widetilde{I}^{\frac{1}{2}}\widetilde\xi} \big|=\underset{\xi\in K}{\sup}\big|V_n(\xi)-\frac{1}{2}{\xi^TJ\xi} \big|\leq \widetilde{\varepsilon}_0;\\
 &  \underset{\widetilde\xi\in B_{R}^d}{\sup}\big\| \widetilde{I}^{\frac{1}{2}}\widetilde{\nabla}V_n( \widetilde{I}^{\frac{1}{2}}\xi)- \widetilde{I}^{\frac{1}{2}}J \widetilde{I}^{\frac{1}{2}}\xi\|=\underset{\xi\in K}{\sup}\big\| \widetilde{I}^{\frac{1}{2}}(\widetilde{\nabla}V_n(\xi)-J\xi)\|\leq\widetilde{\varepsilon}_1\mnorm{ \widetilde{I}^{\frac{1}{2}}}_{  \rm  op},
       \end{aligned}
 \end{equation*}
 and $\wt \pi_{\rm loc}(\wt \xi\in B_{R/2}^d)=\pi_n(\|\sqrt{n}\wt I^{-\frac{1}{2}}(\theta-\wh \theta)\|\leq R/2)\geq 1-\exp(-4\wt\varepsilon_0)\cdot\frac{h\rho_1\varepsilon^2}{M_0^2}$.  We then verify the log-isoperimetric inequality in the following lemma.
\begin{lemma}\label{lemmalogiso}
 Let $\wt K=B_{R/2}^d$, consider any measurable partition form $\wt K=S_1\cup S_2\cup S_3$ such that $\inf_{x\in S_1,z\in  S_2}\|x-z\|\geq t$,  we have 
  \begin{equation*}
      \wt\pi_{\rm loc}|_{\wt K}(S_3)\geq \frac{\sqrt{\rho_1}}{2}t\exp(-4\widetilde{\varepsilon}_0) \min\{\wt\pi_{\rm loc}|_{\wt K}(S_1),\wt\pi_{\rm loc}|_{\wt K}(S_2)\}\log^{\frac{1}{2}}\Big(1+\frac{1}{\min\{\wt\pi_{\rm loc}|_{\wt K}(S_1),\wt \pi_{\rm loc}|_{\wt K}(S_2)\}}\Big).
  \end{equation*}
  \end{lemma}
 
 \noindent We then show that $\|\wt T(x,\cdot)-\wt T(y,\cdot)\|_{\rm TV}$ can be bounded with high probability in the following lemma.
 
 \begin{lemma}\label{boundTV}
   There exists a set $E$ so that $\wt\pi_{\rm loc}(E)\geq 1-\exp(-4\wt\varepsilon_0)\cdot\frac{2\varepsilon^2 h\rho_1}{M_0^2}$ and for any $x,z\in E$ with 
   $\|x-z\|\leq \frac{\sqrt{h}}{3}$,  we have $\|\wt T(x,\cdot)-\wt T(z,\cdot)\|_{\rm TV}\leq 1-\frac{\exp(-2\wt\varepsilon_0)}{4}$.
 \end{lemma}

Thus the first and second assumptions in Lemma~\ref{lemma:conductance} holds with $\lambda=\frac{\sqrt{\rho_1}}{2}\exp(-4\wt\varepsilon_0)$, $\psi=\frac{\sqrt{h}}{3}$ and $\omega=\frac{\exp(-2\wt\varepsilon_0)}{4}$. Moreover, for the third assumption in Lemma~\ref{lemma:conductance}, by $h\rho_1\leq c_0d^{-\frac{1}{3}}$, for small enough $c_0$, we have 
\begin{equation*}
   \exp(-4\wt\varepsilon_0)\cdot \frac{2\varepsilon^2h\rho_1}{M_0^2}\leq \frac{\sqrt{2h}}{24}\frac{\sqrt{\rho_1}}{2}\exp(-4\wt\varepsilon_0)\frac{\varepsilon^2}{256M_0^2}.
\end{equation*}
 Thus all the assumptions in Lemma~\ref{lemma:conductance} are satisfied.  The desired result then follows from equation~\eqref{eqn:mixingfinal}.
 
 \subsection{Proof of Theorem~\ref{thmala:lower} (MALA mixing time lower bound)}\label{proofthmala:lower}
 Without loss of generality,  we assume $\wt I =I_d$. Otherwise, similar as the proof of Theorem~\ref{thmala}, we could transform the measures $\mu_k$ and $\ov \pi$ by the scale matrix $\wt I^{-\frac{1}{2}}$, and study the convergence of the transformed measures.  We utilize the following lower bound on the $\chi^2$-divergence via Dirichlet form.
 \begin{lemma}\label{lowerbound:spectral}(Corollary 7 of~\cite{wu2022minimax})
    Le $T$ be the transition kernel of  a reversible Markov chain with invariant distribution $\ov\pi$. For any $\varepsilon>0$ and any initial distribution $\mu_0\ll \ov \pi$ satisfying $\chi^2(\mu_0,\ov\pi)<\infty$, let $h_0=\frac{\,\dd\mu_0}{\,\dd \ov\pi}$, if  $\mathcal{E}(h_0, h_0)/\chi^2(\mu_0,\ov\pi)\leq \frac{1}{4}$ with $\mathcal{E}(\cdot,\cdot)$ being the Dirichlet form associated with $T$, then its mixing time in $\chi^2$-divergence has a lower bound 
    \begin{equation*}
        \tau_{\rm mix}(\varepsilon,\mu_0)\geq \frac{1}{4}\left(\frac{\mathcal{E}(h_0, h_0)}{\chi^2(\mu_0,\ov\pi)}\right)^{-1}\log\left(\frac{\chi^2(\mu_0,\ov\pi)}{\varepsilon^2}\right).
    \end{equation*}
 \end{lemma}
\noindent Then, we state the following lemma for bounding  $\mathcal{E}(h_0, h_0)/\chi^2(\mu_0,\ov\pi)$. 
 \begin{lemma}\label{lemma:controlsgap}
      Consider the target distribution $\ov \pi=N_d(0,J^{-1})$ with $J={\rm diag}(\rho_2,\rho_2,\cdots,\rho_2,\rho_1)$ and $1\leq \kappa=\frac{\rho_2}{\rho_1}\leq c_1\cdot d^{c_2}$, then
     \begin{enumerate}
      \item There exists a $2$-warm initial distribution $\mu_0$ with $\chi^2(\mu_0,\ov\pi)\geq \frac{1}{5}$ so that for any $h\in (0,\frac{1}{\rho_1})$,  denote $h_0=\frac{\,\dd\mu_0}{\,\dd \ov\pi}$, then for any $\zeta\in [0,1]$, the term $\mathcal{E}(h_0, h_0)/\chi^2(\mu_0,\ov\pi)$ under the $\zeta$-lazy version MALA transition kernel with step size $h$ satisfies
         \begin{equation*}
             \frac{\mathcal{E}(h_0,h_0)}{\chi^2(\mu_0,\ov\pi)}\leq 60\rho_1h.
         \end{equation*}
         \item When $M_0\geq 2$, there exists an $M_0$-warm initial distribution $\mu_0'$ with $\chi^2(\mu_0',\ov\pi)=M_0-1$ and  a constant $N$ that depends only on $c_1,c_2$ so that when $d\geq N$, denote $h_0=\frac{\dd \mu_0'}{\dd \ov\pi}$, for any $h\in (\frac{8(\log (d\kappa))^{\frac{1}{3}}}{\rho_2d^{\frac{1}{3}}},\infty)$, for any $\zeta\in [0,1]$, the term $\mathcal{E}(h_0, h_0)/\chi^2(\mu_0',\ov\pi)$ under the $\zeta$-lazy version MALA transition kernel with step size $h$ satisfies
         \begin{equation*}
             \frac{\mathcal{E}(h_0,h_0)}{\chi^2(\mu_0',\ov\pi)}\leq  \frac{8}{\kappa d}.
         \end{equation*}
        \end{enumerate}
 \end{lemma}
 
\noindent So when $d\geq N\vee 3$, if $h>\frac{8(\log (d\kappa))^{\frac{1}{3}}}{\rho_2d^{\frac{1}{3}}}$, we have 
\begin{equation*}
    \underset{2-\text{warm }\mu_0}{\sup}\tau_{\rm mix}(\varepsilon,\mu_0)\geq \frac{\kappa d}{46}\log (\frac{1}{\varepsilon^2})\geq   \frac{\kappa d^{\frac{1}{3}}}{46}\log (\frac{1}{\varepsilon^2});
\end{equation*}
when $h\leq \frac{8(\log (d\kappa))^{\frac{1}{3}}}{\rho_2d^{\frac{1}{3}}}$, we have $\rho_1 h<1$ and thus,
\begin{equation*}
    \underset{2-\text{warm }\mu_0}{\sup}\tau_{\rm mix}(\varepsilon,\mu_0)\geq \frac{1}{240}\rho_1^{-1}h^{-1}\log (\frac{1}{5\varepsilon^2}) \geq \frac{\kappa d^{\frac{1}{3}}}{1920(\log (d\kappa))^{\frac{1}{3}}}\log (\frac{1}{5\varepsilon^2}).
\end{equation*}
Proof is completed.
 
 \subsection{Proof of Theorem~\ref{th:Gibbsmixing} (Complexity of MALA for Bayesian sampling)}\label{sec:proofmixing}
Without loss of generality, we can assume the learning rate $\alpha=1$, as otherwise we can take $\ell(X,\theta)=\alpha\cdot \ell(X,\theta)$. We only need to verify that the Assumptions in Theorem~\ref{thmala} holds for the Bayesian Gibbs posterior.
 We state the following Lemmas to verify Condition A. 
 \begin{lemma}\label{th1}
 Let $\kappa_2= {\frac{\beta_1}{\gamma_3+\beta_1(1+\gamma_4)+2\gamma_0-\gamma_4}}\wedge {\frac{1}{1+2\gamma+2\gamma_2+4\gamma_0}}\wedge  {\frac{1}{2+2(\gamma+\gamma_0+\gamma_1)}}$. Under Conditions B.1-B.4, if $d\leq c(\frac{n}{\log n})^{\kappa_2}$ for a small enough constant $c$, then there exist $(n,d)$-independent constants $c_1,C,C_1$ so that it holds with probability at least $1-c_1n^{-2}$
 that for any $\xi\in \mb R^d$ with $1\leq \|\xi\|\leq C\sqrt{n}$,
 \begin{align*}
 &\Big|V_n(\xi)-\frac{\xi^T\m H_{\theta^*}\xi}{2}\Big| \leq  C_1\,\Big(d^{1+\gamma}\|\xi\|  \frac{\log n}{\sqrt{n}} +d^{\gamma_2}\|\xi\|^3 \frac{1}{\sqrt{n}} +d^{\frac{1+\gamma_4}{2}+\gamma_2}\|\xi\|^2\sqrt\frac{\log  n}{n}   \\
 &\qquad\qquad+ d^{\frac{1+\gamma_3}{2}} \|\xi\|^{1+\beta_1}\frac{\sqrt{\log n}}{n^{\beta_1/2}}\Big);\\
  &\Big\|\widetilde{\nabla} V_n(\xi)-\m H_{\theta^*}\xi\Big\|\leq  C_1\,\Big(d^{1+\gamma}  \frac{\log n}{\sqrt{n}} +d^{\gamma_2}\|\xi\|^2 \frac{1}{\sqrt{n}} +d^{\frac{1+\gamma_4}{2}+\gamma_2}\|\xi\|\sqrt\frac{\log  n}{n}  \\
  &+ d^{\frac{1+\gamma_3}{2}} \|\xi\|^{\beta_1}\frac{\sqrt{\log n}}{n^{\beta_1/2}}\Big)\text{ with }\widetilde{\nabla}V_n(\xi)= 
       \frac{1}{\sqrt n}\sum_{i=1}^n g\Big(X_i,\frac{\xi}{\sqrt n}+\widehat\theta\,\Big)-
       \frac{1}{\sqrt n}\nabla(\log \pi)\Big(\frac{\xi}{\sqrt n}+\widehat\theta\,\Big).
     \end{align*}
 \end{lemma}
 We provide in the following lemma a tail inequality for the Gibbs posterior $\pi_n$.
 \begin{lemma}\label{lemmatail}
 Under Condition B.1-B.4. when $d\leq c\frac{n^{\kappa_3}}{\log n}$ for a small enough constant $c$, where
  \begin{equation*}
 \begin{aligned}
      \kappa_3&= {\frac{\beta_1}{1+\gamma_3+[(2\gamma_0)\vee ((1+\gamma_0)(1+\beta_1))]}}\wedge {\frac{1}{3+\gamma_0+((2\gamma)\vee (\gamma_4+2\gamma_2+\gamma_0)\vee(2\gamma_2+2\gamma_0))}}\\
      &\wedge  {\frac{1}{1+2\gamma+6\gamma_0+4\gamma_2+\gamma_4}}\wedge {\frac{1}{2\gamma+2\gamma_0+2\gamma_1+(2\vee (1+\gamma_4))}},
       \end{aligned}
 \end{equation*}
  then there exist $(n,d)$-independent constants $c_1,c_2,c_3$ so that it holds with probability at least $1-c_1n^{-2}$ that 
 \begin{equation*}
     \pi_n\Big(\sqrt{n}\|\wt{I}^{-\frac{1}{2}}(\theta-\wh\theta)\|\geq \mnorm{\wt I^{-\frac{1}{2}}}_{\rm op}\vee\frac{3(\sqrt{d}+t)}{\sqrt{\lambda_{\min}(\wt J)}}\Big)\leq \exp(-t^2)+c_2\exp\big(-c_3 n\,d^{-\gamma_0}(d^{-\gamma_1}\wedge d^{-2\gamma_0-2\gamma_2})\big),
 \end{equation*}
 where $\wt J=\wt{I}^{\frac{1}{2}}\m H_{\theta^*}\wt{I}^{\frac{1}{2}}$.
 \end{lemma}
 By Condition B.1, we have 
 $\mnorm{\m H_{\theta^*}}_{\rm op}\leq C\,d$  and $\mnorm{\m H_{\theta^*}^{-1}}_{\rm op}\leq C\,d^{\gamma_0}$. Moreover, since $\mnorm{\wt I^{-1}}_{\rm op}\mnorm{\wt I}_{\rm op}\leq C \mnorm{\m H_{\theta^*}}_{\rm op}\mnorm{\m H_{\theta^*}^{-1}}_{\rm op}$ and
 ${\mnorm{\wt I}_{\rm op}\mnorm{(\wt I^{\frac{1}{2}}H_{\theta^*}\wt I^{\frac{1}{2}})^{-1}}_{\rm op}}\leq C\, \mnorm{\m H_{\theta^*}^{-1}}_{\rm op}$, we can obtain that there exists constants $C_2,C_3$ so that for any $R=\mnorm{\wt I^{-\frac{1}{2}}}_{\rm op}\vee \frac{3(\sqrt{d}+t)}{\sqrt{\lambda_{\min}(\wt J)}}$ with $t\geq 0$,  and set $K=\{x:\|\wt I^{-1/2}x\|\leq R\}$, we have $K\subseteq \{x\,:\, \|x\|\leq  C_2 d^{\frac{1+\gamma_0}{2}}+C_3 td^{\frac{\gamma_0}{2}}\}$. Then by Lemma~\ref{th1}, for any $t=C_1\, (d^{\frac{  \gamma_5}{2}}+\sqrt{\log n})$ (note that $\gamma_5\geq 1$),  we can find a constant $c$ so that when $d\leq c\frac{n^{\kappa_1}}{\log n}$, we have 
 \begin{equation*}
  \mnorm{\wt I}_{\rm op}\,R^2  \underset{\xi \in K}{\sup} \|\widetilde{\nabla} V_n(\xi)-\m H_{\theta^*}\xi\Big\|^2\leq d^{\frac{1}{3}}.
 \end{equation*}
 So in this case the step size parameter $\wt h=h/n$ in Theorem~\ref{thmala} satisfies 
 \begin{equation*}
     h\leq c_0\cdot\bigg[\rho_2\Big(2d^{\frac{1}{3}}+d^{\frac{1}{4}}\big(\log \frac{M_0d\kappa}{\varepsilon}\big)^{\frac{1}{4}}+\big(\log \frac{M_0d\kappa}{\varepsilon}\big)^{\frac{1}{2}}\Big)\bigg]^{-1}
 \end{equation*}
 Then by the assumption $\log (\frac{M_0\rho_2}{\varepsilon\rho_1})\leq C_1\,(d^{\gamma_5}+\log n)$, using Lemma~\ref{lemmatail} and $d\leq c\frac{n^{\kappa_1}}{\log n}$, we can obtain that there exists a large enough $C_1$ so that for $t=C_1\, (d^{\frac{  \gamma_5}{2}}+\sqrt{\log n})$, 
  $\pi_n(K)=\pi_n\Big(\sqrt{n}\|\wt{I}^{-\frac{1}{2}}(\theta-\wh\theta)\|\geq \mnorm{\wt I^{-\frac{1}{2}}}_{\rm op}\vee\frac{3(\sqrt{d}+t)}{\sqrt{\lambda_{\min}(\wt J)}}\Big)\geq 1-\frac{h\rho_1\varepsilon^2}{M_0^2}$.
 So the Assumptions in Theorem~\ref{thmala} are satisfied.

 \section{Proof of Lemmas for Theorem~\ref{thmala} and Theorem~\ref{thmala:lower}}
 \subsection{Proof of Lemma~\ref{lemma:mixingtime}}
  Fix an arbitrary $\varepsilon>0$. Suppose $\tau_{\rm mix}(\sqrt{2}\varepsilon,\mu_0)>N=\int_{\frac{4}{M_0}}^{\frac{1}{2}} \frac{16\dd v}{\zeta\cdot v\Phi_s^2(v)}+\int_{\frac{1}{2}}^{\frac{4}{\varepsilon}} \frac{64\dd v}{\zeta\cdot v\Phi_s^2(\frac{1}{2})}$. Then for any $k\leq N$, $\chi^2(\mu_k,\mu)> 2\varepsilon^2$, where we use $\chi^2(\cdot,\cdot)$ to denote the $\chi^2$ divergence, $\mu_k$ to denote the distribution in $k$ step of the Markov chain and $\mu\in \m P(\mb R^d)$ to denote the stationary distribution. Then we will prove by contradiction that if $N<\tau_{\rm mix}(\sqrt{2}\varepsilon,\mu_0)$, then when $k=N$,  $\chi^2(\mu_k,\mu)\leq 2\varepsilon^2$, which oppositely implies $N\geq \tau_{\rm mix}(\sqrt{2}\varepsilon,\mu_0)$.
  Our proof is based on the strategy used in~\cite{JMLR:v21:19-441}.   We first introduce the following  related notations. For a measurable set $S\subseteq \mb R^d$ and positive numbers $\varepsilon, M_0$ , the $(\varepsilon,M_0)$-spectral gap for the set $S$ is defined as 
 \begin{equation*}
     \Lambda_{\varepsilon,M_0}(S):= \inf _{g \in c_{\varepsilon,M_0}^{+}(S)} \frac{\mathcal{E}(g, g)}{\operatorname{Var}_{\mu}(g)}
 \end{equation*}
 where 
 \begin{equation*}
     c_{\varepsilon,M_0}^{+}(S):=\left\{g \in L_{2}\left(\mu\right) \mid \operatorname{supp}(g)=\{x \,:\, g(x)>0 \} \subset S,\, 0\leq g \leq M_0,\, {\rm Var}_{\mu}(g)\geq \varepsilon^2\right\},
 \end{equation*}
 and
 \begin{equation*}
     \mathcal{E}(g, g)=\frac{1}{2}\int (g(x)-g(y))^2 T(x,\dd y) \mu(\dd x),
 \end{equation*}
 with $T(x,\dd y)$ denoting the Markov transition kernel.
 Moreover, we can define the $(\varepsilon,M_0,s)$-spectral profile $\overline{\Lambda}_s^{\varepsilon,M_0}$ as 
 \begin{equation*}
     \overline{\Lambda}_s^{\varepsilon,M_0}(v):=\underset{\mu(S)\in (s,v]}{\inf} \Lambda_{\varepsilon,M_0}(S).
 \end{equation*}
 Define the ratio density 
 \begin{equation*}
     h_k(x)=\frac{\mu_k(x)}{\mu(x)}.
 \end{equation*}
 Note that 
 \begin{equation*}
     \mb E_{\mu}[h_k]=1 \quad \text{and} \quad \chi^2(\mu_k,\mu)={\rm Var}_{\mu}(h_k),
 \end{equation*}
 and $h_{k}(x)\leq M_0$ for all $k\geq 0$ (see for example, equation (64) of~\cite{JMLR:v21:19-441}).
 By tracking the proof of Lemma 11 in~\cite{JMLR:v21:19-441}, it suffices to show that for any $k \leq N$, 
 \begin{equation}\label{Dirichletform}
     2\,\m{E}(h_k,h_k)\geq  {\rm Var}_{\mu}(h_k)\overline{\Lambda}^{\varepsilon,M_0}_s\Big(\frac{4}{{\rm Var}_{\mu}(h_k)}\Big),
 \end{equation}
 and 
 \begin{equation}\label{profilecond}
 \overline\Lambda_{s}^{\varepsilon,M_0}(v) \geq\left\{\begin{array}{ll}\frac{\Phi_{s}^{2}(v)}{16} & \text { for all } v \in\left[\frac{4}{M_0}, \frac{1}{2}\right]; \\ \frac{\Phi_{s}^{2}\left( \frac{1}{2}\right)}{64} & \text { for all } v \in (\frac{1}{2}, \infty),\end{array}\right.
  \end{equation}
  with $s=\frac{\varepsilon^2}{16M_0^2}$. We first prove claim~\eqref{Dirichletform}. Define $\gamma_k= \frac{{\rm Var}_{\mu}(h_k)}{4\mb E_{\mu}[h_k]}=\frac{{\rm Var}_{\mu}(h_k)}{4}$. Then for any $k\leq N$,
 \begin{equation*}
 \begin{aligned}
     \operatorname{Var}_{\mu}\left((h_k-\gamma_k)_{+}  \right)&=\mb{E}_{\mu}\left[((h_k-\gamma_k)_{+})^2 \right]-(\mb{E}_{\mu}\left[(h_k-\gamma_k)_{+} \right])^2\\
     &\overset{(i)}{\geq}\mb{E}_{\mu}[h_k^2]-2\gamma_k\mb{E}_{\mu}[h_k]-(\mb{E}_{\mu}[h_k])^2\\
     &=  \operatorname{Var}_{\mu}(h_k)-2\gamma_k\mb{E}_{\mu}[h_k]\\
     &= \frac{1}{2}\operatorname{Var}_{\mu}(h_k)\geq \varepsilon^2,
      \end{aligned}
 \end{equation*}
  where $(x)_{+}=\max\{0,x\}$, $(i)$ is due to $((a-b)_{+})^2\leq a^2-2ab$, $(a-b)_{+}\leq a$, and the last inequality is due to the assumption that $N<\tau_{\rm mix}(\sqrt{2}\varepsilon,\mu_0)$; moreover, since for any $x\in \mb R^d$,  $0\leq h_k(x)\leq M_0$, we can get $(h_k-\gamma_k)_{+}\in c_{\varepsilon,M_0 }^{+}(\{h_k>\gamma_k\})$, which leads to 
 \begin{equation}\label{eqn:Dirichlet1}
     \begin{aligned}
         \m{E}(h_k,h_k)\overset{(ii)}{\geq} \m{E}((h_k-\gamma_k)_{+},(h_k-\gamma_k)_{+})\geq \operatorname{Var}_{\mu}\left((h_k-\gamma_k)_{+}  \right) \cdot\inf _{f \in c_{\varepsilon,M_0 }^{+}(\{h_k>\gamma_k\})} \frac{\mathcal{E}(f, f)}{\operatorname{Var}_{\mu}(f)},
     \end{aligned}
 \end{equation}
where $(ii)$ follows from the fact that $(a-b)^2=(a-c-(b-c))^2\geq ((a-c)_{+}-(b-c)_{+})^2$. Furthermore, We have for any $k\leq N$,
 \begin{equation*}
   M_0^2 \,\mu(h_k\geq \gamma_k)\geq \mb{E}_{\mu}[((h_k-\gamma_k)_{+})^2]\geq  \operatorname{Var}_{\mu}\left((h_k-\gamma_k)_{+}  \right)\geq \frac{1}{2}\operatorname{Var}_{\mu}(h_k)\geq \varepsilon^2.
 \end{equation*}
   On the other hand, by applying Markov's inequality, we also have 
 \begin{equation*}
     \mu(h_k\geq \gamma_k)\leq \frac{\mb{E}_{\mu}[h_k]}{\gamma_k}=\frac{4}{{\rm Var}_{\mu}(h_k)}.
 \end{equation*}
 Thus by equation~\eqref{eqn:Dirichlet1}, we can get for $s=\frac{\varepsilon^2}{16M_0^2}$,
 \begin{equation*}
      \m{E}(h_k,h_k)\geq  \frac{1}{2}\operatorname{Var}_{\mu}(h_k) \overline\Lambda_{s}^{\varepsilon,M_0}\Big(\frac{4}{{\rm Var}_{\mu}(h_k)}\Big).
 \end{equation*}
 Then we prove claim~\eqref{profilecond}. 
 For $v\in \big[\frac{4}{M_0},\frac{1}{2}\big]$, fix any $A\subset \mb R^d$ with $s<\mu(A)\leq v$  and $g\in c^{+}_{\frac{\varepsilon}{2},M_0}(A)$. Then by 
 \begin{equation*}
 \begin{aligned}
     &\mathbb{E}_{\mu} \Big[\int(g^2(x)-g^2(y))_{+} T(x, \dd y)\Big]\\
     &= \mathbb{E}_{\mu} \Big[\int (g^2(x)-g^2(y))\mathbf{1}(g^2(x)>g^2(y)) T(x, \dd y)\Big]\\
      &= \mathbb{E}_{\mu} \Big[\int\int_{0}^{+\infty} \mathbf{1}(g^2(y)\leq t<g^2(x)) \,\dd t \,T(x, \dd y) \Big]\\
     &= \int_{0}^{+\infty}\mathbb{E}_{\mu}\Big[\int \mathbf{1}(g^2(y)\leq t<g^2(x)) T(x,\dd y)\Big]\,\dd t,
      \end{aligned}
 \end{equation*}
 let $H_t=\{x\in \mb R^d\,: \,g^2(x)>t\}$, we have
 \begin{equation*}
     \begin{aligned}
     &\int\int |g^2(x)-g^2(y)| T(x,\dd y) \mu(\dd x)\\
     &\geq \int\int (g^2(x)-g^2(y))_{+}T(x,\dd y) \mu(\dd x) \\
     &=\int_{0}^{+\infty}\mathbb{E}_{\mu}\Big[\int \mathbf{1}(g^2(y)\leq t<g^2(x)) T(x,\dd y)\Big]\,\dd t\\
     &=\int_{0}^{+\infty}\int_{x\in H_t}T(x,H_t^c) \mu(\dd x) \,\dd t.
     \end{aligned}
 \end{equation*}
 Let $t^*=\sup\{t\geq 0:  \mu(H_t)> s\}$, note that $t^*$ always exists as otherwise, $ \mu(g(x)=0)\geq 1-s$ and thus ${\rm Var}_{ \mu}(g)\leq M_0^2s= \frac{\varepsilon^2}{16}$, which is contradictory to the requirement that ${\rm Var}_{ \mu}(g)\geq \frac{\varepsilon^2}{4}$. Then 
 \begin{equation*}
     \begin{aligned}
     &\int\int |g^2(x)-g^2(y)| T(x,\dd y) \mu(\dd x)\\
     &\geq \int_{0}^{t^*}\int_{x\in H_t}T(x,H_t^c)  \mu(\dd x)\,\dd t+\int_{t^*}^{+\infty}\int_{x\in H_t}T(x,H_t^c)  \mu(\dd x)\,\dd t\\
     &\geq \int_{0}^{t^*} ( \mu(H_t)-s)
     \,\dd t\cdot \Phi_s( \mu(A)) \\
     &= \Big(\mb E_{ \mu}[g^2]-\int_{t^*}^{M_0^2} \mu(H_t)\,\dd t-st^*\Big)\cdot \Phi_s( \mu(A)) \\
     &\overset{(ii)}{\geq} \Big(\mb E_{ \mu}[g^2]- \frac{\varepsilon^2}{8}\Big)\cdot \Phi_s( \mu(A))\\
     &\overset{(iii)}{\geq}\frac{1}{2}\mb E_{ \mu}[g^2] \Phi_s( \mu(A)),
     \end{aligned}
 \end{equation*}
 where $(ii)$ uses the fact that $t^*\leq M_0^2$ and when $t>t^*$, $ \mu(H_t)\leq s=\frac{\varepsilon^2}{16M_0^2}$ and $(iii)$ uses $\mb E_{ \mu}[g^2]\geq {\rm Var}_{ \mu}(g)\geq \frac{\varepsilon^2}{4}$.  Moreover, since
 \begin{equation*}
     \begin{aligned}
         &\int\int |g^2(x)-g^2(y)| T(x,\dd y) \mu(\dd x)\\
         &\leq \sqrt{\int\int (g(x)-g(y))^2 T(x,\dd y) \mu(\dd x) }\cdot \sqrt{\int\int (g(x)+g(y))^2 T(x,\dd y) \mu(\dd x)}\\
         &\leq \sqrt{2\m{E}(g,g)}\cdot  \sqrt{\int\int (2g^2(x)+2g^2(y)) T(x,\dd y) \mu(\dd x)}\\
         &= \sqrt{2\,\m{E}(g,g)}\cdot\sqrt{4\,\mb{E}_{ \mu}[g^2]},
     \end{aligned}
 \end{equation*}
 we have 
 \begin{equation*}
     \begin{aligned}
         &\qquad\frac{1}{2}\mb E_{ \mu}[g^2]\cdot \Phi_s( \mu(A))\leq \sqrt{2\,\m{E}(g,g)}\cdot\sqrt{4\,\mb{E}_{ \mu}[g^2]}\\
         &\Rightarrow \frac{\m{E}(g,g)}{\operatorname{Var}_{ \mu}(g)}\geq \frac{\Phi^2_s( \mu(A))}{16} .
     \end{aligned}
 \end{equation*}
 Taking infimum over  $A\subset \mb R^d$ with $s<\mu(A)\leq v$  and $g\in c^{+}_{\frac{\varepsilon}{2},M_0}(A)$, we have
 \begin{equation}\label{eqn:spectralbound}
     \overline{\Lambda}_s^{\varepsilon,M_0}(v)\geq  \ \overline{\Lambda}_s^{\frac{\varepsilon}{2},M_0}(v)\geq\underset{s<\mu(A)\leq v}{\inf} \frac{ \Phi_s^2(\mu(A))}{16}\geq \frac{\Phi_s^2(v)}{16}.
 \end{equation}
 For the case $v>\frac{1}{2}$, consider any $A\subset \mb R^d$ with $\mu(A)>\frac{1}{2}$ and $g\in c^{+}_{\varepsilon,M_0}(A)$. Let $0\leq\gamma\leq M_0$ be the number such that 
 \begin{equation*}
     s<\mu(\{g>\gamma\})\vee\mu(\{g<\gamma\})\leq\frac{1}{2}.
 \end{equation*}
   $\gamma$ always exists as otherwise, there exists $0\leq \widetilde \gamma\leq M_0$ such that $\mu\{g=\widetilde\gamma\}\geq 1-2s$, which leads to ${\rm Var}_{\mu}(g)\leq \mb E_{\mu}[(g-\widetilde\gamma)^2]\leq 2M_0^2s<\varepsilon^2$, and this causes contradiction. We first consider the case that $\mu(\{g>\gamma\})\wedge\mu(\{g<\gamma\})>s$.  We have
      \begin{equation*}
       \begin{aligned}
           \m{E}(g,g)= \m{E}((g-\gamma),(g-\gamma))\geq  \m{E}((g-\gamma)_{+},(g-\gamma)_{+})+ \m{E}((g-\gamma)_{-},(g-\gamma)_{-}).\\
       \end{aligned}
   \end{equation*}
 Since for any function $h\geq 0$ with $\mu({\rm supp}(h))\leq \frac{1}{2}$, using Cauchy-Schwarz inequality, it holds that 
 \begin{equation*}
     \mb E_{\mu}[h^2]=\int_{x\in \rm supp(h)} h^2(x) \mu(x) dx\geq \frac{(\mb E_{\mu}[h])^2}{\mu({\rm supp}(h))}\geq 2\,\mb (E_{\mu}[h])^2,
 \end{equation*}
 which leads to 
 \begin{equation*}
     {\rm Var}_{\mu}(h)\geq \frac{1}{2}\mb E_{\mu}[h^2].
 \end{equation*}
 Since 
 $\varepsilon^2\leq {\rm Var}_{\mu}(g)\leq \mb {E}_{\mu}[(g-\gamma)^2]$ and $\mb{E}_{\mu}[(g-\gamma)^2]=\mb {E}_{\mu}[(g-\gamma)_{+}^2]+\mb {E}_{\mu}[(g-\gamma)_{-}^2]$, w.l.o.g, we can assume $  \mb {E}_{\mu}[(g-\gamma)_{+}^2]\geq \frac{\mb{E}_{\mu}[(g-\gamma)^2]}{2}\geq  \frac{\varepsilon^2}{2}$. Then taking $h=(g-\gamma)_{+}$, we can obtain 
 \begin{equation}\label{eqn:spectralbound2}
 \begin{aligned}
     \m{E}(g,g)&\geq  \m{E}((g-\gamma)_{+},(g-\gamma)_{+})\\
     &\geq \mb{E}_{\mu}[(g-\gamma)^2_{+}]\cdot \frac{ \m{E}((g-\gamma)_{+},(g-\gamma)_{+})}{2\,{\rm Var}_{\mu}((g-\gamma)_{+})}\\
     &\overset{(i)}{\geq} \frac{1}{4}{\rm Var}_{\mu}(g)\cdot \underset{\mu(S)\in (s,\frac{1}{2}]}{\inf}\underset{f\in c^{+}_{\frac{\varepsilon}{2},M_0}(S)}{\inf} \frac{\m{E}(f,f)}{{\rm Var}_{\mu}(f)}\\
     &\overset{(ii)}{\geq} \frac{1}{64} {\rm Var}_{\mu}(g)\Phi_s^2(\frac{1}{2}),
      \end{aligned}
 \end{equation}
 where $(i)$ uses $\mb {E}_{\mu}[(g-\gamma)_{+}^2]\geq \frac{\mb{E}_{\mu}[(g-\gamma)^2]}{2}\geq \frac{{\rm Var}_{\mu}(g)}{2}$ and ${\rm Var}_{\mu}((g-\gamma)_{+})\geq \frac{1}{2}\mb{E}_{\mu}[(g-\gamma)^2_{+}]\geq \frac{\varepsilon^2}{4}$, and $(ii)$ uses~\eqref{eqn:spectralbound}. Then we consider the case that $\mu(\{g>\gamma\})\wedge\mu(\{g<\gamma\})\leq s< \mu(\{g>\gamma\})\vee\mu(\{g<\gamma\})$. W.l.o.g, we can assume $\mu(\{g>\gamma\})>s$.  Then we can obtain 
 \begin{equation*}
 \begin{aligned}
     \mb {E}_{\mu}[(g-\gamma)_{+}^2]&= \mb {E}_{\mu}[(g-\gamma)^2]- \mb {E}_{\mu}[(g-\gamma)_{-}^2]\\
     &\geq \mb {E}_{\mu}[(g-\gamma)^2]-M_0^2s=\mb {E}_{\mu}[(g-\gamma)^2]-\frac{\varepsilon^2}{8}\geq \frac{\mb {E}_{\mu}[(g-\gamma)^2]}{2}, \end{aligned}
 \end{equation*}
 where the last inequality is due to $\mb {E}_{\mu}[(g-\gamma)^2]\geq {\rm Var}_{\mu}(g)\geq \varepsilon^2$. We can then obtain the desired result by taking infimum over $A\subset \mb R^d$ with $\mu(A)>\frac{1}{2}$ and $g\in c^{+}_{\varepsilon,M_0}(A)$ in~\eqref{eqn:spectralbound2}.

\subsection{Proof of Lemma~\ref{lemma:conductance}}
 The proof follows from the standard conductance argument in~\cite{pmlr-v134-chewi21a,10.2307/30243694,JMLR:v20:19-306,JMLR:v21:19-441}. Let $s= \frac{\varepsilon^2}{16M_0^2}$,  and let $S$ be any measurable set of $\mb R^d$ with $s\leq \mu(S)\leq v\leq \frac{1}{2}$. Define the following subsets:
 \begin{equation*}
     \begin{aligned}
     & S_1:=\{x\in S|T(x,S^c)\leq \frac{\omega}{2}\},\\
      & S_2:=\{x\in S^c|T(x,S)\leq \frac{\omega}{2}\},\\
      &  S_3:= (S_1\cup S_2)^c,
     \end{aligned}
 \end{equation*}
  Then same as the analysis in~\cite{pmlr-v134-chewi21a}, if $\mu(S_1)\leq \mu(S)/2$ or $\mu(S_2)<\mu(S^c)/2$, then by the fact that $\mu$ is stationary w.r.t the transition kernel $T$, we have 
  \begin{equation*}
  \begin{aligned}
      \int_S T(x,S^c)\mu(\dd x)&=\int T(x,S)\mu(\dd x)-\int_{S}T(x,S)\mu(\dd x)\\
      &=\int_{S^c} T(x,S)\mu(\dd x)\geq\frac{\omega}{2}\cdot \max\{\mu(S\cap S_1^c),\mu(S^c\cap S_2^c)\}\\
      &\geq \frac{\omega\cdot \mu(S)}{4}.
        \end{aligned}
  \end{equation*}
  Then when $\mu(S_1)\wedge \mu(S_2)\geq \frac{\mu(S)}{2}$, consider $x\in E\cap S_1$ and $z\in E\cap S_2$, then $\|T_x-T_z\|_{  \rm TV}\geq T(z,S^c)-T(x,S^c)\geq 1-\omega$, thus $\|x-z\|\geq \psi$, which implies that $\inf_{x\in E\cap S_1,z\in E \cap S_2}\|x-z\|\geq  \psi$.    Then consider sets $E\cap K\cap S_1$ and $E\cap K\cap S_2$ in the log-isoperimetric inequality  of $\mu|_{K}$, we can obtain that
  \begin{equation*}
  \begin{aligned}
      \mu|_{K}(((E\cap K\cap S_1)\cup (E\cap K\cap S_2))^c)&\geq \lambda\cdot\psi\cdot\min\{\mu|_K(E\cap K\cap S_1),\mu|_K(E\cap K\cap S_2)\}\\
      &\cdot \log^{\frac{1}{2}}\Big(1+\frac{1}{\min\{\mu|_K(E\cap K\cap S_1),\mu|_K(E \cap K\cap S_2)\}}\Big)\\
      &\geq\lambda\cdot\psi\cdot\min\{\mu(E\cap K\cap S_1),\mu(E\cap K\cap S_2)\}\\
      &\cdot \log^{\frac{1}{2}}\Big(1+\frac{1}{\min\{\mu(E\cap K\cap S_1),\mu(E \cap K\cap S_2)\}}\Big),\\
        \end{aligned}
  \end{equation*}
  where the last inequality is due to the fact that the function $x\log^{\frac{1}{2}} (1+\frac{1}{x})$ is an increasing function.  W.l.o.g, we can assume $\mu(E\cap K\cap S_1)\leq\mu(E\cap K\cap S_2)$, then by $((E\cap K\cap S_1)\cup (E\cap K\cap  S_2))^c\subseteq E^c\cup K^c \cup S_3$ and $\mu(E^c)\leq  (\lambda\psi\wedge 1)\frac{\varepsilon^2}{256 M_0^2}=\frac{(\lambda\psi\wedge 1)s}{16}$, $\mu(K^c)\leq  \frac{(\lambda\psi\wedge 1)s}{16}$, we can obtain
  \begin{equation*}
      \begin{aligned}
      &\mu(S_3)+\frac{16\lambda\psi s}{127}\\
      &\geq \mu(S_3)+ \frac{\mu(K^c)+\mu(E^c)}{\mu(K)}\\
      &\geq  \frac{\mu(S_3)+\mu(E^c)}{\mu(K)}\\
      &\geq   \mu|_{K}(((E\cap K\cap S_1)\cup (E\cap K\cap S_2))^c)\\
      &\geq \lambda\cdot\psi\cdot\mu(E\cap K\cap S_1) \cdot \log^{\frac{1}{2}}\Big(1+\frac{1}{\mu(E\cap K\cap S_1)}\Big),\\
      &\overset{\rm {(i)}}{\geq} \lambda\cdot\psi\cdot\big(\frac{\mu(S)}{4}+\frac{s}{4}-\frac{s}{8}\big)
      \log^{\frac{1}{2}}\Big(1+\frac{1}{\frac{\mu(S)}{4}+\frac{s}{4}-\frac{s}{8}\big)}\Big)\\
      &\geq \lambda\cdot\psi\cdot\frac{\mu(S)}{4} 
      \log^{\frac{1}{2}}\Big(1+\frac{4}{\mu(S)}\Big),
      \end{aligned}
  \end{equation*}
  where (i) uses $ \mu(E\cap K\cap S_1)\geq \mu(S_1)-\mu(E^c)-\mu(K^c)$, $\mu(S_1)\geq \frac{\mu(S)}{2}\geq \frac{s}{2}$ and the function $x\log^{\frac{1}{2}} (1+\frac{1}{x})$ is an increasing function. Then by $\mu(S)\geq s$, we can obtain 
  \begin{equation*}
      \mu(S_3)\geq \lambda\cdot\psi\cdot\frac{\mu(S)}{9} 
      \log^{\frac{1}{2}}\Big(1+\frac{4}{\mu(S)}\Big),
  \end{equation*}
  hence 
  \begin{equation*}
  \begin{aligned}
      \int_S T(x,S^c) \mu(\dd x)&\geq \frac{1}{2}\left(  \int_S T(x,S^c) \mu(\dd x)+  \int_{S^c} T(x,S) \mu(\dd x)\right)\\
      &\geq \frac{\omega}{4}\mu(S_3)\geq \frac{\omega\cdot\lambda\cdot\psi}{36} \cdot{\mu(S)} 
      \log^{\frac{1}{2}}\Big(1+\frac{4}{\mu(S)}\Big),
       \end{aligned}
  \end{equation*}
  which leads to 
  \begin{equation*}
  \begin{aligned}
     \frac{\int_S T(x,S^c) \mu(\dd x)}{\mu(S)}\geq \frac{\omega\cdot\lambda\cdot\psi}{36} \cdot
      \log^{\frac{1}{2}}\Big(1+\frac{4}{\mu(S)}\Big)\geq   \frac{\omega\cdot\lambda\cdot\psi}{36} \cdot
      \log^{\frac{1}{2}}\Big(1+\frac{1}{v}\Big).
       \end{aligned}
  \end{equation*}
   Then combining with the result for the first case, we can obtain a lower bound of $$\frac{\omega}{4}\,\min\Big\{1,\frac{\lambda\cdot\psi}{9} \sqrt{\log\big(1+\frac{1}{v}\big)}\Big\}$$ on $s$-conductance profile $\Phi_s(v)$ with $s=\frac{\varepsilon^2}{16M_0^2}$.

 \subsection{Proof of Lemma~\ref{transform}}
 Recall the transition kernel associated with $\mu_k$,
\begin{equation*}
    T(\theta,\dd y)=\left[1-(1-\zeta)\cdot\int A(\theta,y) Q(\theta,y)\,\dd y\right]\delta_{\theta}(\dd y)+(1-\zeta)\cdot Q(\theta,y) A(\theta,y)\,\dd y
 \end{equation*}
 with
 \begin{equation*}
     A(\theta,y)=1\wedge \frac{\pi_{n}(y)Q(y,\theta)}{\pi_{n}(\theta)Q(\theta,y)}; \quad Q(\theta,\cdot)=N_d\Big(\theta-\frac{h}{\sqrt{n}}\widetilde I 
\widetilde\nabla V_n\big(\sqrt{n}(\theta-\wh\theta)\big),\frac{2h}{n}\,\widetilde I\Big).
 \end{equation*}
 Then given $\xi\in \mb R^d$, the distribution of $G_{\#}T(\xi,\cdot)$ is 
 \begin{equation*}
 \begin{aligned}
      &T^*(\theta,\dd z)\\
      &=\left[1-(1-\zeta)\cdot\int  Q^*(\theta,z) A\big(\theta, \wh\theta+\widetilde{I}^{\frac{1}{2}}\frac{z}{\sqrt{n}}\big)\,\dd z\right]\delta_{\sqrt{n}\wt I^{-\frac{1}{2}}(\theta-\wh\theta)}\big(\dd z\big)\\
    &+(1-\zeta)\cdot Q^*(\theta,z) A\big(\theta, \wh\theta+\widetilde{I}^{\frac{1}{2}}\frac{z}{\sqrt{n}}\big)\,\dd z,
     \end{aligned}
 \end{equation*}
 where $Q^*(\theta,\cdot)$ is the density function of $N_d( \sqrt{n}\widetilde{I}^{-\frac{1}{2}}(\theta-\wh\theta)-h \widetilde{I}^{\frac{1}{2}}\widetilde{\nabla}V_n\big(\sqrt{n}(\theta-\wh\theta)), 2hI_d\big)$. Then by the fact that 
 \begin{equation*}
 \begin{aligned}
\frac{Q\big(\wh\theta+\widetilde{I}^{\frac{1}{2}}\frac{z}{\sqrt{n}}, \wh\theta+\widetilde{I}^{\frac{1}{2}}\frac{\xi}{\sqrt{n}}\big)}{Q\big(\wh\theta+\widetilde{I}^{\frac{1}{2}}\frac{\xi}{\sqrt{n}},  \wh\theta+\widetilde{I}^{\frac{1}{2}}\frac{z}{\sqrt{n}}\big)}&=\exp\left(-\frac{1}{4h}\left(\| \xi-  z+h \widetilde{I}^{\frac{1}{2}}\widetilde \nabla V_n( \widetilde{I}^{\frac{1}{2}}z)\|^2-\|  z-  \xi+h \widetilde{I}^{\frac{1}{2}}\widetilde \nabla V_n( \widetilde{I}^{\frac{1}{2}} \xi)\|^2\right)\right)\\
      &=\frac{\wt Q (z, \xi)}{\wt Q ( \xi,z)},
     \end{aligned}
 \end{equation*}
 we have
 \begin{equation*}
T^*\Big(\wh\theta+\widetilde{I}^{\frac{1}{2}}\frac{\xi}{\sqrt{n}},\dd z\Big)= \left[1-(1-\zeta)\cdot\int \wt A ( \xi,z) \wt Q (\xi,z)\,\dd z\right]\mathbf{1}_{ \xi}(\dd z)+(1-\zeta)\cdot \wt Q( {\xi},z) \wt A ( \xi,z)\dd z=\wt T (\xi, \dd z).
 \end{equation*}
  Thus when $\widetilde{\mu}_{k-1}=G_{\#}\mu_{k-1}$, we have $\widetilde{\mu}_k=G_{\#}\mu_{k}$. Then combine with the fact that $ \widetilde\mu_0=G_{\#}\mu_0$, we can obtain by induction that  $\widetilde\mu_k=G_{\#}\mu_{k}$ for $k\in \mb N$. 
  \subsection{Proof of Lemma~\ref{lemmalogiso}}
 To begin with, we consider the following lemma stated in~\cite{JMLR:v21:19-441}.
 \begin{lemma}\label{lemmalogiso1}
 (Lemma 16 of~\cite{JMLR:v21:19-441}) Let \(\gamma\) denote the density of the standard Gaussian distribution \(\mathcal{N}\left(0, \sigma^{2} {I}_{d}\right)\), and let \(\mu\) be a distribution with density \(\mu=q \cdot \gamma\), where \(q\) is a log-concave function. Then for any partition \(S_{1}, S_{2}, S_{3}\) of \(\mathbb{R}^{d}\), we have
$$
\mu\left(S_{3}\right) \geq \frac{d\left(S_{1}, S_{2}\right)}{2 \sigma} \min \left\{\mu\left(S_{1}\right), \mu\left(S_{2}\right)\right\} \log ^{\frac{1}{2}}\left(1+\frac{1}{\min \left\{\mu\left(S_{1}\right), \mu\left(S_{2}\right)\right\}}\right) .
$$
 \end{lemma}
 We first consider the case $\wt J=I_d$ where recall $\wt J=\wt {I}^{\frac{1}{2}}J\wt {I}^{\frac{1}{2}}$. Then define $\overline{\pi}=N(0,I_d)|_{\wt K}$, by the fact that $\wt K=B_{R/2}^d$ is a convex set and $\mathbf{1}_{\wt K}$ is a log-concave function, using lemma~\ref{lemmalogiso1}, we can obtain that for any partition  \(S_{1}, S_{2}, S_{3}\) of \({\wt K}\), we have 
 $$
\overline\pi \left(S_{3}\right) \geq \frac{d\left(S_{1}, S_{2}\right)}{2} \min \left\{\overline\pi\left(S_{1}\right), \overline\pi\left(S_{2}\right)\right\} \log ^{\frac{1}{2}}\left(1+\frac{1}{\min \left\{\overline\pi\left(S_{1}\right), \overline\pi\left(S_{2}\right)\right\}}\right) .
$$
Then recall $\wt \pi_{\rm loc}|_{\wt K}(\xi)=\frac{\mathbf{1}_{\wt K}\exp(-V_n(\wt{I}^{\frac{1}{2}}\xi))}{\int_{\wt K}  \exp(-V_n(\wt{I}^{\frac{1}{2}}\xi))\dd\xi}$, using the fact that $\underset{\widetilde\xi\in B_{R}^d}{\sup}\big|V_n( \widetilde{I}^{\frac{1}{2}}\widetilde\xi)-\frac{1}{2}{\widetilde\xi^T \wt J\widetilde\xi} \big|\leq\widetilde{\varepsilon}_0$, we can obtain that for any measurable set $S\subseteq \wt K$, we have 
\begin{equation*}
    \exp(-2\widetilde{\varepsilon}_0)\leq \frac{\wt \pi_{\rm loc}|_{\wt K}(S)}{\overline\pi(S)}=\frac{\int_{S\cap \wt K}\exp(-V_n(\wt {I}^{\frac{1}{2}}\xi))\dd \xi \int_K\exp(-\frac{1}{2}\xi^T\xi)\dd \xi}{\int_{S\cap K}\exp(-\frac{1}{2}\xi^T\xi)\dd \xi\int_{K}\exp(-V_n(\wt{I}^{\frac{1}{2}}\xi))\dd \xi}\leq  \exp(2\widetilde{\varepsilon}_0).
\end{equation*}
Thus
\begin{equation}\label{eqn:caseI}
    \begin{aligned}
        &\wt \pi_{\rm loc}|_{\wt K}(S_3)\geq  \exp(-2\widetilde{\varepsilon}_0)\overline{\pi}(S_3)\\
        &\geq \frac{d(S_1,S_2)}{2} \exp(-2\widetilde{\varepsilon}_0)\min \left\{\overline\pi\left(S_{1}\right), \overline\pi\left(S_{2}\right)\right\} \log ^{\frac{1}{2}}\left(1+\frac{1}{\min \left\{\overline\pi\left(S_{1}\right), \overline\pi\left(S_{2}\right)\right\}}\right)\\
        &\overset{(i)}{\geq} \frac{d(S_1,S_2)}{2} \exp(-4\widetilde{\varepsilon}_0)\min \left\{ \wt\pi_{\rm loc}|_{\wt K}\left(S_{1}\right), \wt\pi_{\rm loc}|_{\wt K}\left(S_{2}\right)\right\} \log ^{\frac{1}{2}}\left(1+\frac{1}{\exp(-2\widetilde{\varepsilon}_0)\min \left\{ \wt\pi_{\rm loc}|_{\wt K}\left(S_{1}\right),  \wt\pi_{\rm loc}|_{\wt K}\left(S_{2}\right)\right\}}\right)\\
        &\geq  \frac{d(S_1,S_2)}{2} \exp(-4\widetilde{\varepsilon}_0)\min \left\{ \wt\pi_{\rm loc}|_{\wt K}\left(S_{1}\right),  \wt\pi_{\rm loc}|_{\wt K}\left(S_{2}\right)\right\} \log ^{\frac{1}{2}}\left(1+\frac{1}{\min \left\{ \wt\pi_{\rm loc}|_{\wt K}\left(S_{1}\right),  \wt\pi_{\rm loc}|_{\wt K}\left(S_{2}\right)\right\}}\right),\\
    \end{aligned}
\end{equation}
where $(i)$ uses the fact that $x\log^{\frac{1}{2}}(1+\frac{1}{x})$ is an increasing function.  For the general case where $\wt J$ is not necessary an identity matrix,  we can define $K'=\wt J^{\frac{1}{2}}\wt K=\{x=\wt J^{\frac{1}{2}}y\,: \,y\in \wt K\}$, and $\lambda=\wt J^{\frac{1}{2}}\xi$, where $\xi$ is a random variable with density $\pi_{\rm loc}|_{\wt K}$. Thus $\lambda$ has a density
\begin{equation*}
    \pi_{\lambda}(\lambda)=\frac{\mathbf{1}_{K'}(\lambda)\exp(-V_n(\wt I^{\frac{1}{2}}\wt J^{-\frac{1}{2}}\lambda))}{\int_{K'}\exp(-V_n(  \wt{I}^{\frac{1}{2}}J^{-\frac{1}{2}}\lambda))\dd \lambda},
\end{equation*}
Moreover, for any $\lambda\in K'$, it holds that 
\begin{equation*}
    \big|V_n(\wt{I}^{\frac{1}{2}}\wt J^{-\frac{1}{2}}\lambda)-\frac{1}{2}\lambda^T\lambda\big|\leq \widetilde{\varepsilon}_0.
\end{equation*}
Then   for any partition  \(S_{1}, S_{2}, S_{3}\) of \(\wt K\), let 
 \begin{equation*}
     \begin{aligned}
         &\widetilde{S_1}=\wt J^{\frac{1}{2}}S_1;\\
          &\widetilde{S_2}=\wt J^{\frac{1}{2}}S_2;\\
           &\widetilde{S_3}=\wt J^{\frac{1}{2}}S_3.\\
     \end{aligned}
 \end{equation*}
  Then by the positive definiteness of $\wt J$, $(\widetilde{S_1},\widetilde{S_2},\widetilde{S_3})$ forms a partition for $K'$, and 
  \begin{equation*}
      d(\widetilde{S_1},\widetilde{S_2})\geq \sqrt{\rho_1}\, d(S_1,S_2).
  \end{equation*}
 Since $K'$ is a convex set,  by applying $ \pi_{\lambda}$ to statement~\eqref{eqn:caseI}, we can obtain
 \begin{equation*}
     \begin{aligned}
         &\wt\pi_{\rm loc}|_{\wt K}(S_3)= \pi_{\lambda}(\widetilde S_3)
         \geq    \frac{d(\widetilde S_1,\widetilde S_2)}{2} \exp(-4\widetilde{\varepsilon}_0)\min \left\{  \pi_{\lambda}  (\widetilde S_{1} ),  \pi_{\lambda} (\widetilde S_{2} )\right\} \log ^{\frac{1}{2}}\bigg(1+\frac{1}{\min \left\{  \pi_{\lambda} (\widetilde S_{1} ),   \pi_{\lambda} (\widetilde S_{2})\right\}}\bigg)\\
         &\geq \frac{\sqrt{\rho_1}}{2} d(S_1,S_2)\exp(-4\widetilde{\varepsilon}_0)\min \left\{ \pi_{\rm loc}|_{\wt K}\left(S_{1}\right),  \pi_{\rm loc}|_{\wt K}\left(S_{2}\right)\right\} \log ^{\frac{1}{2}}\left(1+\frac{1}{\min \left\{ \wt\pi_{\rm loc}|_{\wt K}\left(S_{1}\right),  \wt\pi_{\rm loc}|_{\wt K}\left(S_{2}\right)\right\}}\right).
     \end{aligned}
 \end{equation*}
 Proof is completed.
 
  \subsection{Proof of Lemma~\ref{boundTV}}
  We first construct the high probability set $E$ as follows:  let
  $$r_d= \left(\sqrt{c'd\left(\log \left(\frac{M_0^2}{\varepsilon^2 h\rho_1}\right)+\wt\varepsilon_0\right)} \rho_2^2\right)\vee \left(c'\left(\log \left(\frac{M_0^2}{\varepsilon^2 h\rho_1}\right)+\wt\varepsilon_0\right)\rho_2^2\right),$$
  and $\wt J=\wt{I}^{\frac{1}{2}}J\wt{I}^{\frac{1}{2}}$.
 We define $E=\{\xi\in B_{R/2}^d: \big|\xi^T\wt J^3\xi-{\rm tr}(\wt J^2)\big|\leq r_d\}\cap \{\xi\in B_{R/2}^d: \big|\xi^T\wt J^2\xi-{\rm tr}(\wt J)\big|\leq r_d/\rho_2\}$.  By the choice of $h$,  when $c_0$ is small enough, it holds that 
  \begin{equation*}
  h\leq \sqrt{c_0}\cdot\left\{\left({\rho_2}^{-\frac{1}{3}} (\rho_2^2d+r_d)^{-\frac{1}{3}}\right)\wedge (r_d)^{-\frac{1}{2}}\right\}.
  \end{equation*}
Now we show that $E$ is indeed a high probability set in the following lemma. Note that all the following lemmas in this subsection are under Assumptions in Theorem~\ref{thmala}.
\begin{lemma}\label{lemmaprobE}
Consider $E=\{\xi\in B_{R/2}^d: \big|\xi^T\wt J^3\xi-{\rm tr}(\wt J^2)\big|\leq r_d\}\cap \{\xi\in B_{R/2}^d: \big|\xi^T\wt J^2\xi-{\rm tr}(\wt J)\big|\leq r_d/\rho_2\}$. If $r_d= \Big(\sqrt{c'd\log \left(\frac{M_0^2}{\varepsilon^2 h\rho_1}\right)} \rho_2^2\Big)\vee \Big(c'\log \left(\frac{M_0^2}{\varepsilon^2 h\rho_1}\right)\rho_2^2\Big)$ for a sufficiently large enough constant $c'$, then $\wt \pi_{\rm loc}(E)\geq 1-\exp(-4\wt\varepsilon_0)\cdot\frac{2\varepsilon^2 h\rho_1}{M_0^2}$.
\end{lemma}
 We now show that for any $x,z\in E$ with $\|x-z\|\leq\frac{\sqrt{h}}{3}$, the total variation distance between $\wt T_x=\wt T(x,\cdot)$ and $\wt T_z=\wt T(z,\cdot)$ can be upper bounded by $1-\frac{\exp(-2\wt\varepsilon_0)}{4}$.  For any $x,z\in E$, we consider the following decomposition:
\begin{equation*} 
    \begin{aligned}
     &\|\wt T_{x}-\wt T_{z}\|_{TV}\\
     &=\frac{1}{2}\int |\wt T(x,y)-\wt T(z,y)|\,\dd y\\
     &=\frac{1}{2}\wt T_{x}(\{x\})+\frac{1}{2}\wt T_{z}(\{z\})+\frac{1}{2}\int_{\mathbb{R}^d\backslash\{x,z\}}|\wt T(x,y)-\wt T(z,y)|\,\dd y\\
   &=\frac{1}{2}-\frac{1-\zeta}{2} \int_{\mb R^d} \wt Q(x,y)\wt A(x,y)\,\dd y+\frac{1}{2}-\frac{1-\zeta}{2} \int_{\mb R^d} \wt Q(z,y)\wt A(z,y)\,\dd y\\
   &\quad+\frac{1-\zeta}{2}\int_{\mb R^d} |\wt Q(x,y)\wt A(x,y)-\wt Q(z,y)\wt A(z,y)|\,\dd y\\
   &=1-(1-\zeta)\int_{\mb R^d} \min\left(\wt A(x,y)\wt Q(x,y), \wt A(z,y)\wt Q(z,y)\right)\,\dd y  \\
   &\leq 1-(1-\zeta)\int_{ B_{R}^d}\min\left(\wt A(x,y)\wt Q(x,y), \wt A(z,y)\wt Q(z,y)\right)\,\dd y  
    \end{aligned}
\end{equation*}
Recall that 
\begin{equation*}
    \wt A(x,y)=1\wedge \frac{\wt \pi_{\rm loc}(y)\wt Q(y,x)}{\wt \pi_{\rm loc}(x)\wt Q(x,y)},
\end{equation*}
where $\wt\pi_{\rm loc}(x)\propto \exp(-V_n(\wt I^{\frac{1}{2}}x))$ and 
 \begin{equation*}
 \underset{x\in B_{R}^d}{\sup}\big|V_n( \widetilde{I}^{\frac{1}{2}}x)-\frac{1}{2}{x^T \wt Jx} \big|\leq \widetilde{\varepsilon}_0.
 \end{equation*}
Define $\overline{\pi}$ as the density function  of $N_d(0,\wt J^{-1})$, we have
\begin{equation*}
    \frac{\wt \pi_{\rm loc}(y)}{\wt \pi_{\rm loc}(x)}=\frac{ \exp(-V_n(\wt I^{\frac{1}{2}}y))}{ \exp(-V_n(\wt I^{\frac{1}{2}}y))}\geq \exp(-2\wt\varepsilon_0)\cdot\frac{ \exp(-\frac{1}{2}{y^T \wt Jy} )}{ \exp(-\frac{1}{2}{x^T \wt Jx} )}=\exp(-2\wt\varepsilon_0)\cdot\frac{\ov \pi(y)}{\ov \pi(x)}.
\end{equation*}
Therefore, denote 
\begin{equation*}
    \ov A(x,y)=1\wedge \frac{\ov \pi(y)\wt Q(y,x)}{\ov \pi(x)\wt Q(x,y)},
\end{equation*}
we have 
\begin{equation*}
      \wt A(x,y)\geq  1\wedge \frac{\exp(-2\wt\varepsilon_0)\cdot\ov \pi(y)\wt Q(y,x)}{\ov \pi(x)\wt Q(x,y)}\geq  \exp(-2\wt\varepsilon_0)\cdot     \ov A(x,y).
\end{equation*}
We can then derive 
\begin{equation}
    \begin{aligned}
      &\|\wt T_{x}-\wt T_{z}\|_{TV}\\
      &\leq  1-(1-\zeta)\int_{ B_{R}^d}\min\left(\wt A(x,y)\wt Q(x,y), \wt A(z,y)\wt Q(z,y)\right)\,\dd y \\
      &\leq   1-(1-\zeta) \exp(-2\wt\varepsilon_0)\cdot\int_{ B_{R}^d}\min\left(\ov A(x,y)\wt Q(x,y), \ov A(z,y)\wt Q(z,y)\right)\,\dd y\\
      &= 1-\frac{1}{2}(1-\zeta) \exp(-2\wt\varepsilon_0) \cdot \bigg(\int_{ B_{R}^d}\ov A(x,y)\wt Q(x,y)\,\dd y+\int_{ B_{R}^d}\ov A(z,y)\wt Q(z,y)\,\dd y\\
      &\qquad-\int_{ B_{R}^d}\big|\ov A(x,y)\wt Q(x,y)-\wt A(z,y)\wt Q(z,y)\big|\,\dd y\bigg)
    \end{aligned}
     \end{equation}
 Then consider the inequality:
\begin{equation*}
\begin{aligned}
 \int_{B_R^d}|\wt Q(x,y)\ov A(x,y)-\wt Q(z,y)\ov A(z,y)|\,\dd y&\leq \int_{B_R^d} \wt Q(x,y)(1-\ov A(x,y))\,\dd y\\
 &+ \int_{B_R^d} \wt Q(z,y)(1-\ov A(z,y))\,\dd y+2\|\wt Q_x-\wt Q_z\|_{  \rm TV},
 \end{aligned}
\end{equation*}
 where we use $\wt Q_x$ to denote the probability measure with density function $\wt Q(x,\cdot)$. Moreover, consider the equation:
 \begin{equation*}
 \begin{aligned}
     \int_{ B_{R}^d}\ov A(x,y)\wt Q(x,y)\,\dd y &=  \int_{ B_{R}^d}(\ov A(x,y)-1)\wt Q(x,y)\,\dd y+ \int_{ B_{R}^d} \wt Q(x,y)\,\dd y\\
     &= 1-\int_{B_R^d} \wt Q(x,y)(1-\ov A(x,y))\,\dd y-\int_{ (B_{R}^d)^c} \wt Q(x,y)\,\dd y.
      \end{aligned}
 \end{equation*}
Combined with~\eqref{decompTXz}, we can obtain 
 \begin{equation}\label{decompTXz}
     \begin{aligned}
        &\|\wt T_{x}-\wt T_{z}\|_{TV}\\
       &\leq 1-(1-\zeta) \exp(-2\wt\varepsilon_0) \cdot \bigg(1-\int_{B_R^d} \wt Q(x,y)(1-\ov A(x,y))\,\dd y-\int_{B_R^d} \wt Q(z,y)(1-\ov A(z,y))\,\dd y\\
       &\qquad -\|\wt Q_x-\wt Q_z\|_{  \rm TV}-\frac{1}{2}\int_{ (B_{R}^d)^c} \wt Q(x,y)\,\dd y-\frac{1}{2}\int_{ (B_{R}^d)^c} \wt Q(z,y)\,\dd y\bigg)
     \end{aligned}
 \end{equation}
 Consider the proposal distribution of MALA for sampling from  the Gaussian $\overline{\pi}:\,=N_d(0,\wt J^{-1})$,
  \begin{equation*}
 Q^{\Delta}_x(\cdot)= N_d(x-h\wt Jx,2hI_d),   
 \end{equation*}
 whose density is denoted as $Q^{\Delta}(x,\cdot)$. Then $\|\wt Q_x-\wt Q_z\|_{  \rm TV}\leq \|\wt Q_x-Q^{\Delta}_x\|_{  \rm TV}+ \|  Q_x^{\Delta}-  Q_{z}^{\Delta}\|_{  \rm TV} + \|\wt Q_z-Q^{\Delta}_z\|_{  \rm TV}$ can be upper bounded by Pinsker's inequality, that is, 
for any $x\in B_R^d$, 
 \begin{equation*}
     \|\wt Q_x- Q^{\Delta}_x\|_{  \rm TV}\leq \frac{1}{2}\sqrt{\frac{h^2\|\widetilde{I}^{\frac{1}{2}}\widetilde{\nabla}V_n( \widetilde{I}^{\frac{1}{2}}x)-\wt Jx\|^2}{2h}}\leq \frac{\sqrt{h}\widetilde{\varepsilon}_1\mnorm{\wt I^{\frac{1}{2}}}_{\rm op}}{2\sqrt{2}},
 \end{equation*}
 and for any $x,z \in B_R^d$
  \begin{equation*}
     \|Q_x^{\Delta}-Q_{z}^{\Delta}\|_{  \rm TV}\leq \frac{1}{2}\sqrt{\frac{\|(I-h\wt J)(x-z)\|^2}{2h}}\leq \frac{\|x-z\|}{2\sqrt{2h}}.
 \end{equation*}
Therefore, when $\|x-z\|\leq \frac{\sqrt{h}}{3}$ and $ \sqrt{h}\widetilde{\varepsilon}_1\mnorm{\wt I^{\frac{1}{2}}}_{\rm op}\leq \frac{\sqrt{2}}{36}$,  we have 
\begin{equation*}
    \|\wt Q_x-\wt Q_z\|_{  \rm TV}\leq  \frac{\|x-z\|}{2\sqrt{2h}}+\frac{\sqrt{h}\widetilde{\varepsilon}_1\mnorm{\wt I^{\frac{1}{2}}}_{\rm op}}{\sqrt{2}}<\frac{1}{6}.
\end{equation*}

For the term of $\int_{B_R^d}\wt Q(x,y)(1-\ov A(x,y))\,\dd y$, we use Condition A by comparing $Q_x$ with $Q_x^{\Delta}$, leading to the following decomposition:
\begin{equation*}
    \begin{aligned}
     &\int_{B_R^d}\wt Q(x,y)(1-\ov A(x,y))\,\dd y\\
     &\leq \int_{B_R^d}\Big|\wt Q(x,y)-\frac{\ov\pi(y)\wt Q(y,x)}{ \ov\pi(x)}\Big|\,\dd y\\
     &\leq 2\|\wt Q_x-Q_x^{\Delta}\|_{  \rm TV}+\underbrace{\int\left|Q^{\Delta}(x,y)-\frac{\ov\pi(y)Q^{\Delta}(y,x)}{\ov\pi(x)}\right|\dd y}_{\rm (A)}+\underbrace{\int_{B_R^d}\left|\frac{\ov\pi(y)Q^{\Delta}(y,x)}{\ov\pi(x)}-\frac{\ov\pi(y)\wt Q(y,x)}{\ov\pi(x)}\right|\,\dd y}_{\rm (B)}.
    \end{aligned}
\end{equation*}
 We then state the following lemma for bounding the term (A).
\begin{lemma}\label{lemma3}
 Consider the choice of (rescaled) step size $h$ in Theorem~\ref{thmala}, then when $c_0$ is small enough and $x\in E$, it holds that 
\begin{equation*}
    \int\left|Q^{\Delta}(x,y)-\frac{\ov\pi(y)Q^{\Delta}(y,x)}{\ov\pi(x)}\right|\,\dd y\leq \frac{1}{24}.
\end{equation*}
\end{lemma}
Our proof of Lemma~\ref{lemma3} is technically similar to that of Proposition 38 in~\cite{pmlr-v134-chewi21a} for bounding the mixing time of MALA with a standard Gaussian target (i.e.~$\overline{\pi}=N_d(0,I_d)$). The non-trivial part in our analysis lies in keeping track of the dependence on the maximal and minimal eigenvalues of $J$. We then bound the term (B)  by the following lemma. 
\begin{lemma}\label{boundCD}
Consider the choice of (rescaled) step size $h$ in Theorem~\ref{thmala},  then when $c_0$ is small enough, for any $x\in E$, it holds that 
\begin{equation*}
\begin{aligned}
     &\int_{B_R^d} \left|Q^{\Delta}(y,x)-Q(y,x)\right|\frac{\ov{\pi}(y)}{\ov{\pi}(x)}\,\dd y\leq  \frac{1}{72}.\\
    \end{aligned}
\end{equation*}
\end{lemma}
Thus when $\|x-z\|\leq \frac{\sqrt{h}}{3}$ and $ \sqrt{h}\widetilde{\varepsilon}_1\mnorm{\wt I^{\frac{1}{2}}}_{\rm op}\leq \frac{\sqrt{2}}{36}$, 
\begin{equation}\label{eqn:boundQA}
    \begin{aligned}
     &\int_{B_R^d} \wt Q(x,y)(1-\ov A(x,y))\,\dd y\\
         &\leq 2 \|\wt Q_x-Q_x^{\Delta}\|_{  \rm TV}+\frac{1}{24}+\frac{1}{72}\\
         &\leq \sqrt{\frac{h}{2}}\widetilde{\varepsilon}_1 \mnorm{\wt I^{\frac{1}{2}}}_{\rm op}+\frac{1}{18}\leq \frac{1}{12}.
    \end{aligned}
\end{equation}
 Finally, since for any $x\in E\subset B_{R/2}^d$, 
\begin{equation*}
\begin{aligned}
    \int_{\mb (B_R^d)^c} \wt Q(x,y)\,\dd y&\leq   \int_{\mb (B_R^d)^c}  Q^\Delta(x,y)\,\dd y+2\, \|\wt Q_x-Q_x^\Delta\|_{\rm TV}\\
    &\leq \mb{E}_{u\in  N_d(0,I_d)}\Big[\mathbf{1}\big(\|u\|\geq \frac{R}{2\sqrt{2h}}\big)\Big]+\frac{\sqrt{h}\wt\varepsilon_1 \mnorm{\wt I^{\frac{1}{2}}}_{\rm op}}{\sqrt{2}}.
\end{aligned}
\end{equation*}
Since $R\geq 8\sqrt{d/\lambda_{\min}(\wt J)}$, when the constant $c_0$ in $h$ is small enough,  we can obtain
\begin{equation*}
   \int_{\mb (B_R^d)^c} \wt Q(x,y)\,\dd y\leq \frac{1}{6}.
\end{equation*}
Then combined with the bound in equation~\eqref{eqn:boundQA} and decomposition~\eqref{decompTXz}, we can obtain that when $c_0$ is small enough, for any $x,z\in E$ with $\|x-z\|<\frac{\sqrt{h}}{3}$ and $\zeta\in (0,\frac{1}{2}]$, it holds that 
  \begin{equation*} 
     \begin{aligned}
        &\|\wt T_{x}-\wt T_{z}\|_{TV}\\
       &\leq 1-(1-\zeta) \exp(-2\wt\varepsilon_0) \cdot \bigg(1-\int_{B_R^d} \wt Q(x,y)(1-\ov A(x,y))\,\dd y-\int_{B_R^d} \wt Q(z,y)(1-\ov A(z,y))\,\dd y\\
       &\qquad -\|\wt Q_x-\wt Q_z\|_{  \rm TV}-\frac{1}{2}\int_{ (B_{R}^d)^c} \wt Q(x,y)\,\dd y-\frac{1}{2}\int_{ (B_{R}^d)^c} \wt Q(z,y)\,\dd y\bigg)\\
       &\leq  1-\frac{1-\zeta}{2} \exp(-2\wt\varepsilon_0)  \\
       &\leq 1-\frac{\exp(-2\wt\varepsilon_0)}{4}.
     \end{aligned}
 \end{equation*}

 \subsection{Proof of Lemma~\ref{lemmaprobE}}
We can write $\wt\pi_{\rm loc}$ as
  \begin{equation*}
      \wt\pi_{\rm loc}(\xi)=\frac{\frac{\sqrt{{\rm det}(\wt J)}}{(2\pi)^{\frac{d}{2}}}\exp(-V_n(\wt I^{\frac{1}{2}}\xi))}{\bigintss \frac{\sqrt{{\rm det}(\wt J)}}{(2\pi)^{\frac{d}{2}}}\exp(-V_n(\wt I^{\frac{1}{2}}\xi))\,\dd\xi}.
  \end{equation*}
Then
 \begin{equation*}
 \begin{aligned}
   1-\wt\pi_{\rm loc}(E)&\leq \frac{\bigintss_{ \big\{\xi\in B_{R/2}^d\,:\,|\xi^T\wt J^3\xi-{\rm tr}(\wt J^2)|> r_d\big\}}\frac{\sqrt{{\rm det}(\wt J)}}{(2\pi)^{\frac{d}{2}}}\exp(-V_n(\wt I^{\frac{1}{2}}\xi))\,\dd\xi}{\bigintss \frac{\sqrt{{\rm det}(\wt J)}}{(2\pi)^{\frac{d}{2}}}\exp(-V_n(\wt I^{\frac{1}{2}}\xi))\,\dd\xi}\\
   &+\frac{\bigintss_{ \big\{\xi\in B_{R/2}^d\,:\,|\xi^T\wt J^2\xi-{\rm tr}(\wt J)|> r_d/\rho_2\big\}}\frac{\sqrt{{\rm det}(\wt J)}}{(2\pi)^{\frac{d}{2}}}\exp(-V_n(\wt I^{\frac{1}{2}}\xi))\,\dd\xi}{\bigintss \frac{\sqrt{{\rm det}(\wt J)}}{(2\pi)^{\frac{d}{2}}}\exp(-V_n(\wt I^{\frac{1}{2}}\xi))\,\dd\xi}\\
   &+\wt\pi_{\rm loc}(\|\xi\|>R/2).
      \end{aligned}
 \end{equation*}
  Then for the denominator, as
  \begin{equation*}
 \begin{aligned}
 \underset{\widetilde\xi\in B_{R}^d}{\sup}\big|V_n( \widetilde{I}^{\frac{1}{2}}\widetilde\xi)-\frac{1}{2}{\widetilde\xi^T \widetilde{I}^{\frac{1}{2}}J \widetilde{I}^{\frac{1}{2}}\widetilde\xi} \big|\leq \widetilde{\varepsilon}_0,
       \end{aligned}
 \end{equation*}
 when $R\geq 8 (\frac{d}{\lambda_{\min}(\wt J)})^{\frac{1}{2}}$, we can obtain that 
  \begin{equation*}
  \begin{aligned}
     & \bigintss \frac{\sqrt{{\rm det}(\wt J)}}{(2\pi)^{\frac{d}{2}}}\exp(V_n(\wt I^{-\frac{1}{2}}\xi))\,\dd\xi\\
     &\geq \bigintss_{B_R^d}\frac{\sqrt{{\rm det}(\wt J)}}{(2\pi)^{\frac{d}{2}}}\exp(-\frac{\xi^T\wt J\xi}{2})\exp(\frac{\xi^T\wt J\xi}{2}-V_n(I^{\frac{1}{2}}\xi))\,\dd\xi\\
     &\geq \exp(-\widetilde{\varepsilon}_0)\bigintss_{B_R^d}\frac{\sqrt{{\rm det}(\wt J)}}{(2\pi)^{\frac{d}{2}}}\exp(-\frac{\xi^T\wt J\xi}{2})\,\dd\xi \\
     &\geq \frac{1}{2}\exp(-\wt\varepsilon_0).
       \end{aligned}
  \end{equation*}
  Furthermore, by Bernstein's inequality (see for example, Theorem 2.8.2 of~\cite{vershynin_2018}), for $x\sim N_d(0,\Sigma)$, it holds that 
  \begin{equation}\label{Bernstein}
    \mb P(\left|\|x\|^2-{\rm tr}(\Sigma)\right|\geq t)\leq 2\exp(-\frac{1}{8}(\frac{t^2}{\mnorm{\Sigma}_{  \rm F}^2}\wedge \frac{t}{\mnorm{\Sigma}_{  \rm  op}}))
  \end{equation}
    We can then obtain
 \begin{equation*}
 \begin{aligned}
   \pi_{\rm loc}(E)  &\geq 1-2\exp(2\widetilde{\varepsilon}_0) \int_{  \big\{|\xi^T\wt J^3\xi-{\rm tr}(\wt J^2)|> r_d\big\}}\frac{\sqrt{{\rm det}(\wt J)}}{(2\pi)^{\frac{d}{2}}}\exp(-\frac{\xi^T\wt J\xi}{2})\,\dd\xi\\
   &-2\exp(2\widetilde{\varepsilon}_0) \int_{  \big\{|\xi^T\wt J^2\xi-{\rm tr}(\wt J)|> r_d/\rho_2\big\}}\frac{\sqrt{{\rm det}(\wt J)}}{(2\pi)^{\frac{d}{2}}}\exp(-\frac{\xi^T\wt J\xi}{2})\,\dd\xi-\frac{\varepsilon^2 h\rho_1}{M_0^2}\\
     &\geq 1-\exp(-4\wt\varepsilon_0)\cdot\frac{2\varepsilon^2 h\rho_1}{M_0^2},
      \end{aligned}
 \end{equation*}
  where the last inequality is due to the Bernstein's inequality in~\eqref{Bernstein}.  
 
\subsection{Proof of Lemma~\ref{lemma3}}
Recall $\overline{\pi}=N_d(0,\wt J^{-1})$ and $Q^{\Delta}(x,\cdot)$ be the density of $N_d(x-h\wt Jx,2hI_d)$, we have
 \begin{equation*}
     \begin{aligned}
&\int\left|Q^{\Delta}(x,y)-\frac{\ov\pi(y)Q^{\Delta}(y,x)}{\ov\pi(x)}\right|\,\dd y\\
&=\int \frac{1}{(4\pi h)^{\frac{d}{2}}}\left|\exp\left(-\frac{\|y-x+h\wt Jx\|^2}{4h}\right)-\exp\left(\frac{x^T\wt Jx-y^T\wt Jy}{2}\right)\exp\left(-\frac{\|x-y+h\wt Jy\|^2}{4h}\right)\right|\,\dd y\\
&=\int  \frac{1}{(4\pi h)^{\frac{d}{2}}}\exp\left(-\frac{\|y-x+h\wt Jx\|^2}{4h}\right)\left|1-\exp\left(\frac{h^2\|\wt Jx\|^2-h^2\|\wt Jy\|^2}{4h}\right)\right|\,\dd y,
\end{aligned}
 \end{equation*}
let $u=\frac{y-x+h\wt Jx}{\sqrt{2h}}$ in the above integral, then consider $u\sim N_d(0,I_d)$  and let
$$\m A=\left\{u\in \mb R^d: \frac{1}{4}\left|2h^2\|\wt Ju\|^2+2\sqrt{2}h^{\frac{3}{2}}x^T\wt J^2u-2\sqrt{2}h^{\frac{5}{2}}x^T\wt J^3u+h^3x^T\wt J^4x-2h^2x^T\wt J^3x\right|\leq\frac{1}{49}\right\}.$$  
We can then  obtain
\begin{equation}\label{eqn:2}
     \begin{aligned}
 &\int\left|Q^{\Delta}(x,y)-\frac{\ov\pi(y)Q^{\Delta}(y,x)}{\ov\pi(x)}\right|\,\dd y\\   
&=\mathbb{E}_{u}\left[\left|1-\exp\left(\frac{-h^2\|\sqrt{2h}\wt Ju+\wt Jx-h\wt J^2x\|^2+h^2\|\wt Jx\|^2}{4h}\right)\right|\right]\\
&=\mathbb{E}_{u}\left[\left|1-\exp\left(-\frac{1}{4}\left(2h^2\|\wt Ju\|^2+2\sqrt{2}h^{\frac{3}{2}}x^T\wt J^2u-2\sqrt{2}h^{\frac{5}{2}}x^T\wt J^3u+h^3x^T\wt J^4x-2h^2x^T\wt J^3x\right)\right)\right|\right]\\
&\leq \bigg\{\mathbb{E}_{u}\bigg[\bigg|1-\exp\Big(-\frac{1}{4}\Big(2h^2\|\wt Ju\|^2+2\sqrt{2}h^{\frac{3}{2}}x^T\wt J^2u-2\sqrt{2}h^{\frac{5}{2}}x^T\wt J^3u+h^3x^T\wt J^4x-2h^2x^T\wt J^3x\Big)\Big)\bigg|\\
&\cdot\mathbf{1}_{\m A}(u)\bigg]\bigg\}+\big\{\mathbb{E}_{u}\left[\mathbf{1}_{\m A^c}(u)\right]\big\}+\bigg\{\exp\left(-\frac{1}{4}h^3x^T\wt J^4x \right)\sqrt{\mathbb{E}_{u}\left[\mathbf{1}_{\m A^c}(u)\right]}\\
&\cdot\left(\mathbb{E}_{u} \left[\exp(-3h^2(u^T\wt J^2u-x^T\wt J^3x))\right]\cdot\mathbb{E}_{u} \left[\exp(3\sqrt{2}h^{\frac{3}{2}}x^T\wt J^2u)\right]\cdot\mathbb{E}_{u}\left[ \exp(3\sqrt{2}h^{\frac{5}{2}}x^T\wt J^3u)\right]\right)^{\frac{1}{6}}\bigg\},\\
     \end{aligned}
 \end{equation}
 where the last inequality uses H\"{o}lder inequality. The first term of the right hand side of equation~\eqref{eqn:2} can be upper bound by $\exp(1/49)-1\leq 1/48$. For the second and third term, by (1) $h\leq\sqrt{c_0}{ \rho_2^{-\frac{1}{3}}({\rm tr}(\wt J^2)+r_d)^{-\frac{1}{3}}}$ and $h\leq \sqrt{c_0} r_d^{-\frac{1}{2}}$ with  $r_d= \bigg\{\Big(\sqrt{c'\log \frac{M_0^2}{\varepsilon^2 h\rho_1}}\mnorm{\wt J^2}_{  \rm F} \Big)\vee \Big(c'\log \frac{M_0^2 }{\varepsilon^2 h\rho_1}\rho_2^2\Big)\bigg\}\wedge  (\rho_2^3\|K\|^2)$ and $\|K\|\geq C(\frac{d}{\rho_1})^{\frac{1}{2}}$; (2) $x\in E=\{x\in K: \big|x^T\wt J^3x-{\rm tr}(\wt J^2)\big|\leq r_d\}$,  it holds that
 \begin{equation*}
 \begin{aligned}
  &h^3x^T\wt J^4x\leq  h^3\rho_2x^T\wt J^3x\leq h^3\rho_2(r_d+{\rm tr}(\wt J^2))\leq  c_0^{\frac{3}{2}}.\\
   \end{aligned}
 \end{equation*}
 Moreover, since for a Gaussian random variable $\bar{u}\sim N(0,\sigma^2)$, it holds that 
 \begin{equation*}
 \begin{aligned}
  & \mathbb{E}\exp(t \bar u)=\exp(\frac{\sigma^2t^2}{2})\\
    & \mathbb{E}\exp(-t^2 \bar u^2)=\frac{1}{\sqrt{1+2t^2\sigma^2}} \quad |t|<\sqrt{\frac{1}{2\sigma^2}}.
      \end{aligned}
 \end{equation*}
 We can get
 \begin{equation*}
 \begin{aligned}
     &\mathbb{E}_{u} \left[\exp(t^2h^2(x^T\wt J^3x-\|\wt Ju\|^2))\right]\\&\leq \exp(t^2h^2(x^T\wt J^3x-{\rm tr}(\wt J^2)))\prod_{j=1}^d \frac{1/\sqrt{1+2t^2 h^2 \lambda_j(\wt J^2)}}{\exp\big(-t^2 h^2\lambda_j(\wt J^2) \big)}\\
     &\leq \exp(t^2h^2r_d)\cdot\prod_{j=1}^d\big(1+C\,t^4h^4(\lambda_j(\wt J^2))^2\big)\\
     &\leq   \exp(t^2c_0)\exp(Ct^4h^4\mnorm{\wt J^2}_{\rm F}^2)\\
     &\leq \exp(t^2c_0+t^4C\,c_0^2), \quad |t|\leq \sqrt{\frac{1}{4h^2\rho_2(\wt J^2)}},
     \end{aligned}
 \end{equation*}
 where the last inequality uses $h\leq\sqrt{c_0}{ \rho_2^{-\frac{1}{3}}({\rm tr}(\wt J^2)+r_d)^{-\frac{1}{3}}}\leq \sqrt{c_0}{ \rho_2^{-\frac{1}{3}}({\rm tr}(\wt J^2) )^{-\frac{1}{3}}}\leq \sqrt{c_0}\|\wt J^2\|_{F}^{-\frac{1}{2}}$,  and
  \begin{equation*}
    \mathbb{E}_{u} \left[\exp(th^{\frac{3}{2}}x^T\wt J^2u)\right]\leq \exp\left(\frac{1}{2}t^2h^3\|x^T\wt J^2\|^2\right)\leq\exp\left(\frac{1}{2}t^2h^3\rho_2({\rm tr}(\wt J^2)+r_d)\right)\leq \exp\left(\frac{1}{2}c_0^{\frac{3}{2}}t^2\right);
 \end{equation*}
 
 \begin{equation*}
     \mathbb{E}_{u}\left[ \exp( th^{\frac{5}{2}}x^T\wt J^3u)\right]\leq \exp\left(\frac{1}{2}t^2h^5\|x^T\wt J^3\|^2\right)\leq \exp\left(\frac{1}{2}t^2h^5\rho_2^3({\rm tr}(\wt J^2)+r_d)\right)\leq\exp(\frac{1}{2}c_0^{\frac{5}{2}}t^2),
 \end{equation*}
 where the last inequality uses $h\leq\sqrt{c_0}{ \rho_2^{-\frac{1}{3}}({\rm tr}(\wt J^2)+r_d)^{-\frac{1}{3}}}\leq \sqrt{c_0}{ \rho_2^{-1}}$. Then by Markov inequality,  we can obtain that
 \begin{equation*}
     \mb P_u\left(|h^{\frac{3}{2}}x^T\wt J^2u|\geq \frac{1}{96\sqrt{2}}\right)\leq 2\,\underset{t>0}{\inf}\exp\Big(\frac{1}{2}c_0^{\frac{3}{2}}t^2-\frac{t}{96\sqrt{2}}\Big)= 2\exp\left(-\frac{1}{2\cdot(96\sqrt{2})^2c_0^{\frac{3}{2}}}\right);
 \end{equation*}
 \begin{equation*}
     \mb P_u\left(|h^{\frac{5}{2}}x^T\wt J^3u|\geq \frac{1}{96\sqrt{2}}\right)\leq 2\,\underset{t>0}{\inf}\exp\Big(\frac{1}{2}c_0^{\frac{5}{2}}t^2-\frac{t}{96\sqrt{2}}\Big)= 2\exp\left(-\frac{1}{2\cdot(96\sqrt{2})^2c_0^{\frac{5}{2}}}\right).
 \end{equation*}
 Also, by Bernstein's inequality in~\eqref{Bernstein}, we have
 \begin{equation*}
 \begin{aligned}
     \mb P_u\left(h^2\left|\|\wt Ju\|^2-x^T\wt J^3x\right|\geq \frac{1}{96}\right)&\leq P_u\left(\left|\|\wt Ju\|^2-{\rm tr}(\wt J^2)\right|\geq \frac{1}{96h^2}-r_d\right)\\
     &\leq P_u\left(\left|\|\wt Ju\|^2-{\rm tr}(\wt J^2)\right|\geq \frac{1}{h^2}(\frac{1}{96}-c_0)\right)\\
     &\leq  2\exp\left(-\frac{1}{c'}\left(\frac{\frac{1}{96}-c_0}{h^2\rho_2^2}\wedge \frac{(\frac{1}{96}-c_0)^2}{h^4\|\wt J^2\|^2_F}\right)\right)\\
     &\leq 2\exp\left(-\frac{1}{c'}\left(\frac{\frac{1}{96}-c_0}{c_0}\wedge \frac{(\frac{1}{96}-c_0)^2}{c_0^2}\right)\right),
 \end{aligned}
 \end{equation*}
 where the last inequality uses $h\leq \sqrt{c_0}\|\wt J^2\|_{F}^{-\frac{1}{2}}$. Therefore, when $c_0$ is small enough, we have
 \begin{equation*}
 \begin{aligned}
     &\mathbb{E}_{u}\left[\mathbf{1}_{\m A^c}(u)\right]\\
     &\leq \mb P_u\Big(h^2|\|\wt Ju\|^2-x^T\wt J^3x|\geq \frac{1}{96}\Big)+ \mb P_u\Big(|h^{\frac{3}{2}}x^T\wt J^2u|\geq \frac{1}{96\sqrt{2}}\Big)+ \mb P_u\Big(|h^{\frac{5}{2}}x^T\wt J^3u|\geq \frac{1}{96\sqrt{2}}\Big)\\
     &\leq 2\exp\Big(-\frac{\frac{1}{96}-c_0}{c'c_0}\Big) +2\exp\Big(-\frac{1}{2\cdot(96\sqrt{2})^2c_0^{\frac{3}{2}}}\Big)+2\exp\Big(-\frac{1}{2\cdot(96\sqrt{2})^2c_0^{\frac{5}{2}}}\Big)
      \end{aligned}
 \end{equation*}
and  
 \begin{equation*}
     \begin{aligned}
      &\mathbb{E}_{u}\left[\mathbf{1}_{\m A^c}(u)\right]+\exp\left(-\frac{1}{4}h^3x^T\wt J^4x\right)\sqrt{\mathbb{E}_{u}\left[ \mathbf{1}_{\m A^c}(u)\right]}\\
&\cdot \left(\mathbb{E}_{u} \left[\exp(-3h^2(\|\wt Ju\|^2-x^T\wt J^3x))\right]\cdot\mathbb{E}_{u} \left[\exp(3\sqrt{2}h^{\frac{3}{2}}x^T\wt J^2u)\right]\cdot\mathbb{E}_{u}\left[ \exp(3\sqrt{2}h^{\frac{5}{2}}x^T\wt J^3u)\right]\right)^{\frac{1}{6}}\\
&\leq \frac{1}{48}.
     \end{aligned}
 \end{equation*}
 We can then obtain the desired result by combining all pieces.

 \subsection{Proof of Lemma~\ref{boundCD}}
We first write
\begin{equation*}
   \begin{aligned}
   & \int_{B_R^d} \left|Q^{\Delta}(y,x)-\wt Q(y,x)\right|\frac{\ov{\pi}(y)}{\ov{\pi}(x)}\,\dd y\\
    &=\int_{B_R^d} \left|1-\frac{\wt Q(y,x)}{Q^{\Delta}(y,x)}\right|\frac{\ov{\pi}(y)}{\ov{\pi}(x)}Q^{\Delta}(y,x)\,\dd y\\
    &=\int_{B_R^d} \left|1-\exp\left(\frac{-\|x-y+h\widetilde{I}^{\frac{1}{2}}\widetilde{\nabla}V_n( \widetilde{I}^{\frac{1}{2}}y)\|^2+\|x-y+h\wt Jy\|^2}{4h}\right)\right|\frac{\ov{\pi}(y)}{\ov{\pi}(x)}Q^{\Delta}(y,x)\,\dd y.\\
     \end{aligned} 
\end{equation*}
Since  $h\leq\sqrt{c_0}{ \rho_2^{-\frac{1}{3}}({\rm tr}(\wt J^2)+r_d)^{-\frac{1}{3}}}\leq \sqrt{c_0}{ \rho_2^{-1}}$ and $h\rho_2\mnorm{\wt I}_{\rm op} R^2 \wt\varepsilon_1^2\leq c_0$, when $c_0$ is sufficiently small, we have for any $x\in E$ and $y\in {B_R^d}$,
\begin{equation*}
 \begin{aligned}
  & \frac{\left|-\|x-y+h\widetilde{\nabla}\widetilde{I}^{\frac{1}{2}}\widetilde{\nabla}V_n( \widetilde{I}^{\frac{1}{2}}y)\|^2+\|x-y+h\wt Jy\|^2\right|}{4h} \\
  &=\frac{\left|h(\wt Jy+\widetilde{I}^{\frac{1}{2}}\widetilde{\nabla}V_n( \widetilde{I}^{\frac{1}{2}}y))^T(\wt Jy-\widetilde{I}^{\frac{1}{2}}\widetilde{\nabla}V_n( \widetilde{I}^{\frac{1}{2}}y))+2(x-y)^T(\wt Jy-\widetilde{I}^{\frac{1}{2}}\widetilde{\nabla}V_n( \widetilde{I}^{\frac{1}{2}}y))\right|}{4}\\
  &\leq \frac{h\big(2\rho_2R+ \mnorm{\wt I^{\frac{1}{2}}}_{\rm op}\widetilde{\varepsilon}_1\big) \mnorm{\wt I^{\frac{1}{2}}}_{\rm op}\widetilde{\varepsilon}_1+2\|x-y\| \mnorm{\wt I^{\frac{1}{2}}}_{\rm op}\widetilde{\varepsilon}_1}{4}\\
  &\leq \frac{\sqrt{c_0}}{4}\big(3+\frac{2\|x-y\|}{R\sqrt{h\rho_2} }\big).
 \end{aligned}   
\end{equation*}
Thus we can bound 
\begin{equation*}
     \int_{B_R^d} \left|Q^{\Delta}(y,x)-\wt Q(y,x)\right|\frac{\ov{\pi}(y)}{\ov{\pi}(x)}\,\dd y\leq \int_{B_R^d} \Big(\exp\Big(\frac{\sqrt{c_0}}{4}\big(3+\frac{2\|x-y\|}{R\sqrt{h\rho_2} }\big)\Big)-1\Big) \frac{\ov{\pi}(y)}{\ov{\pi}(x)}Q^{\Delta}(y,x)\,\dd y.
\end{equation*}
Furthermore, by Lemma~\ref{lemma3}, we can get
\begin{equation*}
  \int \frac{\ov{\pi}(y)}{\ov{\pi}(x)}Q^{\Delta}(y,x)\,\dd y\leq \int Q^{\Delta}(x,y)\,\dd y+  \int\left|Q^{\Delta}(x,y)-\frac{\ov\pi(y)Q^{\Delta}(y,x)}{\ov\pi(x)}\right|\,\dd y\leq \frac{25}{24},
\end{equation*}
which leads to
\begin{equation*}
    \begin{aligned}
   &\int_{B_R^d} \left|Q^{\Delta}(y,x)-\wt Q(y,x)\right|\frac{\ov{\pi}(y)}{\ov{\pi}(x)}\,\dd y\\
   &\leq \frac{25}{24}\int_{B_R^d} \Big(\exp\Big(\frac{\sqrt{c_0}}{4}\big(3+\frac{2\|x-y\|}{R\sqrt{h\rho_2} }\big)\Big)-1\Big) \frac{\frac{\ov{\pi}(y)}{\ov{\pi}(x)}Q^{\Delta}(y,x)}{  \int \frac{\ov{\pi}(y)}{\ov{\pi}(x)}Q^{\Delta}(y,x)\,\dd y}\,\dd y\\
     &=\frac{25}{24}\int_{B_R^d} \Big(\exp\Big(\frac{\sqrt{c_0}}{4}\big(3+\frac{2\|x-y\|}{R\sqrt{h\rho_2} }\big)\Big)-1\Big) N_d\big((I+h^2\wt J)^{-1}(x-h\wt Jx),2h(I+h^2\wt J)^{-1}\big)\,\dd y,
    \end{aligned}
\end{equation*}
where the last inequality is due to $\frac{\ov{\pi}(y)}{\ov{\pi}(x)}Q^{\Delta}(y,x)\propto \exp(-\frac{y^T(I+h^2J)y-2y^T(x-h\wt Jx)}{4h})$.  Consider $u\sim  N_d(0,I_d)$, for sufficiently small $c_0$, we have 
\begin{equation*}
    \begin{aligned}
 &\int_{B_R^d} \Big(\exp\Big(\frac{\sqrt{c_0}}{4}\big(3+\frac{2\|x-y\|}{R\sqrt{h\rho_2} }\big)\Big)-1\Big)  N_d\big((I+h^2\wt J)^{-1}(x-h\wt Jx),2h(I+h^2\wt J)^{-1}\big)\,\dd y\\
  &\leq \mb{E}_{u\sim  N_d(0,I_d)}\bigg[\exp\Big(\frac{\sqrt{c_0}}{4}\Big(3+\frac{2\|(I+h^2\wt J)^{-1}(x-h\wt Jx)-x+\sqrt{2h}(I+h^2\wt J)^{-\frac{1}{2}}u\|}{R\sqrt{h\rho_2} }\Big)\Big)-1\bigg]\\
  &\overset{(i)}{\leq } \mb{E}_{u\sim  N_d(0,I_d)}\bigg[\exp\Big(\frac{\sqrt{c_0}}{4}\Big(3+2c_0^{\frac{1}{4}}+\frac{2\sqrt{2}}{R\sqrt{\rho_2}}\|u\| \Big)\Big)-1\bigg]\\
  &\leq \frac{81}{80} \cdot\mb{E}_{u\sim  N_d(0,I_d)}\bigg[\exp\Big(\frac{\sqrt{2c_0}}{2R\sqrt{\rho_2}}\|u\|\Big)\bigg]-1\\
   &\leq \frac{81}{80} \cdot\sqrt{\mb{E}_{u\sim  N_d(0,I_d)}\bigg[\exp\Big(\frac{c_0}{2R^2\rho_2}\|u\|^2\Big)\bigg]}-1\\
   &\leq  \frac{81}{80}\exp(\frac{c_0d}{2R^2\rho_2})-1\\
   &\overset{(ii)}{\leq} \frac{1}{75},
    \end{aligned}
\end{equation*}
where $(i)$ is due to $\|(I+h^2\wt J)^{-1}(x-h\wt Jx)-x\|\leq h^2\rho_2\|x\|+h\rho_2\|x\|\leq 2\sqrt{h\rho_2}c_0^{\frac{1}{4}}\|x\|\leq \sqrt{h\rho_2}c_0^{\frac{1}{4}}R$, and $(ii)$ is due to $R\geq 8\sqrt{d/\lambda_{\min}(\wt J)}$. Thus we can obtain $ \int_{B_R^d} \left|Q^{\Delta}(y,x)-\wt Q(y,x)\right|\frac{\ov{\pi}(y)}{\ov{\pi}(x)}\,\dd y\leq \frac{1}{72}$.  
\subsection{Proof of Lemma~\ref{lemma:controlsgap}}
 \subsubsection{Proof of statement (1) of Lemma~\ref{lemma:controlsgap}} Define the following compact supported function $k:\mb R\to \mb R$:
    \begin{equation*}
        k(t)=\left\{
        \begin{array}{cc}
          2(t-t^3)  & t\in (-1,1), \\
            0 & \text{otherwise}.
        \end{array}
        \right.
    \end{equation*}
    Then consider a initial distribution with density function $\mu_0(x)=(1+k(\sqrt{\rho_1}x_d))\cdot\ov\pi(x)$. This constriction guarantees that
    \begin{equation*}
    \begin{aligned}
         \chi^2(\mu_0,\ov\pi)&=\sqrt{\frac{\rho_1}{2\pi}}\int_{-\sqrt{\frac{1}{\rho_1}}}^{\sqrt{\frac{1}{\rho_1}}} k^2(\sqrt{\rho_1}x_d)\exp(-\frac{\rho_1}{2}x_d^2)\,\dd x_d\\
         &=\sqrt{\frac{1}{2\pi}}\int_{-1}^{1} k^2(t) \exp(-\frac{1}{2}t^2)\,\dd t\in (0.2,0.21),
    \end{aligned}
    \end{equation*}
    \begin{equation*}
      \underset{x\in \mb R^d}{\sup} \frac{\mu_0(x)}{\ov \pi(x)}= 1+\underset{t\in (-1,1)}{\sup}k(t)<2,
    \end{equation*}
    and 
    \begin{equation*}
        |h_0(x)-h_0(y)|=\big|\frac{\mu_0(x)}{\ov\pi(x)}-\frac{\mu_0(y)}{\ov\pi(y)}\big|=|k(\sqrt{\rho_1}x_d)-k(\sqrt{\rho_1}y_d)|\leq 2\sqrt{\rho_1}|x_d-y_d|.
    \end{equation*}
    Therefore, the spectral gap of this initialization is controlled by 
    \begin{equation*}
    \begin{aligned}
         \frac{\mathcal{E}(h_0,h_0)}{\chi^2(\mu_0',\ov\pi)}&\leq 10\rho_1\cdot\mb{E}_{x\in \ov \pi,y \in T(x,\cdot)}\big[(x_d-y_d)^2\big]\\
         &\leq 10\rho_1\cdot \mb{E}_{x\in \ov \pi,y \in  N_d(x-h\rho_1 x,2hI_d)}\big[(x_d-y_d)^2\big]\\
         &= 10\rho_1\cdot \mb{E}_{x_d\in N(0,1/\rho_1),\xi\in N(0,1)}\big[(h\rho_1 x_d-\sqrt{2h}\xi)^2\big]\\
         &\leq 20m^2h^2+ 40 mh\leq 60mh.\\
    \end{aligned}
        \end{equation*}
        
   \subsubsection{Proof of Statement (2) of Lemma~\ref{lemma:controlsgap}}

    Denote sets
    \begin{equation*}
    \begin{aligned}
        & K_2=\Big\{x\in \mb R^d\,:\,\big|x^TJ^3x-{\rm tr}(J^2)\big|\leq (5\mnorm{J^2}_{\rm F})\vee (24\mnorm{J^2}_{\rm op})\Big\};\\
              & K_3=\Big\{x\in \mb R^d\,:\,\big|x^TJ^4x-{\rm tr}(J^3)\big|\leq (5\mnorm{J^3}_{\rm F})\vee (24\mnorm{J^3}_{\rm op})\Big\};\\
      & K_4=\Big\{x\in \mb R^d\,:\,\big|x^TJ^6x-{\rm tr}(J^5)\big|\leq (5\mnorm{J^5}_{\rm F})\vee (24\mnorm{J^5}_{\rm op})\Big\},\\
    \end{aligned}
      \end{equation*}
    To control the probability of the above events, we utilize the following Bernstein's inequality: for $x\in  N_d(0,\Sigma)$,
  \begin{equation}\label{Bernstein1}
    \mb P(\left|\|x\|^2-{\rm tr}(\Sigma)\right|\geq t)\leq 2\exp\Big(-\frac{1}{8}\big(\frac{t^2}{\mnorm{\Sigma}_{  \rm F}^2}\wedge \frac{t}{\mnorm{\Sigma}_{  \rm  op}}\big)\Big),
  \end{equation}
  which leads to 
  \begin{equation*}
       \mb P\left(\left|\|x\|^2-{\rm tr}(\Sigma)\right|\geq \big(\sqrt{8\lambda}\mnorm{\Sigma}_F\big)\vee \big(8\lambda\mnorm{\Sigma}_{\rm op}\big)\right)\leq 2\exp(-\lambda).
  \end{equation*}
 Therefore, for $x\sim N_d(0,J^{-1})$, the probability of events $x\in K_2\cap K_3\cap K_4$ is inside the interval of $(0.7,1)$. Then let $K_1\subset \mb R^d$ be an arbitrary measurable set so that the probability of events $x\in K=K_1\cap K_2\cap K_3\cap K_4$ is equal $\frac{1}{M_0}$ (notice that $M_0\geq 2$ and $\frac{1}{M_0}\leq \frac{1}{2}<0.7$, therefore such a set $K_1$ exists). Then consider a initial distribution with density function $\mu_0'(x)=\frac{\ov \pi(x)\textbf{1}_K(x)}{\mb{E}_{\ov \pi}[\mathbf{1}_K(x)]}$, it holds that 
 \begin{equation*}
     \chi^2(\mu_0',\ov\pi)=\frac{1}{\mb{E}_{\ov \pi}[\mathbf{1}_K(x)]}-1=M_0-1,
 \end{equation*}
 and 
 \begin{equation*}
         \underset{x\in \mb R^d}{\sup} \frac{\mu_0(x)}{\ov \pi(x)}= \frac{1}{\mb{E}_{\ov \pi}[\mathbf{1}_K(x)]}=M_0.
 \end{equation*}
 Then denote for bounding the spectral gap $\frac{\mathcal{E}(h_0,h_0)}{\chi^2(\mu_0',\ov\pi)}$ with $h_0=\frac{\dd \mu_0'}{\dd \ov\pi}$, we claim it suffices to show the following claim: denote $Q^{\Delta}(x,)$ to be the density function of $ N_d(x-hJx,2hI_d)$, then for any $x\in K$, there exists a set $G_x\subset \mb R^d$ so that
  \begin{equation}\label{claimsgap1}
     \begin{aligned}
         \frac{\ov \pi(y)Q^{\Delta}(y,x)}{\ov \pi(x)Q^{\Delta}(x,y)}\leq \exp(-16\log (\kappa d)), \quad \forall\, y\in G_x,
     \end{aligned}
 \end{equation}
 and 
 \begin{equation}\label{claimsgap2}
     \begin{aligned}
         \int_{G_x}Q^{\Delta}(x,y)\,\dd y\geq 1-\frac{3}{\kappa d}.
     \end{aligned}
 \end{equation}
 Indeed, under claim~\eqref{claimsgap1} and~\eqref{claimsgap2}, we have 
 \begin{equation*}
 \begin{aligned}
      \frac{\mathcal{E}(h_0,h_0)}{\chi^2(\mu_0',\ov\pi)}&=\frac{M_0^2\cdot\mb{E}_{x\in \ov\pi, y\in T(x,\cdot)}\big[(\mathbf{1}_K(x)-\mathbf{1}_K(y))^2\big]}{2(M_0-1)}\\
      &\leq\frac{M_0^2}{M_0-1}\int_{K}\int_{K^c} \min\Big\{1,\frac{\ov \pi(y)Q^{\Delta}(y,x)}{\ov \pi(x)Q^{\Delta}(x,y)}\Big\}\ov\pi(x)Q^{\Delta}(x,y)\,\dd y\,\dd x\\
      &\leq \frac{M_0^2}{M_0-1}\int_{x\in K}\Big(\int_{G_x} \frac{\ov \pi(y)Q^{\Delta}(y,x)}{\ov \pi(x)Q^{\Delta}(x,y)} Q^{\Delta}(x,y)\,\dd y+\int_{G_x^c} Q^{\Delta}(x,y)\,\dd y\Big)\ov\pi(x)\,\dd x\\
      &\leq \frac{M_0}{M_0-1}\underset{x\in K}{\sup} \Big(\underset{y\in G_x}{\sup} \frac{\ov \pi(y)Q^{\Delta}(y,x)}{\ov \pi(x)Q^{\Delta}(x,y)}+\int_{G_x^c} Q^{\Delta}(x,y)\Big) \,\dd y\\
      &\leq \frac{8}{\kappa d},
     \end{aligned}
 \end{equation*}
 where the last inequality uses claim~\eqref{claimsgap1} and~\eqref{claimsgap2}. Now we show the desired claim. First note that
 \begin{equation*}
     \frac{\ov \pi(y)Q^{\Delta}(y,x)}{\ov \pi(x)Q^{\Delta}(x,y)}=\exp\Big(\frac{h\cdot(x^TJ^2x-y^TJ^2y)}{4}\Big).
 \end{equation*}
 Let $u=\frac{y-x+hJx}{\sqrt{2h}}$, then for $y\in  N_d(x-hJx,2hI_d)$, we have $u\in  N_d(0,I_d)$. Therefore, it suffice to show that for any $x\in K$, there exists  a set $G_x'\in \mb R^d$ so that $\mb{E}_{ N_d(0, I_d)}[\mathbf{1}_{G_x'}(u)]\geq 1-\frac{3}{\kappa d}$ and 
 \begin{equation*}
     \begin{aligned}
          \frac{\ov \pi(\sqrt{2h}u+x-hJx)Q^{\Delta}(\sqrt{2h}u+x-hJx,x)}{\ov \pi(x)Q^{\Delta}(x,\sqrt{2h}u+x-hJx)}\leq \exp(-16\log (\kappa d)), \quad \forall\, u\in G_x'.
     \end{aligned}
 \end{equation*}
  Denote the sets
  \begin{equation*}
      \begin{aligned}
          &\m G_x^1=\{u\in \mb R^d\,:\,\|Ju\|^2-x^TJ^3x\geq -\big(\sqrt{8\log (\kappa d)}+5\big)\rho_2^2\sqrt{d}\};\\
           &\m G^2_x=\{u\in \mb R^d\,:\,x^TJ^2u\geq -\sqrt{\log (\kappa d)\cdot\big(10\rho_2^3\sqrt{d}+2\rho_2^3 d\big)}\};\\
           &\m G^3_x=\{u\in \mb R^d\,:\,x^TJ^3u\leq \sqrt{\log (\kappa d)\cdot\big(10\rho_2^5\sqrt{d}+2\rho_2^5 d\big)}\}.
      \end{aligned}
  \end{equation*}
 Then under $ G'_x=\m G_x^1\cap \m G_x^2\cap \m G_x^3$, we have 
 \begin{equation*}
 \begin{aligned}
       &\frac{\ov \pi(\sqrt{2h}u+x-hJx)Q^{\Delta}(\sqrt{2h}u+x-hJx,x)}{\ov \pi(x)Q^{\Delta}(x,\sqrt{2h}u+x-hJx)}\\
       &=\exp\left(-\frac{1}{4}\left(2h^2\| Ju\|^2+2\sqrt{2}h^{\frac{3}{2}}x^T  J^2u-2\sqrt{2}h^{\frac{5}{2}}x^T J^3u+h^3x^T  J^4x-2h^2x^T J^3x\right)\right)\\
       &\leq \exp\bigg(-\frac{1}{4}\Big(h^3x^TJ^4x-2h^2\big(\sqrt{8\log (\kappa d)}+5\big)\rho_2^2\sqrt{d}\\
       &\qquad -2\sqrt{2}\sqrt{\log (\kappa d)}\big(h^{\frac{3}{2}}\sqrt{ 10\rho_2^3\sqrt{d}+2\rho_2^3 d}+h^{\frac{5}{2}}\sqrt{ 10\rho_2^5\sqrt{d}+2\rho_2^5 d}\big)\Big)\bigg).\\
 \end{aligned}
 \end{equation*}
 Then there exists a universal constant $N_1$ so that when $d\geq N_1$, for any $x\in K$ and $u\in G_x'$, 
 \begin{equation*}
 \begin{aligned}
       &\frac{\ov \pi(\sqrt{2h}u+x-hJx)Q^{\Delta}(\sqrt{2h}u+x-hJx,x)}{\ov \pi(x)Q^{\Delta}(x,\sqrt{2h}u+x-hJx)}\\
       &\leq \exp\Big(-\frac{1}{4}\big(h^3x^TJ^4x-6h^2\sqrt{\log (\kappa d)}\rho_2^2\sqrt{d} -5\sqrt{\log (\kappa d)}\big(h^{\frac{3}{2}}\sqrt{\rho_2^3 d}+h^{\frac{5}{2}}\sqrt{\rho_2^5 d}\big)\big)\Big)\\
       &\leq \exp\Big(-\frac{1}{4}\big(h^3x^TJ^4x-6h^2\sqrt{\log (\kappa d)}\rho_2^2\sqrt{d} -5\sqrt{\log (\kappa d)}\big(h^{\frac{3}{2}}\sqrt{\rho_2^3 d}+h^{\frac{5}{2}}\sqrt{\rho_2^5 d}\big)\big)\Big)\\
       &\leq \exp\big(-32\log (\kappa d)+96 \log^{\frac{7}{6}} (\kappa d)d^{-\frac{1}{6}}+227\log^{\frac{4}{3}} (\kappa d)d^{-\frac{1}{3}}\big).
 \end{aligned}
 \end{equation*}
 Therefore, use $\kappa\leq c_1\cdot d^{c_2}$ there exists $N_2$ that depends only on $c_1,c_2$ so that when $d\geq N_2$,  for any $x\in K$ and $u\in G_x'$
 \begin{equation*}
 \begin{aligned}
       &\frac{\ov \pi(\sqrt{2h}u+x-hJx)Q^{\Delta}(\sqrt{2h}u+x-hJx,x)}{\ov \pi(x)Q^{\Delta}(x,\sqrt{2h}u+x-hJx)}\leq \exp\big(-16\log (\kappa d)\big).
 \end{aligned}
 \end{equation*}
 Now we control the probability $u\in G_x'$.  Firstly by Bernstein's inequality, for $u\in N_d(0,I_d)$, we have 
 \begin{equation*}
     \mb P\left(u^TJ^2u- {\rm tr}(J^2)\geq -\big((\sqrt{8\log (\kappa d)}\mnorm{J^2}_{\rm F})\vee (8\log (\kappa d)\mnorm{J^2}_{\rm op})\big)\right)\leq 1-\frac{1}{d\kappa}.
 \end{equation*}
 So there exists a universal constant $N_3$ so that when $d\geq N_3$, it holds with probability at least $1-\frac{1}{d\kappa}$ that
 \begin{equation*}
 \begin{aligned}
      \|Ju\|^2-x^TJ^3x&\geq {\rm tr}(J^2)-(\sqrt{8\log (\kappa d)}\mnorm{J^2}_{\rm F})\vee (8\log (\kappa d)\mnorm{J^2}_{\rm op})-{\rm tr}(J^2)- (5\mnorm{J^2}_{\rm F})\vee (24\mnorm{J^2}_{\rm op})\\
     &\geq -\big(\sqrt{8\log (\kappa d)}+5\big)\rho_2^2\sqrt{d}.
 \end{aligned}
 \end{equation*}
 Moreover, since for any $t\in \mb R$ and $x\in K$,
 \begin{equation*}
 \begin{aligned}
       \mb{E}[\exp(tx^TJ^2u)]&=\exp\left(\frac{1}{2}t^2 x^TJ^4x\right)\\
       &\leq \exp\Big(\frac{1}{2}\Big({\rm tr}(J^3)+\big((5\mnorm{J^3}_{\rm F})\vee (24\mnorm{J^3}_{\rm op})\big)\Big)t^2\Big),
 \end{aligned}
 \end{equation*}
 and 
  \begin{equation*}
 \begin{aligned}
       \mb{E}[\exp(tx^TJ^3u)]&=\exp\left(\frac{1}{2}t^2 x^TJ^6x\right)\\
       &\leq \exp\Big(\frac{1}{2}\Big({\rm tr}(J^5)+\big((5\mnorm{J^5}_{\rm F})\vee (24\mnorm{J^5}_{\rm op})\big)\Big)t^2\Big),
 \end{aligned}
 \end{equation*}
 by Markov inequality, there exists a universal constant $N_4$ so that when $d\geq N_4$,  it holds with probability at least $1-\frac{2}{\kappa d}$ that 
 \begin{equation*}
     x^TJ^2u\geq -\sqrt{\log (\kappa d)2\rho_2^3 d+10\rho_2^3\sqrt{d}}
 \end{equation*}
 and
  \begin{equation*}
     x^TJ^3u\leq \sqrt{\log (\kappa d)2\rho_2^5 d+10\rho_2^5\sqrt{d}}.
 \end{equation*}
 We can then obtain the desired result by combining all pieces.

 \section{Proof of Lemmas for Theorem~\ref{th:Gibbsmixing}}
  \subsection{Proof of Lemma~\ref{th1}}
Without loss of generality, we can assume the learning rate $\alpha=1$, as otherwise we can take $\ell(X,\theta)=\alpha\cdot \ell(X,\theta)$. To begin with, we provide in the following lemma some localized ``maximal” type inequalities that control the supreme of  empirical processes to deal with the non-smoothness of the loss function. All the following lemmas in this subsection are under Condition B.1-B.4.
 \begin{lemma}\label{lemma1}
 There exist positive constants $c$ and $r$ such that it holds with probability larger than $1-n^{-2}$ that ,
  \begin{enumerate}
      \item  For any $\theta,\theta'\in B_r(\theta^\ast)$, $\Big\|\frac{1}{n}\sum_{i=1}^ng(X_i,\theta)-\frac{1}{n}\sum_{i=1}^ng(X_i,\theta')-\mathbb{E}[g(X,\theta)]+\mathbb{E}[g(X,\theta')]\Big\|\leq C\,\bigg(\sqrt{\frac{\log n}{n}}d^{\frac{1+\gamma_3}{2}}\|\theta-\theta'\|^{\beta_1}+\frac{\log n}{n}d^{1+\gamma}\bigg)$.
      \item For any  $\theta,\theta'\in   \Theta$, $\Big|\frac{1}{n}\sum_{i=1}^n\ell(X_i,\theta)-\frac{1}{n}\sum_{i=1}^n\ell(X_i,\theta')-\mathbb{E}[\ell(X,\theta)]+\mathbb{E}[\ell(X,\theta')]\Big|\leq \\ C\,\bigg( \sqrt{\frac{\log n}{n}} d^{\frac{1}{2}+\gamma}\|\theta-\theta'\|+\frac{\log n}{n}d^{\frac{3}{2}+\gamma}\bigg)$.
     \item For any $\theta,\theta'\in B_r(\theta^\ast)$, $\Big|\frac{1}{n}\sum_{i=1}^n\ell(X_i,\theta)-\frac{1}{n}\sum_{i=1}^n\ell(X_i,\theta')-\frac{1}{n}\sum_{i=1}^ng(X_i,\theta')(\theta-\theta')-\mathbb{E}[\ell(X,\theta)]+\mathbb{E}[\ell(X,\theta')]+\mb E[g(X,\theta')(\theta-\theta')]\Big|\leq C\,\bigg(\sqrt{\frac{\log n}{n}}d^{\frac{1+\gamma_3}{2}}\|\theta-\theta'\|^{\beta_1+1}+\frac{\log n}{n}d^{ 1+\gamma}\|\theta-\theta'\|+(\frac{\log n}{n})^2\bigg)$.
  \end{enumerate}
 \end{lemma}
 
\noindent Recall $V_n(\xi)=n\,(\m R_n(\widehat{\theta}+\frac{\xi}{\sqrt{n}})-\m R_n(\widehat{\theta}))+\log \pi(\widehat{\theta}+\frac{\xi}{\sqrt{n}})-\log \pi(\widehat\theta)$, in order to bound the difference between $V_n(\xi)$ and $\frac{\xi^T\m H_{\theta^\ast}\xi}{2}$ using Lemma~\ref{lemma1}, we should first prove that $\widehat{\theta}$ is close to $\theta^\ast$.   Define a first order approximate to $\widehat{\theta}$:  $\widehat{\theta}^{\diamond}=\theta^\ast-\frac{1}{n}\sum_{i=1}^n \m H_{\theta^\ast}^{-1}g(X_i,\theta^\ast)$, we have the following lemma for bounding the difference between  $\widehat{\theta}^{\diamond}$ and $\theta^\ast$.
 
 \begin{lemma}\label{lemmathetadiamond}
  It holds with probability larger than $1-n^{-2}$ that 
  \begin{equation*}
      \|\widehat{\theta}^{\diamond}-\theta^\ast\|\leq C\, d^{\frac{1+\gamma_4}{2}}\sqrt\frac{\log n}{n}+C\,d^{1+\gamma_0+\gamma} \frac{\log n}{n}.
    \end{equation*}
 \end{lemma}
 And we resort to the following lemma that provides an upper bound on the $\ell_2$ distance between $\widehat\theta$ and $\widehat\theta^{\diamond}$.
\begin{lemma}\label{lemma2}
There exists a small enough positive constant $c$ such that when $d\leq c\big((\frac{n}{\log n})^{\frac{1}{2+2(\gamma+\gamma_0+\gamma_1)}}\wedge (\frac{n}{\log n})^{\frac{1}{1+2\gamma+2\gamma_2+4\gamma_0}}\big)$, then it holds with probability larger than $1-c\cdot n^{-2}$ that 
 \begin{equation}\label{diffhattheta}
     \left\|\frac{1}{n}\sum_{i=1}^n g(X_i,\widehat{\theta})\right\|\leq C\, d^{1+\gamma} \frac{\log n}{n}+ C\, d^{\frac{1+\gamma_3}{2}}\big(\frac{\log n}{n}\big)^{\frac{1}{2}+\beta_1};
\end{equation}
\begin{equation}
    \|\widehat\theta-\widehat\theta^\diamond\|\leq  C\,d^{\frac{1+\gamma_3}{2}+\beta_1(\frac{1+\gamma_4}{2})+\gamma_0}\big(\frac{\log n}{n}\big)^{\frac{1+\beta_1}{2}}+C\, d^{1+\gamma\vee(\gamma_2+\gamma_4)+\gamma_0}\frac{\log n}{n}+C\,\Big(d^{\frac{1+\gamma_3}{2}+\gamma_0}\sqrt\frac{\log n}{n}\Big)^{\frac{1}{1-\beta_1}}.
\end{equation}
\end{lemma}
\noindent By $\sup_{x\in \m X} \|g(X,\theta^\ast)\|\leq C\,d^{\gamma}$, we have $\mnorm{\m H_{\theta^\ast}^{-1}\mb E[g(X,\theta^\ast)g(X,\theta^\ast)^T]\m H_{\theta}^{-1}}_{  \rm  op}\leq C_1\, d^{2\gamma}\mnorm{\m H_{\theta^\ast}^{-1} }_{  \rm  op}^2 \leq C_2\, d^{2\gamma+2\gamma_0}$, which leads to  $\gamma_4\leq 2\gamma_0+2\gamma$. Then by Lemma~\ref{lemmathetadiamond} and Lemma~\ref{lemma2}, when $$d\leq c\,\Big(\big(\frac{n}{\log n}\big)^{\frac{\beta_1}{\gamma_3+\beta_1(1+\gamma_4)+2\gamma_0-\gamma_4}}\wedge\big(\frac{n}{\log n}\big)^{\frac{1}{1+2\gamma+2\gamma_2+4\gamma_0}}\wedge \big(\frac{n}{\log n}\big)^{\frac{1}{2+2(\gamma+\gamma_0+\gamma_1)}}\Big),$$  it holds with probability larger than $1-c_1\cdot n^{-2}$ that
\begin{equation}\label{diffhattrue}
    \|\widehat\theta-\theta^\ast\|\leq C\, d^{\frac{1+\gamma_4}{2}}\sqrt{\frac{\log n}{n}}.
\end{equation}
We can now derive (high probability) upper bound to the term of $|V_n(\xi)-\frac{\xi^T\m H_{\theta^\ast}\xi}{2}|$ over $1\leq \|\xi\|\leq C\sqrt{n}$. Consider the following decomposition:
\begin{equation*}
    \begin{aligned}
     &\left|V_n(\xi)-\frac{\xi^T\m H_{\theta^\ast}\xi}{2}\right|\\
     &\leq \left|n\,\big(\m R_n(\widehat{\theta}+\frac{\xi}{\sqrt{n}})-\m R_n(\widehat{\theta})\big)-\frac{\xi^T\m H_{\theta^\ast}\xi}{2}\right| +\left|\log \pi(\widehat{\theta}+\frac{\xi}{\sqrt{n}})-\log \pi(\widehat\theta)\right|\\
     &\leq  n\,\left|\frac{1}{n}\sum_{i=1}^ng(X_i,\widehat{\theta})\frac{\xi}{\sqrt n}\right|+n\,\bigg|\m R_n(\widehat{\theta}+\frac{\xi}{\sqrt{n}})-\m R_n(\widehat{\theta})-\frac{1}{n}\sum_{i=1}^ng(X_i,\widehat{\theta})\frac{\xi}{\sqrt n}-\Big(\m R(\widehat{\theta}+\frac{\xi}{\sqrt{n}})-\m R(\widehat{\theta})\\
     &-\mb Eg(X,\widehat{\theta})\frac{\xi}{\sqrt n}\Big)\bigg|+n\left|\m R(\widehat{\theta}+\frac{\xi}{\sqrt{n}})-\m R(\widehat{\theta})-\mb E[g(X,\widehat{\theta})]\frac{\xi}{\sqrt n}-\frac{\xi^T\m H_{\theta^\ast}\xi}{2}\right|+C\,\sqrt{d}\cdot\frac{\|\xi\|}{\sqrt n}.
    \end{aligned}
\end{equation*}
The first term can be bounded by Lemma~\ref{lemma2}, that is
\begin{equation*}
    \left|\frac{1}{n}\sum_{i=1}^ng(X_i,\widehat{\theta})\frac{\xi}{\sqrt n}\right|\leq C\, \frac{\|\xi\|}{\sqrt n}\left[d^{1+\gamma} \frac{\log n}{n}+  d^{\frac{1+\gamma_3}{2}}\big(\frac{\log n}{n}\big)^{\frac{1}{2}+\beta_1}\right];
\end{equation*}
for the second term, by the third statement of Lemma~\ref{lemma1}, we can obtain that 
\begin{equation*}
    \begin{aligned}
     &\left|\m R_n(\widehat{\theta}+\frac{\xi}{\sqrt{n}})-\m R_n(\widehat{\theta})-\frac{1}{n}\sum_{i=1}^ng(X_i,\widehat{\theta})\frac{\xi}{\sqrt n}-\left(\m R(\widehat{\theta}+\frac{\xi}{\sqrt{n}})-\m R(\widehat{\theta})-\mb Eg(X,\widehat{\theta})\frac{\xi}{\sqrt n}\right)\right|\\
     &\leq C\, \bigg[d^{1+\gamma}\frac{\log n}{n}\frac{\|\xi\|}{\sqrt n}+\sqrt\frac{\log n}{n} d^{\frac{1+\gamma_3}{2}}\big(\frac{\|\xi\|}{\sqrt n}\big)^{1+\beta_1}+\big(\frac{\log n}{n}\big)^2\bigg];
    \end{aligned}
\end{equation*}
for the third term, by the twice differentiability  of $\m R(\theta)$ and Lipschitzness of $\m H_{\theta}$, we can obtain that 
\begin{equation*}
    \begin{aligned}
     &\left|\m R(\widehat{\theta}+\frac{\xi}{\sqrt{n}})-\m R(\widehat{\theta})-\mb E[g(X,\widehat{\theta})]\frac{\xi}{\sqrt n}-\frac{\xi^T\m H_{\theta^\ast}\xi}{2n}\right|\\
     &\leq \frac{\|\xi\|^2}{2n} \underset{\xi\in K}{\sup}  \mnorm{\m H_{\widehat\theta+\frac{\xi}{\sqrt{n}}}-\m H_{\theta^\ast}}_{  \rm  op}\\
     &\leq C\, \frac{\|\xi\|^2}{n}d^{\gamma_2}\left(d^{\frac{1+\gamma_4}{2}}\sqrt\frac{\log n}{n}+\frac{\|\xi\|}{\sqrt n}\right)\\
     &= C\,\frac{\|\xi\|^3}{n^{\frac{3}{2}}}d^{\gamma_2}+C\,\frac{\|K\|^2}{n} \sqrt\frac{\log n}{n}d^{\frac{1+\gamma_4}{2}+\gamma_2}.
    \end{aligned}
\end{equation*}
 Therefore, by combining all these result, when $1\leq \|\xi\|\leq c\sqrt{n}$ for a small enough $c$, we can obtain that 
\begin{equation*}
    \begin{aligned}
      \left|V_n(\xi)-\frac{\xi^T\m H_{\theta^\ast}\xi}{2}\right|&\leq C\,d^{1+\gamma}\|\xi\|  \frac{\log n}{\sqrt{n}} + C\,d^{\frac{1+\gamma_3}{2}} \|\xi\|^{1+\beta_1}n^{-\frac{\beta_1}{2}}\sqrt{\log n}\\
      &+C\,d^{\frac{1+\gamma_4}{2}+\gamma_2}\|\xi\|^2\sqrt\frac{\log  n}{n}+C\,d^{\gamma_2}\|\xi\|^3 n^{-\frac{1}{2}}.
    \end{aligned}
\end{equation*}
For the second statement, since when $1\leq \|\xi\|\leq c\sqrt{n}$ for a small enough $c$,
\begin{equation*}
    \begin{aligned}
      &\|\widetilde{\nabla}V_n(\xi)-\m H_{\theta^\ast}\xi\|\\
      &=\bigg\|\frac{1}{\sqrt{n}}\sum_{i=1}^n g(X_i,\frac{\xi}{\sqrt{n}}+\widehat\theta)-\frac{1}{\sqrt{n}}\nabla [\log \pi](\frac{\xi}{\sqrt{n}}+\widehat\theta)-\m H_{\theta^*}\xi\bigg\|\\
      &\leq \left\|\frac{1}{\sqrt{n}}\sum_{i=1}^n g(X_i,\widehat\theta)\right\|+\sqrt{n}\left\|\frac{1}{n}\sum_{i=1}^n g(X_i,\frac{\xi}{\sqrt{n}}+\widehat\theta)-\frac{1}{n}\sum_{i=1}^n g(X_i,\widehat\theta)-\mathbb{E}[g(X,\frac{\xi}{\sqrt{n}}+\widehat\theta)]+\mathbb{E}[g(X,\widehat{\theta})]\right\|\\
      &+\sqrt{n}\left\|\mathbb{E}[g(X,\frac{\xi}{\sqrt{n}}+\widehat\theta)]-\mathbb{E}[g(X,\widehat{\theta})]-\m H_{\theta^\ast}\xi\right\|+\left\|\frac{1}{\sqrt{n}}\nabla [\log \pi](\frac{\xi}{\sqrt{n}}+\widehat\theta)\right\|.
    \end{aligned}
\end{equation*}
Then by the first statement of Lemma~\ref{lemma1}, Lemma~\ref{lemma2}, the twice-differentiability of $\m R(\theta)$ and Lipschitz continuity of $\m H_{\theta}$. Similar to analysis for the first statement,  we can obtain that for any $1\leq \|\xi\|\leq c\sqrt{n}$,
\begin{equation*}
    \begin{aligned}
      &\|\widetilde{\nabla}V_n(\xi)-\m \m H_{\theta^\ast}\xi\|\\
       &\leq C\,\sqrt{n}\left[d^{1+\gamma} \frac{\log n}{n}+  d^{\frac{1+\gamma_3}{2}}\big(\frac{\log n}{n}\big)^{\frac{1}{2}+\beta_1}\right]+C\,\left(d^{\frac{1+\gamma_3}{2}}\sqrt{\log n}\big(\frac{\|\xi\|}{\sqrt{n}}\big)^{\beta_1}+d^{1+\gamma}\frac{\log n}{\sqrt n}\right)\\
       &+C\,\left(d^{\gamma_2}\frac{\xi\|^2}{\sqrt{n}}+d^{\frac{1+\gamma_4}{2}+\gamma_2}\sqrt{\log n}\frac{\|\xi\|}{\sqrt{n}}\right)+C\,\sqrt{\frac{d}{n}}\\
       &\leq C\,  d^{1+\gamma}  \frac{\log n}{\sqrt{n}} + C\,d^{\frac{1+\gamma_3}{2}} \|\xi\|^{\beta_1}n^{-\frac{\beta_1}{2}}\sqrt{\log n}+C\,d^{\frac{1+\gamma_4}{2}+\gamma_2}\|\xi\|\sqrt\frac{\log  n}{n}+C\,d^{\gamma_2}\|\xi\|^2 n^{-\frac{1}{2}}.
    \end{aligned}
\end{equation*}

 \subsection{Proof of Lemma~\ref{lemmatail}}
Without loss of generality, we can assume the learning rate $\alpha=1$, as otherwise we can take $\ell(X,\theta)=\alpha\cdot \ell(X,\theta)$. Denote $K=\big\{\xi:\|\wt I^{-1/2}\xi\|\leq \mnorm{\wt I^{-\frac{1}{2}}}_{\rm op}\vee \frac{3(\sqrt{d}+t)}{\sqrt{\lambda_{\min}(\wt J)}}\big\}$. Then 
 \begin{equation*}
     \pi_n(\sqrt{n}(\theta-\wh\theta)\in K^c)=\frac{\int_{K^c}\exp(-V_n(\xi))\,\dd \xi\cdot (2\pi)^{-\frac{d}{2}}{\rm det}(\m H_{\theta^*})}{\int\exp(-V_n(\xi))\,\dd \xi\cdot (2\pi)^{-\frac{d}{2}}{\rm det}(\m H_{\theta^*})}
 \end{equation*}
 Denote $K_1=K^c\cap \{\xi\,:\,\|\xi\|\leq c_1d^{-\gamma_0-\gamma_2}\sqrt{n}\}$ and $K_2=K^c\cap \{\xi\,:\,\|\xi\|\geq c_1d^{-\gamma_0-\gamma_2}\sqrt{n}\}$. 
 When $\xi \in K_1$, we have $\|\xi\|\geq \frac{\mnorm{\wt I^{-\frac{1}{2}}}_{\rm op}}{\mnorm{\wt I^{-\frac{1}{2}}}_{\rm op}}=1$. So by Lemma~\ref{th1} and the fact that 
 \begin{equation*}
     \xi^T \m H_{\theta^*} \xi= (\wt{I}^{-\frac{1}{2}}\xi)^T \wt{I}^{\frac{1}{2}}\m H_{\theta^*}\wt{I}^{\frac{1}{2}} \wt{I}^{-\frac{1}{2}}\xi\geq \lambda_{\min}(\wt J)\|\wt{I}^{-\frac{1}{2}}\xi\|^2\geq  9(\sqrt{d}+t)^2;
 \end{equation*}
 and
  \begin{equation*}
     \xi^T \m H_{\theta^*} \xi \geq \lambda_{\min}(\m H_{\theta^*})\|\xi\|^2\geq d^{-\gamma_0}\|\xi\|^2,
 \end{equation*}
 we can verify that when  $d\leq c\frac{n^{\kappa_3}}{\log n}$
 for small enough $c$ and $K_1=K^c\cap \{\xi\,:\,\|\xi\|\leq c_1d^{-\gamma_0-\gamma_2}\sqrt{n}\}$ for a small  enough $c_1$,  it holds that
 \begin{equation*}
     V_n(\xi)\geq \frac{\xi^T\m H_{\theta^*}\xi}{4},\quad \xi \in K_1.
 \end{equation*}
 So we have
 \begin{equation*}
     \begin{aligned}
         &\int_{K_1}\exp(-V_n(\xi))\,\dd \xi\cdot (2\pi)^{-\frac{d}{2}}{\rm det}(\m H_{\theta^*})\\
         &\leq 2^{\frac{d}{2}}(2\pi)^{-\frac{d}{2}}{\rm det}\big( \frac{\m H_{\theta^*}}{2}\big)\int_{K_1} \exp\big(-\frac{\xi^T\m H_{\theta^*}\xi}{4})\,\dd \xi\\
         &\leq 2^{\frac{d}{2}}\cdot \mb P_{\chi^2(d)}(\|x\|\geq 4(\sqrt{d}+t)^2)\\
         &\leq \exp(-t^2-\frac{1}{4}),
     \end{aligned}
 \end{equation*}
 where the last inequality uses the tail inequality of $\chi^2$ distribution with $d$ degree of freedom (see for example, Lemma 1 of~\cite{10.1214/aos/1015957395}).
 
 \quad\\
 \noindent For $\xi \in K_2$ and $\theta=\wh \theta+\frac{\xi}{\sqrt{n}}$, we have 
 \begin{equation*}
     \|\hat\theta-\theta\|\geq c_1 d^{-\gamma_0-\gamma_2}.
 \end{equation*}
 Moreover, by equation~\eqref{diffhattrue} which states that $\|\widehat\theta-\theta^\ast\|\lesssim d^{\frac{1+\gamma_4}{2}}\sqrt{\frac{\log n}{n}}$, when $d\leq c\frac{n^{\kappa_3}}{\log n}$
 for small enough $c$, we have 
 \begin{equation*}
     \|\theta-\theta^*\|\geq \frac{c_1}{2} d^{-\gamma_0-\gamma_2}.
 \end{equation*}
 Therefore, by the second statement of Lemma~\ref{lemma1}, we can conclude  
 \begin{equation*} 
     \m R(\theta)-\m R(\wh \theta)= \m R(\theta)-\m R(\theta^*)+\m R(\theta^*)-\m R(\wh \theta)\geq C d^{-\gamma_0}(d^{-\gamma_1}\wedge \|\theta-\theta^*\|^2)-C_1 d^{\gamma+\frac{1+\gamma_4}{2}}\sqrt{\frac{\log n}{n}};
 \end{equation*}
 and 
 \begin{equation*} 
     |\m R_n(\theta)-\m R_n(\wh \theta)-\m R(\theta)+\m R(\wh \theta)|\leq C_2\sqrt{\frac{\log n}{n}}d^{\frac{1}{2}+\gamma}\|\theta-\wh \theta\|+C_2 \frac{\log n}{n} d^{\frac{3}{2}+\gamma}.\\
 \end{equation*}
 Then if (1) $c_1 d^{-\gamma_0-\gamma_2}\sqrt{n}\leq \|\hat\theta-\theta\|\leq d^{-\frac{\gamma_1}{2}}$, we have 
 \begin{equation*}
 \begin{aligned}
     \m R_n(\theta)-\m R_n(\wh \theta)&\geq  \m R(\theta)-\m R(\wh \theta)- |\m R_n(\theta)-\m R_n(\wh \theta)-\m R(\theta)+\m R(\wh \theta)|\\
     &\geq C\, c_1 d^{-3\gamma_0-2\gamma_2}-C_1 d^{\gamma+\frac{1+\gamma_4}{2}}\sqrt{\frac{\log n}{n}}-C_2\sqrt{\frac{\log n}{n}}d^{\frac{1}{2}+\gamma-\frac{\gamma_1}{2}}-C_2 \frac{\log n}{n} d^{\frac{3}{2}+\gamma}\\
     &\geq \frac{C\, c_1 }{2} d^{-3\gamma_0-2\gamma_2},
      \end{aligned}
 \end{equation*}
 where the last inequality uses $d\leq c\frac{n^{\kappa_3}}{\log n}$ for small enough $c$;  when (2) $\|\hat\theta-\theta\|\geq d^{-\frac{\gamma_1}{2}}$, then by $\Theta\subset [-C,C]^d$, we can get
 \begin{equation*}
 \begin{aligned}
     \m R_n(\theta)-\m R_n(\wh \theta)&\geq  \m R(\theta)-\m R(\wh \theta)- |\m R_n(\theta)-\m R_n(\wh \theta)-\m R(\theta)+\m R(\wh \theta)|\\
     &\geq C\, c_1 d^{-\gamma_1-\gamma_0}-C_1 d^{\gamma+\frac{1+\gamma_4}{2}}\sqrt{\frac{\log n}{n}}-C_2\sqrt{\frac{\log n}{n}}d^{1+\gamma}-C_2 \frac{\log n}{n} d^{\frac{3}{2}+\gamma}\\
     &\geq \frac{C\, c_1 }{2} d^{-\gamma_1-\gamma_0},
      \end{aligned}
 \end{equation*}
 where the last inequality uses $d\leq c\frac{n^{\kappa_3}}{\log n}$ for small enough $c$. So we can obtain that when $\xi \in K_2$, 
 \begin{equation*}
     V_n(\xi)=n\, \Big(\m R_n(\wh \theta+\frac{\xi}{\sqrt{n}})-\m R_n(\wh \theta)\Big)-\Big(\pi(\wh \theta+\frac{\xi}{\sqrt{n}})-\pi(\wh \theta)\Big)\geq  \frac{C\, c_1 }{4} \cdot n\cdot d^{-\gamma_0}\big(d^{-\gamma_1}\wedge d^{-2\gamma_0-2\gamma_2}\big).
 \end{equation*}
 Thus using  $d\leq c\frac{n^{\kappa_3}}{\log n}$, we have 
 \begin{equation*}
     \begin{aligned}
         &\int_{K_2}\exp(-V_n(\xi))\,\dd \xi\cdot (2\pi)^{-\frac{d}{2}}{\rm det}(\m H_{\theta^*})\\
         &\leq\exp\Big(-\frac{d}{2}\log (2\pi)+\frac{d}{2}\log \big(\mnorm{\m H_{\theta^*}}_{\rm op}\big)\Big)\cdot \exp\Big( \frac{C\, c_1 }{4} \cdot n\cdot d^{-\gamma_0}\big(d^{-\gamma_1}\wedge d^{-2\gamma_0-2\gamma_2}\big)\Big)\\
         &\leq \exp\Big( \frac{C\, c_1 }{8} \cdot n\cdot d^{-\gamma_0}\big(d^{-\gamma_1}\wedge d^{-2\gamma_0-2\gamma_2}\big)\Big).\\
     \end{aligned}
 \end{equation*}
 It remains to bound the denominator $\int\exp(-V_n(\xi))\,\dd \xi\cdot (2\pi)^{-\frac{d}{2}}{\rm det}(\m H_{\theta^*})$, we have
 \begin{equation*}
     \begin{aligned}
         &\int\exp(-V_n(\xi))\,\dd \xi\cdot (2\pi)^{-\frac{d}{2}}{\rm det}(\m H_{\theta^*})\\
         &\geq   (2\pi)^{-\frac{d}{2}}{\rm det}(\m H_{\theta^*})\int_{\|\xi\|\leq 4\sqrt{d/\lambda_{\min}(\m H_{\theta^*})}}\exp\big(-\frac{\xi^T\m H_{\theta^*}\xi}{2}\big)\,\dd \xi\\
         &\qquad\cdot \underset{\|\xi\|\leq 4\sqrt{d/\lambda_{\min}(\m H_{\theta^*})}}{\sup} \exp\big(\frac{\xi^T\m H_{\theta^*}\xi}{2}-V_n(\xi)\big)\\
         &\geq \exp(-\frac{1}{4}),
     \end{aligned}
 \end{equation*}
 where the last inequality uses $\lambda_{\min}(\m H_{\theta^*})\geq C\, d^{-\gamma_0}$, $d\leq c\frac{n^{\kappa_3}}{\log n}$ and the statements of Lemma~\ref{th1}. We can then obtain the desired results by combining all pieces.
 \subsection{Proof of Lemma~\ref{lemma1}}
 We first prove the first statement.  It's equivalent to show that it holds with probability larger than $1-\frac{1}{3n^2}$ that 
 for any $\theta,\theta'\in B_r(\theta^\ast)$ and $v\in\mb S^{d-1}$,
 \begin{equation*}
 \begin{aligned}
  &\Big|\frac{1}{n}\sum_{i=1}^nv^Tg(X_i,\theta)-\frac{1}{n}\sum_{i=1}^nv^Tg(X_i,\theta')-\mathbb{E}[v^Tg(X,\theta)]+\mathbb{E}[v^Tg(X,\theta')]\Big|\\
  &\leq c\Big(\sqrt{\frac{\log n}{n}}d^{\frac{1+\gamma_3}{2}}\|\theta-\theta'\|^{\beta_1}+\frac{\log n}{n}d^{1+\gamma}\Big).
   \end{aligned}
 \end{equation*}
Consider a minimal $\frac{3}{n}$-covering set $\m A$ of $\mb S^{d-1}$ such that $\m A\subset\mb S^{d-1}$, then $\log |\m A|\leq d\log n$. For any $v\in \m A$, define the function class
 \begin{equation*}
     \begin{aligned}
       \m G_v=\{d^{-\gamma}(v^Tg(\cdot,\theta)-v^Tg(\cdot,\theta')): \theta,\theta'\in B_r(\theta^\ast)\}.
     \end{aligned}
 \end{equation*}
 Let $\overline{\m G}_v=\{af: a\in [0,1], f\in \m G_v\}$ be the star hull of $\m G_v$. Then since $ {\sup}_{x\in \m X,\theta\in B_r(\theta^\ast)}\|g(x,\theta)\|\leq C\,d^{\gamma}$, it holds that $\sup_{f\in \overline{\m G}_v,x\in \m X}|f(x)|\leq 2C $. Consider the local Rademacher complexity associated with $\overline{\m G}_v$,
 \begin{equation*}
     \overline{R}_n(\delta;\overline{\m G}_v)=\mathbb{E}_{X^{(n)}}\mathbb{E}_{\varepsilon}\Bigg[\underset{f\in \overline{\m G}_v\atop \mb E f^2\leq \delta^2}{\sup}\bigg|\frac{1}{n}\sum_{i=1}^n \varepsilon_if(X_i)\bigg|\Bigg],
 \end{equation*}
 where $\varepsilon_i$ are i.i.d. samples from Rademacher distribution, i.e., $\mathbb{P}(\varepsilon_i=1)=\mathbb{P}(\varepsilon_i=-1)=0.5$. We will use the following uniform law, which is a special case of Theorem 14.20 of~\cite{wainwright_2019}, to prove the  desired result.
 \begin{lemma}\label{uniformlaw}~(\cite{wainwright_2019}, Theorem 14.20)
 Given a uniformly 1-bounded function class $\m F$ that is star shaped around $0$, let $(\delta^\ast)^2\geq \frac{c}{n}$ be any solution to the inequality  $\overline{R}_n(\delta;\m F)\leq \delta^2$, then we have 
 \begin{equation*}
     \underset{f\in \m F}{\sup}\frac{\Big|\frac{1}{n}\sum_{i=1}^n f(X_i)-\mb E [f(X)]\Big|}{\sqrt{\mb E[f(X)^2]}+\delta^\ast}\leq 10\delta^\ast
 \end{equation*}
 with probability greater than $1-c_1\exp(-c_2\, n\cdot(\delta^*)^2)$.
 \end{lemma}
 Next we will use Dudley's inequality (see for example, Theorem 5.22 of~\cite{wainwright_2019}) to determine the critical radius $\delta^\ast$ in Lemma~\ref{uniformlaw}. For $f,f':\m X\to \mb R$ , define the pseudometric
 \begin{equation*}
     d_n(f,f')=\sqrt{\frac{1}{n}\sum_{i=1}^n (f(X_i)-f'(X_i))^2}.
 \end{equation*}
 Then by uniformly boundness of functions in class $\overline{\m G}_v$ , we can obtain that 
 \begin{equation*}
     \begin{aligned}
       &\log \mathbf N(\overline{\m G}_v,d_n,\varepsilon)\\
       &\leq \log \frac{4C}{\varepsilon}+\log \mathbf N(\m G_v, d_n,\frac{\varepsilon}{2})\\
       &\leq \log \frac{4C}{\varepsilon}+\log \mathbf N(\m B_r(\theta^\ast), d_n^{g},\frac{d^{\gamma}\varepsilon}{2})\\
       &\leq C_1\, d\log\frac{n}{\varepsilon},
     \end{aligned}
 \end{equation*}
 where  recall that $\mathbf N(\mathcal{F},d_n,\varepsilon)$ denote the $\varepsilon$-covering number of class $\m F$ w.r.t pseudo-metric $d_n$. Let 
 \begin{equation*}
     \begin{aligned}
       r_n^2&=\underset{f,f'\in \overline{\m G}_v\atop \mb E [f^2], \mathbb{E}[{f'}^2]\leq \delta^2}{\sup} d_n^2(f,f')\\
       &\leq 4 \underset{f\in \m G_v\atop \mb E[f^2]\leq \delta^2}{\sup}\frac{1}{n}\sum_{i=1}^n f^2(X_i)\\
       &\leq 8 \underset{f\in \m G_v\atop \mb E[f^2]\leq \delta^2}{\sup}\frac{1}{n}\sum_{i=1}^n (f(X_i)-\mb E f(X))^2+8\delta^2.
     \end{aligned}
 \end{equation*}
 Then by (3.84) of~\cite{wainwright_2019}, we can obtain that $\mb E [r_n^2]\leq C\, \delta^2+C\,\m R_n(\delta)$. Choose $\delta^\ast=c\, d^{\frac{1}{2}}\sqrt\frac{\log n}{n}$, then by Dudley's inequality, 
\begin{equation*}
\begin{aligned}
    \overline{\m R}_n(\delta^\ast)&\leq C\, \frac{1}{\sqrt{n}}\mathbb{E}\int_{0}^{r_n} d^{\frac{1}{2}}\sqrt{\log \frac{n}{\varepsilon}}d\varepsilon\\
    &=C\,\frac{1}{\sqrt{n}}\mathbb{E}\int_{0}^{1} r_nd^{\frac{1}{2}}\sqrt{\log \frac{n}{\varepsilon r_n}}d\varepsilon\\
    &=C\,\mathbb{E}\left[\frac{1}{\sqrt{n}}\int_{0}^{1} r_nd^{\frac{1}{2}}\sqrt{\log \frac{n}{\varepsilon r_n}}d\varepsilon\cdot \mathbf{1}(r_n<n^{-\frac{1}{2}})\right]+C\,\mathbb{E}\left[\frac{1}{\sqrt{n}}\int_{0}^{1} r_nd^{\frac{1}{2}}\sqrt{\log \frac{n}{\varepsilon r_n}}d\varepsilon\cdot \mathbf{1}(r_n>n^{-\frac{1}{2}})\right]\\
    &\leq C\, d^{\frac{1}{2}}\frac{\sqrt{\log n}}{n}+C\,\mathbb{E}\bigg[\frac{1}{\sqrt n}\int_{0}^1 r_nd^{\frac{1}{2}}\sqrt{\log \frac{n^{\frac{3}{2}}}{\varepsilon}}d\varepsilon\bigg]\\
    &\leq C_1\,  \sqrt{\frac{\log n}{n}}d^{\frac{1}{2}}\sqrt{{\delta^\ast}^2+\overline{\m R}_n(\delta^\ast)}.
    \end{aligned}
\end{equation*}
 Then if $\overline{\m R}_n(\delta^\ast)>(\delta^\ast)^2$, we can obtain that $\overline{\m R}_n(\delta^\ast)\leq 2C_1^2\, d\frac{\log n}{n}\leq 2C_1^2c^{-2}\, {\delta^\ast}^2$. thus when $c$ is large enough, $\delta^*$ solves the inequality $\overline{\m R}_n(\delta^\ast)\leq(\delta^\ast)^2$. Then by Lemma~\ref{uniformlaw} and the assumption that $\sup_{v\in\mb S^{d-1}}\mb E\big[(v^Tg(X,\theta)-v^Tg(X,\theta'))\big]^2\leq C\, d^{\gamma_3}\|\theta-\theta'\|^{2\beta_1}$, there exists a constant $C$ such that it holds with probability larger than $1-\exp(-4d\log n)$ that for any  $\theta,\theta'\in B_r(\theta^\ast)$, 
 \begin{equation*}
 \begin{aligned}
    &\Big|\frac{1}{n}\sum_{i=1}^nv^Tg(X_i,\theta)-\frac{1}{n}\sum_{i=1}^nv^Tg(X_i,\theta')-\mathbb{E}v^Tg(X,\theta)+\mathbb{E}v^Tg(X,\theta')\Big|\\
    &\leq C\,\Big(\sqrt{\frac{\log n}{n}}d^{\frac{1+\gamma_3}{2}}\|\theta-\theta'\|^{\beta_1}+\frac{\log n}{n}d^{1+\gamma}\Big).
     \end{aligned}
  \end{equation*}
 By the fact that  $\log |\m A|\leq d\log n$, it holds with probability larger than $1-\exp(-3d\log n)$ that for any $v\in \m A$ and $\theta,\theta'\in B_r(\theta^\ast)$, 
 \begin{equation*}
 \begin{aligned}
    &\Big|\frac{1}{n}\sum_{i=1}^nv^Tg(X_i,\theta)-\frac{1}{n}\sum_{i=1}^nv^Tg(X_i,\theta')-\mathbb{E}v^Tg(X,\theta)+\mathbb{E}v^Tg(X,\theta')\Big|\\
    &\leq C\,\Big(\sqrt{\frac{\log n}{n}}d^{\frac{1+\gamma_3}{2}}\|\theta-\theta'\|^{\beta_1}+\frac{\log n}{n}d^{1+\gamma}\Big).
     \end{aligned}
  \end{equation*}
  Moreover, for any $\widetilde v\in\mb S^{d-1}$, there exists $v\in \m A$ so that $\|v-\widetilde v\|\leq \frac{3}{n}$,  hence for any $\theta,\theta'\in B_r(\theta^\ast)$,
  \begin{equation*}
  \begin{aligned}
    &\underset{v\in\mb S^{d-1}}{\sup}\Big|\frac{1}{n}\sum_{i=1}^n\widetilde v^Tg(X_i,\theta)-\frac{1}{n}\sum_{i=1}^n\widetilde v^Tg(X_i,\theta')-\mathbb{E}\widetilde v^Tg(X,\theta)+\mathbb{E}\widetilde v^Tg(X,\theta')\Big|\\
    &=\underset{v\in \m A}{\sup}\Big|\frac{1}{n}\sum_{i=1}^nv^Tg(X_i,\theta)-\frac{1}{n}\sum_{i=1}^nv^Tg(X_i,\theta')-\mathbb{E}v^Tg(X,\theta)+\mathbb{E}v^Tg(X,\theta')\Big|+\m O(\frac{d}{\sqrt{n}}).
      \end{aligned}
  \end{equation*}
 Then, it follows that it holds with probability larger than $1-\exp(3d\log n)\geq 1-\frac{1}{3n^2}$ that
 \begin{equation*}
 \begin{aligned}
  &\bigg\|\frac{1}{n}\sum_{i=1}^n g(X_i,\theta)-\frac{1}{n}\sum_{i=1}^n g(X_i,\theta')-\mathbb{E} g(X,\theta)+\mathbb{E}g(X,\theta')\bigg\|\\
    &=\underset{v\in\mb S^{d-1}}{\sup}\bigg|\frac{1}{n}\sum_{i=1}^nv^Tg(X_i,\theta)-\frac{1}{n}\sum_{i=1}^nv^Tg(X_i,\theta')-\mathbb{E}v^Tg(X,\theta)+\mathbb{E}v^Tg(X,\theta')\bigg|\\
    &\leq C\,\bigg(\sqrt{\frac{\log n}{n}}d^{\frac{1+\gamma_3}{2}}\|\theta-\theta'\|^{\beta_1}+\frac{\log n}{n}d^{1+\gamma}\bigg).
     \end{aligned}
  \end{equation*}
The proof of the first statement is then completed.  For the second statement, by the assumption that for any $\theta,\theta'\in \Theta$ and $x\in \m X$, $|\ell(X,\theta)-\ell(X,\theta')|\leq Cd^{\gamma}\|\theta-\theta'\|$, we can obtain that for any $\theta,\theta'\in \Theta$
\begin{equation*}
    \mb E \big[(\ell(X,\theta)-\ell(X,\theta'))^2\big]\leq C^2 d^{2\gamma} \|\theta-\theta'\|^2,
\end{equation*}
and 
\begin{equation*}
    \underset{x\in \m X}{\sup}|\ell(X,\theta)-\ell(X,\theta')|\leq Cd^{\gamma}(\|\theta\|+\|\theta'\|)\leq C_1\, d^{\frac{1}{2}+\gamma}.
\end{equation*}
We can therefore prove the second statement  using the same strategy as the first statement. For the third statement, define $\delta_n=(\frac{\log n}{n}d^{-\frac{3}{2}})\wedge((\frac{\log n}{n})^{\frac{3}{2}}d^{-\frac{1+\gamma_3}{2}})^{\frac{1}{1+\beta_1}}$. For $k=0,1,\cdots,\lfloor \log_2 \frac{2r}{\delta_n}\rfloor+1$, we define the set
\begin{equation*}
    \m A_k=\left\{\begin{array}{ll}
       \{\theta,\theta'\in B_r(\theta^\ast): \|\theta-\theta'\|\leq \delta_n\}  &  k=0; \\
         \{\theta,\theta'\in B_r(\theta^\ast): 2^{k-1}\delta_n<\|\theta-\theta'\|\leq 2^{k}\delta_n\}  &  k=1,2,\cdots  \lfloor\log_2 \frac{2r}{\delta_n}\rfloor;\\
           \{\theta,\theta'\in B_r(\theta^\ast): 2^{k-1}\delta_n<\|\theta-\theta'\|\leq 2r\}  &  k= \lfloor\log_2 \frac{2r}{\delta_n}\rfloor+1.\\
    \end{array}
    \right.
\end{equation*}
Then $\{\theta,\theta'\in B_r(\theta^\ast)\}=\sum_{k=1}^{\log_2 \lfloor\frac{2r}{\delta_n}\rfloor+1}\m A_k$. Fix an integer $0\leq k\leq \lfloor\log_2 \frac{2r}{\delta_n}\rfloor+1$, we consider the function set 
\begin{equation*}
    \m L_k=\Big\{\frac{1}{2^k\delta_n}d^{-\gamma}(\ell(\cdot,\theta)-\ell(\cdot,\theta')-g(\cdot,\theta')(\theta-\theta')):(\theta,\theta')\in \m A_k\Big\}.
\end{equation*}
Then  there exists a constant $c$ such that for any $f\in \m L_k$, it holds that $\sup_{x\in \m X} |f(x)|\leq c$ and $\mb E [f^2(X)]\leq c\frac{1}{2^{2k}\delta_n^2}d^{-2\gamma}d^{\gamma_3}(2^{k}\delta_n)^{2+2\beta_1}\leq c\,d^{\gamma_3-2\gamma}(2^k\delta_n)^{2\beta_1}\leq 4c\,d^{\gamma_3-2\gamma}(2^{k-1}\delta_n)^{2\beta_1}$. Then consider the star hull $\overline{\m L}_k$ of $\m L_k$, by (1) $d\lesssim n^{\kappa_2}$; (2) the Lipschitzness of $\ell$; (3) the  bound on the $\varepsilon$-covering number of $B_r(\theta^\ast)$w.r.t $d_n^{g}$, it holds that 
\begin{equation*}
\begin{aligned}
     &\log \mathbf N(\overline{\m L}_k,d_n,\varepsilon)\\
     &\leq \log \frac{2c}{\varepsilon}+\log \mathbf N(\m L_k,d_n,\varepsilon)\\
     &\leq C\,  d\log \frac{n}{\varepsilon}.
\end{aligned}
  \end{equation*}
  Then similar as the proof of the first statement, we can use Dudley's inequality and Lemma~\ref{uniformlaw} to obtain that there exists a constant $c$ such that it holds with probability at least $1-\frac{1}{3n^3}$ that for any $(\theta,\theta')\in \m A_k$, 
  \begin{equation*}
  \begin{aligned}
      &\bigg|\frac{1}{n}\sum_{i=1}^n\ell(X_i,\theta)-\frac{1}{n}\sum_{i=1}^n\ell(X_i,\theta')-\frac{1}{n}\sum_{i=1}^ng(X_i,\theta')(\theta-\theta')\\
     & -\Big(\mathbb{E}\ell(X,\theta)-\mathbb{E}\ell(X,\theta')-\mb Eg(X,\theta')(\theta-\theta')\Big)\bigg|\\
      &\leq C\,\bigg(\sqrt{\frac{\log n}{n}}d^{\frac{1+\gamma_3}{2}}\cdot (2^{k-1}\delta_n)^{\beta_1+1}+\frac{\log n}{n}d^{1+\gamma}\cdot(2^{k-1}\delta_n)\bigg)\\
      &\leq C\,\bigg(\sqrt{\frac{\log n}{n}}d^{\frac{1+\gamma_3}{2}} \cdot(\|\theta-\theta'\|+\delta_n)^{\beta_1+1}+\frac{\log n}{n}d^{1+\gamma}\cdot(\|\theta-\theta'\|+\delta_n)\bigg)\\
      &\leq 4C\,\bigg(\sqrt{\frac{\log n}{n}}d^{\frac{1+\gamma_3}{2}}\|\theta-\theta'\|^{\beta_1+1}+\frac{\log n}{n}d^{1+\gamma}\|\theta-\theta'\|+(\frac{\log n}{n})^2\bigg).
      \end{aligned}
  \end{equation*}
  Then by $\log_2\frac{r}{\delta_n}\lesssim \log n$, consider the intersection of the above events for $k=0,1,\cdots,\lfloor\log_2 \frac{r}{\delta_n}\rfloor+1$, we can obtain the desired result.
  
  \subsection{Proof of Lemma~\ref{lemmathetadiamond}}
 Recall $\widehat\theta^{\diamond}=\theta^\ast-n^{-1}\sum_{i=1}^n\m H_{\theta^\ast}^{-1} g(X_i,\theta^\ast)$, then by $\mathbb{E} [g(X,\theta^\ast)=\nabla \m R(\theta^*)=0$,  we have 
     \begin{equation*}
     \begin{aligned}
          \|\widehat{\theta}^{\diamond}-\theta^\ast\|&=\|\frac{1}{n}\sum_{i=1}^n\m H_{\theta^\ast}^{-1} g(X_i,\theta^\ast)-\mb E[\m H_{\theta^\ast}^{-1} g(X,\theta^\ast)]\|\\
          &=\underset{v\in\mb S^{d-1}}{\sup}\left|\frac{1}{n}\sum_{i=1}^n v^T\m H_{\theta^\ast}^{-1} g(X_i,\theta^\ast)-\mb E[v^T \m H_{\theta^\ast}^{-1} g(X,\theta^\ast)]\right|.
     \end{aligned}
          \end{equation*}
 It remains to derive a high probability bound of the supremum of the above empirical process. Consider a minimal $\frac{3}{n}$-covering set $\m A$ of $S^{d-1}$ such that $A\subset\mb S^{d-1}$, then $\log |\m A|\leq d\log n$. Fix an arbitrary $v\in\mb S^{d-1}$, then by the assumption that (1) $\m H_{\theta^\ast}^{-1}\mb E[g(X_i,\theta^\ast)^Tg(X_i,\theta^\ast)]\m H_{\theta^\ast}^{-1}\preceq Cd^{\gamma_4}I_d$; (2)  for any $\theta\in \Theta$, $\m R(\theta)-\m R(\theta^\ast)\geq C'd^{-\gamma_0}(d^ {-\gamma_1}\wedge\|\theta-\theta^\ast\|^2)$, which leads to $\m H_{\theta^\ast}\succeq C'd^{-\gamma_0}I_d$; (3) $\sup_{x\in \m X}\|g(X,\theta^\ast)\|\leq Cd^{\gamma}$, we can obtain 
 \begin{equation*}
     \underset{X\in \m X\atop v\in \mb S^{d-1}}{\sup} |v^T\m H_{\theta^\ast}^{-1} g(X,\theta^\ast)|\leq CC'd^{\gamma+\gamma_0},
 \end{equation*}
  and 
  \begin{equation*}
     \underset{v\in \mb S^{d-1}}{\sup} \mb E[v^T\m H_{\theta^\ast}^{-1} g(X,\theta^\ast)]^2\leq Cd^{\gamma_4}.
 \end{equation*}
  Therefore using Bernstein-type bound (see for example, Proposition 2.10 of~\cite{wainwright_2019}), we can get there exists a constant $c$ such that it holds with probability larger than $1-\exp(3d\log n)$ that,
  \begin{equation*}
      \left|\frac{1}{n}\sum_{i=1}^n v^T\m H_{\theta^\ast}^{-1} g(X_i,\theta^\ast)-\mb Ev^T \m H_{\theta^\ast}^{-1} g(X,\theta^\ast)\right|\leq C\,(d^{\frac{1+\gamma_4}{2}}\sqrt{\frac{\log n}{n}}+d^{1+\gamma+\gamma_0}\frac{\log n}{n}).
  \end{equation*}
Moreover, for any $\widetilde v\in\mb S^{d-1}$, there exists $v\in \m A$ so that $\|v-\widetilde v\|\leq \frac{3}{n}$,  hence for any $\theta,\theta'\in B_r(\theta^\ast)$,
 \begin{equation*}
 \begin{aligned}
     \left|\frac{1}{n}\sum_{i=1}^n \widetilde v^T\m H_{\theta^\ast}^{-1} g(X_i,\theta^\ast)-\mb E\widetilde v^T \m H_{\theta^\ast}^{-1} g(X,\theta^\ast)\right|&\leq  \left|\frac{1}{n}\sum_{i=1}^n v^T\m H_{\theta^\ast}^{-1} g(X_i,\theta^\ast)-\mb Ev^T \m H_{\theta^\ast}^{-1} g(X,\theta^\ast)\right|\\
     &+\m O(d^{\gamma_0+\gamma}\frac{\log n}{n}).
      \end{aligned}
 \end{equation*}
 Thus by a simple union bound, it holds with probability larger than $1-\exp(2d\log n)>1-\frac{1}{n^2}$ that 
  \begin{equation*}
      \underset{v\in \mb S^{d-1}}{\sup}\left|\frac{1}{n}\sum_{i=1}^n v^T\m H_{\theta^\ast}^{-1} g(X_i,\theta^\ast)-\mb Ev^T \m H_{\theta^\ast}^{-1} g(X,\theta^\ast)\right|\leq 2C\,(d^{\frac{1+\gamma_4}{2}}\sqrt{\frac{\log n}{n}}+d^{1+\gamma+\gamma_0}\frac{\log n}{n}).
  \end{equation*}
 We can thus obtain that it holds with probability larger than $1-n^{-2}$ that 
     \begin{equation*}
     \begin{aligned}
          \|\widehat{\theta}^{\diamond}-\theta^\ast\|\leq C\, d^{\frac{1+\gamma_4}{2}}\sqrt\frac{\log n}{n}+C\,d^{1+\gamma+\gamma_0} \frac{\log n}{n}.
     \end{aligned}
          \end{equation*}

  \subsection{Proof of Lemma~\ref{lemma2}}
  Firstly by $\m R_n(\widehat{\theta})\leq \m R_n(\theta^\ast)$ and $\m R(\theta)-\m R(\theta^\ast)\geq C' d^{-\gamma_0}(d^ {-\gamma_1}\wedge \|{\theta}-\theta^\ast\|^2)$, we can obtain that 
  \begin{equation*}
      C' d^{-\gamma_0}(d^ {-\gamma_1}\wedge \|\widehat{\theta}-\theta^\ast\|^2)\leq \m R(\widehat\theta)-\m R(\theta^\ast)\leq \m R(\widehat\theta)-\m R(\theta^\ast)-\m R_n(\widehat{\theta})+\m R_n(\theta^\ast).
  \end{equation*}
 It follows from the second statement of Lemma~\ref{lemma1} that 
  \begin{equation*}
      d^{-\gamma_0}(d^ {-\gamma_1}\wedge \|\widehat{\theta}-\theta^\ast\|^2)\leq C\,\sqrt\frac{\log n}{n} d^{\frac{1}{2}+\gamma}\|\widehat{\theta}-\theta^\ast\|+C\,\frac{\log n}{n} d^{\frac{3}{2}+\gamma}.
  \end{equation*}
 If $\|\widehat\theta-\theta^\ast\|\geq d^{-\frac {\gamma_1}{2}}$, then  
  \begin{equation*}
      d^{-\gamma_0-\gamma_1}\leq C\,\sqrt\frac{\log n}{n} d^{\frac{1}{2}+\gamma}\|\widehat{\theta}-\theta^\ast\|+C\,\frac{\log n}{n} d^{\frac{3}{2}+\gamma}.
  \end{equation*}
 On the other hand,  as $\widehat\theta\in \Theta\subseteq [-C,C]^d$, we have $\|\widehat\theta-\theta^\ast\|\leq 2C\sqrt{d}$, we can then obtain that when $d\leq c(\frac{n}{\log n})^{\frac{1}{2+2(\gamma+\gamma_0+\gamma_1)}}$,
  \begin{equation*}
  \begin{aligned}
      \sqrt\frac{\log n}{n} d^{\frac{1}{2}+\gamma}\|\widehat{\theta}-\theta^\ast\|+\frac{\log n}{n} d^{\frac{3}{2}+\gamma}&\leq 2Cd^{1+\gamma}\sqrt\frac{\log n}{n}+\frac{\log n}{n} d^{\frac{3}{2}+\gamma}\\
      &\leq 2C\sqrt{c}\,d^{-\gamma_0-\gamma_1}+c\,d^{-\frac{1}{2}-\gamma-2(\gamma_0+\gamma_1)},
        \end{aligned}
  \end{equation*}
which will cause contradiction  when $c$ is sufficiently small. Hence  we have $\|\widehat\theta-\theta^\ast\|< d^{-\frac {\gamma_1}{2}}$ and thus
   \begin{equation*}
      d^{-\gamma_0}  \|\widehat{\theta}-\theta^\ast\|^2\leq C\, \sqrt\frac{\log n}{n} d^{\frac{1}{2}+\gamma}\|\widehat{\theta}-\theta^\ast\|+C\,\frac{\log n}{n} d^{\frac{3}{2}+\gamma},
  \end{equation*}
 which leads to $\|\widehat{\theta}-\theta^\ast\|\leq C_1\, \sqrt\frac{\log n}{n} d^{\frac{1}{2}+\gamma+\gamma_0}$. We will first show the first statement of Lemma~\ref{lemma2} and use the statement to improve the dependence of $d$ in the bound of $\sqrt\frac{\log n}{n} d^{\frac{1}{2}+\gamma+\gamma_0}$.\\

\noindent By $\m R_n(\widehat{\theta})\leq \m R_n(\widetilde{\theta})$ for any $\widetilde{\theta}\in B_r(\theta^\ast)$, we can obtain that 
\begin{equation*}
\begin{aligned}
    &-\frac{1}{n}\sum_{i=1}^n g(X_i,\widehat{\theta})(\widetilde{\theta}-\widehat{\theta})\\
   & \leq \m R_n(\widetilde{\theta})-\m R_n(\widehat\theta)-\frac{1}{n}\sum_{i=1}^n g(X_i,\widehat{\theta})(\widetilde{\theta}-\widehat{\theta})\\
&\leq \left|\m R_n(\widetilde \theta)-\m R_n(\widehat\theta)-\frac{1}{n}\sum_{i=1}^n g(X_i,\widehat\theta)(\widetilde \theta-\widehat\theta)-\m R(\widetilde \theta)+\m R(\widehat\theta)+\mb E[g(X,\widehat\theta )(\widetilde \theta-\widehat\theta )]\right|\\
    &+\left| \m R(\widetilde \theta)-\m R(\widehat\theta)-\mb E [g(X,\widehat\theta )(\widetilde \theta-\widehat\theta )]\right|.
    \end{aligned}
\end{equation*}

The first term can be bounded using the third statement of Lemma~\ref{lemma1}, that is 
\begin{equation*}
    \begin{aligned}
         &\left|\m R_n(\widetilde \theta)-\m R_n(\widehat\theta)-\frac{1}{n}\sum_{i=1}^n g(X_i,\widehat\theta)(\widetilde \theta-\widehat\theta)-\m R(\widetilde \theta)+\m R(\widehat\theta)+\mb E [g(X,\widehat\theta )(\widetilde \theta-\widehat\theta )]\right|\\
         &\leq C\, \sqrt{\frac{\log n}{n}}d^{\frac{1+\gamma_3}{2}}\|\widehat\theta-\widetilde{\theta}\|^{\beta_1+1}+C\,\frac{\log n}{n}d^{1+\gamma}\|\widehat\theta-\widetilde\theta\|+C\,(\frac{\log n}{n})^2.
    \end{aligned}
\end{equation*}

The second term can be bounded using the twice differentiability of $\m R$ around $\theta^\ast$, 
\begin{equation*}
    \begin{aligned}
         \left|\m R(\widetilde\theta)-\m R(\widehat\theta)-\mb E [g(X,\widehat\theta )(\widetilde\theta-\widehat\theta)]\right|\leq \frac{1}{2}\underset{c\in[0,1]}{\sup}\mnorm{\m H_{c\widetilde{\theta}+(1-c)\widehat\theta}}_{  \rm  op}\|\widehat\theta-\widetilde\theta\|^2\leq C\,d\|\widehat\theta-\widetilde\theta\|^2.
    \end{aligned}
\end{equation*}
where the last inequality is due to the assumption that the mixed partial derivatives of $\m R(\theta)$ up to order two are uniformly bounded by an $(n,d)$-independent constant on $\m B_r(\theta^*)$. Then we choose $\widetilde\theta=\widehat{\theta}-t\frac{\sum_{i=1}^n g(X_i,\widehat\theta)}{\|\sum_{i=1}^n g(X_i,\widehat\theta)\|}$ for a $t>0$ that will be chosen later.
Thus
\begin{equation*}
    \begin{aligned}
        & C_1\,t\left\|\frac{1}{n}\sum_{i=1}^n g(X_i,\widehat{\theta})\right\|\leq \sqrt{\frac{\log n}{n}}d^{\frac{1+\gamma_3}{2}}t^{\beta_1+1}+\frac{\log n}{n}d^{1+\gamma}t+(\frac{\log n}{n})^2+dt^2\\
        &\Rightarrow C_1\,\left\|\frac{1}{n}\sum_{i=1}^n g(X_i,\widehat{\theta})\right\|\leq \sqrt{\frac{\log n}{n}}d^{\frac{1+\gamma_3}{2}}t^{\beta_1}+\frac{\log n}{n}d^{1+\gamma}+(\frac{\log n}{n})^2/t+dt.
    \end{aligned}
\end{equation*}
Choose $t=\frac{\log n}{n}$, we have it holds with probability at least $1-n^{-2}$ that
\begin{equation*}
     \left\|\frac{1}{n}\sum_{i=1}^n g(X_i,\widehat{\theta})\right\|\leq C\, d^{1+\gamma} \frac{\log n}{n}+ C\, d^{\frac{1+\gamma_3}{2}}\big(\frac{\log n}{n}\big)^{\frac{1}{2}+\beta_1}.
\end{equation*}
For the second statement,  recall $\widehat\theta^\diamond=\theta^\ast-\frac{1}{n}\sum_{i=1}^n \m H_{\theta^\ast}^{-1}g(X_i,\theta)$. By Lemma~\ref{lemmathetadiamond} and the assumption that $d\leq c(\frac{n}{\log n})^{\frac{1}{2+2(\gamma+\gamma_0+\gamma_1)}}$, we can obtain $\|\widehat\theta^\diamond-\theta^\ast\|\leq C\, d^{\frac{1+\gamma_4}{2}}\sqrt\frac{\log n}{n}$.
We claim that it suffices to show that 
\begin{equation}\label{normgdiamond}
    \left\|\frac{1}{n}\sum_{i=1}^ng(X_i,\widehat{\theta}^{\diamond})\right\|\leq C\, d^{\frac{1+\gamma_3}{2}+\beta_1(\frac{1+\gamma_4}{2})}(\frac{\log n}{n})^{\frac{1+\beta_1}{2}}+C\, d^{1+\gamma\vee(\gamma_2+\gamma_4)}\frac{\log n}{n}
\end{equation}
holds with probability at least $1-cn^{-2}$. Indeed, under the above statement, we have
\begin{equation*}
 \begin{aligned}
 & \|\mathbb{E}[g(X,\widehat{\theta})]-\mathbb{E}[g(X,\widehat\theta^{\diamond})]\|\\
  &\leq \bigg\|\frac{1}{n}\sum_{i=1}^ng(X_i,\widehat\theta)-\frac{1}{n}\sum_{i=1}^ng(X_i,\widehat{\theta}^{\diamond})-\mathbb{E}[g(X,\widehat\theta)]+\mathbb{E}[g(X,\widehat{\theta}^{\diamond})]\bigg\|\\
     &\quad+\bigg\|\frac{1}{n}\sum_{i=1}^ng(X_i,\widehat\theta)\bigg\|+\bigg\|\frac{1}{n}\sum_{i=1}^ng(X_i,\widehat{\theta}^{\diamond})\bigg\|\\
     &\leq C\,\bigg( \sqrt\frac{\log n}{n}d^{\frac{1+\gamma_3}{2}}\|\widehat\theta-\widehat{\theta}^{\diamond}\|^{\beta_1}+ d^{\frac{1+\gamma_3}{2}+\beta_1(\frac{1+\gamma_4}{2})}(\frac{\log n}{n})^{\frac{1+\beta_1}{2}}+ d^{1+\gamma\vee(\gamma_2+\gamma_4)}\frac{\log n}{n}\bigg),
      \end{aligned}
 \end{equation*}
where the last inequality follows from the first statement of Lemma~\ref{lemma1}. On the other hand, by the Lipschitzness of $\m H_{\theta}$ around $\theta^\ast$, we can obtain that,
\begin{equation*}
    \begin{aligned}
     &\|\mathbb{E}[g(X,\widehat{\theta})]-\mathbb{E}[g(X,\widehat\theta^{\diamond})]\|\\
     &\geq  \|\m H_{\theta^\ast}(\widehat{\theta}-\widehat\theta^\diamond)\|- \|\mathbb{E}[g(X,\widehat{\theta})]-\mathbb{E}[g(X,\widehat\theta^{\diamond})]-H_{\theta^\ast}(\widehat{\theta}-\widehat\theta^\diamond)\| \\
     &=  \|\m H_{\theta^\ast}(\widehat{\theta}-\widehat\theta^\diamond)\|-  \underset{v\in \mb S^{d-1}}{\sup}\big|\mathbb{E}[v^Tg(X,\widehat{\theta})]-\mathbb{E}[v^Tg(X,\widehat\theta^{\diamond})]-v^TH_{\theta^\ast}(\widehat{\theta}-\widehat\theta^\diamond)\big|  \\
     &\geq  \rho_1(\m H_{\theta^\ast})\|\widehat\theta-\widehat\theta^\diamond\|-\underset{v\in \mb S^{d-1}}{\sup}\underset{t\in (0,1)}{\sup} \big|v^T(\m H_{t\widehat{\theta}^\diamond+(1-t)\widehat\theta}-\m H_{\theta^\ast})(\widehat\theta-\widehat\theta^\diamond)\big|\\
     &\geq \rho_1(\m H_{\theta^\ast})\|\widehat\theta-\widehat\theta^\diamond\|-C\,\left(d^{\frac{1}{2}+\gamma+\gamma_2+\gamma_0}\sqrt{\frac{\log n}{n}}\|\widehat\theta-\widehat\theta^\diamond\|\right),
    \end{aligned}
\end{equation*}
 where the last inequality uses $\|\widehat{\theta}-\theta^\ast\|\leq C_1\, \sqrt\frac{\log n}{n} d^{\frac{1}{2}+\gamma+\gamma_0}$ and $\|\widehat\theta^\diamond-\theta^\ast\|\leq C\, d^{\frac{1+\gamma_4}{2}}\sqrt\frac{\log n}{n}$ with $\gamma_4\leq 2(\gamma_0+\gamma)$. Hence  when $d\leq c(\frac{n}{\log n})^{\frac{1}{1+2\gamma+2\gamma_2+4\gamma_0}}$ for a sufficiently small $c$, we can obtain that 
\begin{equation*}
 \begin{aligned}
  C_1\,d^{-\gamma_0}\|\widehat\theta-\widehat\theta^\diamond\|\leq\sqrt\frac{\log n}{n}d^{\frac{1+\gamma_3}{2}}\|\widehat\theta-\widehat{\theta}^{\diamond}\|^{\beta_1}+ d^{\frac{1+\gamma_3}{2}+\beta_1(\frac{1+\gamma_4}{2})}(\frac{\log n}{n})^{\frac{1+\beta_1}{2}}+ d^{1+\gamma\vee(\gamma_2+\gamma_4)}\frac{\log n}{n},
      \end{aligned}
 \end{equation*}
which leads to 
\begin{equation*}
  \|\widehat\theta-\widehat\theta^\diamond\|\leq C\,\bigg( d^{\frac{1+\gamma_3}{2}+\beta_1(\frac{1+\gamma_4}{2})+\gamma_0}(\frac{\log n}{n})^{\frac{1+\beta_1}{2}}+ d^{1+\gamma\vee(\gamma_2+\gamma_4)+\gamma_0}\frac{\log n}{n}+\Big(d^{\frac{1+\gamma_3}{2}+\gamma_0}\sqrt\frac{\log n}{n}\Big)^{\frac{1}{1-\beta_1}}\bigg).
\end{equation*}
Now we show equation~\eqref{normgdiamond}, using the first statement of Lemma~\ref{lemma1}, we can obtain that 
  \begin{equation*}
  \begin{aligned}
      &\bigg\|\frac{1}{n}\sum_{i=1}^n g(X_i,\widehat{\theta}^{\diamond})-\frac{1}{n}\sum_{i=1}^ng(X_i,\theta^\ast)-\mathbb{E}g(X,\widehat{\theta}^{\diamond})+\mathbb{E}g(X,\theta^\ast)\bigg\|\\
      &\leq C\, (\frac{\log n}{n})^{\frac{1+\beta_1}{2}}d^{\frac{1+\gamma_3}{2}+\beta_1(\frac{1+\gamma_4}{2})}+C\,\frac{\log n}{n}d^{1+\gamma}.
        \end{aligned}
   \end{equation*}
Moreover, by the Lipschitz continuity of $\m H_{\theta}$ around $\theta^\ast$, we can obtain that 
\begin{equation*}
    \begin{aligned}
     \|\mathbb{E}g(X,\widehat{\theta}^{\diamond})-\mathbb{E}g(X,\theta^\ast)-\m H_{\theta^\ast}(\widehat{\theta}^{\diamond}-\theta^\ast)\|&\leq d^{\gamma_2} \|\widehat{\theta}^{\diamond}-\theta^\ast\|^2\leq C\, d^{1+\gamma_4+\gamma_2}\frac{\log n}{n}.
    \end{aligned}
\end{equation*}
Therefore, combined with the fact that $\frac{1}{n}\sum_{i=1}^ng(X_i,\theta^\ast)+\m H_{\theta^\ast}(\widehat{\theta}^{\diamond}-\theta^\ast)=0$, we can obtain that it holds with probability at least $1-cn^{-2}$ that
\begin{equation*} 
    \begin{aligned}
 \left\|\frac{1}{n}\sum_{i=1}^ng(X_i,\widehat{\theta}^{\diamond})\right\|&=\left\|\frac{1}{n}\sum_{i=1}^ng(X_i,\widehat{\theta}^{\diamond})-\frac{1}{n}\sum_{i=1}^ng(X_i,\theta^\ast)-\m H_{\theta^\ast}(\widehat{\theta}^{\diamond}-\theta^\ast)\right\|\\
     &\leq \bigg\|\frac{1}{n}\sum_{i=1}^n g(X_i,\widehat{\theta}^{\diamond})-\frac{1}{n}\sum_{i=1}^ng(X_i,\theta^\ast)-\mathbb{E}g(X,\widehat{\theta}^{\diamond})+\mathbb{E}g(X,\theta^\ast)\bigg\|\\
     &+ \|\mathbb{E}g(X,\widehat{\theta}^{\diamond})-\mathbb{E}g(X,\theta^\ast)-\m H_{\theta^\ast}(\widehat{\theta}^{\diamond}-\theta^\ast)\|\\
     &\leq C\, d^{\frac{1+\gamma_3}{2}+\beta_1(\frac{1+\gamma_4}{2})}(\frac{\log n}{n})^{\frac{1+\beta_1}{2}}+C\, d^{1+\gamma\vee(\gamma_2+\gamma_4)}\frac{\log n}{n}.
    \end{aligned}
\end{equation*}
 
 \section{Proof of Remaining Results}
 \subsection{Proof of Lemma~\ref{lemma:Warm}}
 Let $\pi_{\rm loc}=[\sqrt{n}(\cdot-\wh\theta)]_{\#}\pi_n$ and $\mu_{\rm loc}=[\sqrt{n}(\cdot-\wh\theta)]_{\#}\mu_0$.
 We can bound
 \begin{equation*}
 \begin{aligned}
     M_0&= \underset{A\,:\, \pi_{n}(A)>0}{\sup}\frac{\mu_0(A)}{\pi_{n}(A)}\\
     &\overset{(i)}{=} \underset{A\subset K\,:\, \pi_{\rm loc}(A)>0}{\sup}\frac{\mu_{\rm loc}(A)}{\pi_{\rm loc}(A)}\\
     &=  \underset{A\subset K\,:\, \pi_{\rm loc}(A)>0}{\sup}\frac{\mu_{\rm loc}(A)}{\pi_{\rm loc}|_{K}(A)}\cdot \frac{1}{\pi_{\rm loc}(K)}\\
     &\leq \underset{x\in K}{\sup}\left[\frac{\int_K\exp(-\frac{1}{2}x^TJx)\,\dd x\exp(-\frac{1}{2}x^T\widetilde{I}^{-1}x)}{\int_K\exp(-\frac{1}{2}x^T\widetilde{I}^{-1}x)\,\dd x\exp(-\frac{1}{2}x^TJx)}
     \cdot\frac{\int_K\exp(-V_n(x))\,\dd x\exp(-\frac{1}{2}x^TJx)}{\int_K\exp(-\frac{1}{2}x^TJx)\,\dd x\exp(-V_n(x))}\right]\cdot \frac{1}{\pi_{\rm loc}(K)}\\
     &\leq \underset{x\in K}{\sup} \frac{\int_K\exp(-\frac{1}{2}x^TJx)\,\dd x\exp(-\frac{1}{2}x^T  \widetilde{I}^{-1}x)}{\int_K\exp(-\frac{1}{2}x^T  \widetilde{I}^{-1}x)\,\dd x\exp(-\frac{1}{2}x^TJx)}\cdot\underset{x\in K}{\sup}\frac{\int_K\exp(-V_n(x))\,\dd x\exp(-\frac{1}{2}x^TJx)}{\int_K\exp(-\frac{1}{2}x^TJx)\,\dd x\exp(-V_n(x))}\cdot \frac{1}{\pi_{\rm loc}(K)},
      \end{aligned}
 \end{equation*}
 where $(i)$ uses $\mu_{\rm loc}(K)=0$.
 Since for any function pair $f_1,f_2$, it holds that 
 \begin{equation*}
   \int_K f_1(x) \,\dd x\cdot \underset{x\in K}{\sup} \frac{f_2(x)}{f_1(x)}\geq  \int_K f_1(x)\frac{f_2(x)}{f_1(x)}\,\dd x=\int_K f_2(x) \,\dd x,
 \end{equation*}
 we can obtain that 
 \begin{equation*}
          M_0\leq \underset{x\in K}{\sup}  \exp(|x^T(  \widetilde{I}^{-1}-J)x|)\cdot\underset{x\in K}{\sup}\exp\big(2\big|V_n(x)-\frac{1}{2}x^TJx\big|\big)\cdot\frac{1}{\pi_{\rm loc}(K)}.
 \end{equation*}
   \subsection{Proof of Corollary~\ref{cor:smoothloss}}
   We first verify that under Condition B.3', Condition B.2 and Condition B.3 holds,  where the function $g$ in Condition B.3 is chosen as the gradient $\nabla_{\theta}\ell$. Condition B.2 and B.3.1 directly follows from the assumption that $\|\nabla_{\theta}\ell(x,\theta)\|\leq C d^{\gamma}$. For Condition  B.3.2,  since  $\mnorm{{\rm Hess}_{\theta}(\ell(x,\theta))}_{\rm op}^2\leq  Cd^{\gamma_3}$, we have  for any $x\in \m X$ and $\theta\in \Theta$,
   \begin{equation*}
       \|\nabla_{\theta}\ell(x,\theta)-\nabla_{\theta}\ell(x,\theta')\|\leq  \sqrt{Cd^{\gamma_3}} \|\theta-\theta'\|
   \end{equation*}
   and thus 
   \begin{equation*}
      d_n^g(\theta,\theta') \leq  \sqrt{Cd^{\gamma_3}} \|\theta-\theta'\|.
   \end{equation*}
   Then the covering  number condition for $d_n^g$  follows from the fact that the $\varepsilon$-covering number of unit $d$-ball is bounded by $(\frac{3}{\varepsilon})^d$. Condition B.3.3 directly follows from the assumption that $\mnorm{{\rm Hess}_{\theta}(\ell(x,\theta))}_{\rm op}^2\leq  Cd^{\gamma_3}$.  Condition B.3.4 follows from the assumption that $\m H_{\theta^*}^{-1}\Delta_{\theta^*}\m H_{\theta^*}^{-1}\preceq C\, d^{\gamma_4}I_d $ with $\Delta_{\theta^*}=\mb{E}[\nabla_{\theta} \ell(X,\theta^*)\\\nabla_{\theta} \ell(X,\theta^*)^T]$. Then the first statement directly follows from Theorem~\ref{th:Gibbsmixing}.
   For the second statement, we first verify that $\wt I^{-1}=|S|^{-1}\sum_{i\in S} {\rm Hess}_{\theta}(\ell(X_i,\wh\theta))$ is a reasonable estimator to $\m H_{\theta^*}$ in the following lemma.  
   
    \begin{lemma}\label{cor:pluginest}
Under assumptions in Corollary~\ref{cor:smoothloss}, let $m=|S|$, it holds with probability larger than $1-n^{-2}$ that $$\mnorm{\widetilde{I}^{-1}-\m H_{\theta^*}}_{\rm op}\leq C\, \big(d^{\frac{\gamma_3+1}{2}}\sqrt{\frac{\log n}{m}}\big)\vee \big(d^{\frac{\gamma_3+2}{2}}\frac{\log n}{m}\big)\vee \big(d^{\frac{1+\gamma_4}{2}+\gamma_2}\sqrt{\frac{\log n}{n}}\big).$$ 
 \end{lemma}
   Then since $\mnorm{\m H_{\theta^*}^{-1}}\leq C\, d^{\gamma_0}$, $d\leq c\frac{n^{\kappa_1}}{\log n}$ and $m\geq C_2\,d^{\gamma_3+2\gamma_0+\frac{7}{3}}$ , we have 
   \begin{equation*}
      \mnorm{\widetilde{I}}_{\rm op}\leq 2C d^{\gamma_0},
   \end{equation*}
and
   \begin{equation*}
   \begin{aligned}
       \mnorm{\wt I^{\frac{1}{2}}\m H_{\theta^*}\wt I^{\frac{1}{2}}-I_d}_{\rm op}&\leq \mnorm{\wt I}_{\rm op}\mnorm{\widetilde{I}^{-1}-\m H_{\theta^*}}_{\rm op}\\
       &\leq C\, d^{\gamma_0}\big(d^{\frac{\gamma_3+1}{2}}\sqrt{\frac{\log n}{m}}\big)\vee \big(d^{\frac{\gamma_3+2}{2}}\frac{\log n}{m}\big)\vee \big(d^{\frac{1+\gamma_4}{2}+\gamma_2}\sqrt{\frac{\log n}{n}}\big)\\
       &\leq \frac{1}{2},
          \end{aligned}
   \end{equation*}
   which leads to 
   \begin{equation*}
      \frac{1}{2} I_d \preceq\wt I^{\frac{1}{2}}\m H_{\theta^*}\wt I^{\frac{1}{2}}\preceq 2I_d.
   \end{equation*}
   Then by 
   \begin{equation*}
       \m H_{\theta^*}= \wt I^{-\frac{1}{2}} \big(\wt I^{\frac{1}{2}}\m H_{\theta^*}\wt I^{\frac{1}{2}}\big) \wt I^{-\frac{1}{2}},
   \end{equation*}
   we have 
   \begin{equation*}
       \begin{aligned}
          & \mnorm{\wt I}_{\rm op} \leq 2  \mnorm{\m H_{\theta^*}^{-1}}_{\rm op};\\
          &\mnorm{\wt I^{-1}}_{\rm op} \leq 2  \mnorm{\m H_{\theta^*}}_{\rm op}.
       \end{aligned}
   \end{equation*}
   Thus the requirements for the preconditioning matrix $\wt I$ in Theorem~\ref{th:Gibbsmixing} are satisfied with $\rho_2=2$ and $\rho_1=\frac{1}{2}$.  Finally, we will control the warming parameter using  Lemma~\ref{lemma:Warm}. Recall $\mu_0=N_d(\wh\theta,n^{-1}\wt I)\big|_{\{\theta:\sqrt{n}\wt I^{-\frac{1}{2}}(\theta-\wh\theta)\|\leq 3\sqrt{c_1 d}\}}$, where $c_1$  is a constant so that
   $c_1\geq 9\vee \underset{i\in [d],j\in [d]}{\sup} \frac{\partial^2\m R(\theta^*)}{\partial \theta_i\partial\theta_j} $. By 
   \begin{equation*}
       \mnorm{\wt I^{-\frac{1}{2}}}_{\rm op} \leq \sqrt{2}  \mnorm{\m H_{\theta^*}^{\frac{1}{2}}}_{\rm op}\leq \sqrt{2d \underset{i\in [d],j\in [d]}{\sup} \frac{\partial^2\m R(\theta^*)}{\partial \theta_i\partial\theta_j}}\leq  \sqrt{2c_1d},
   \end{equation*}
   and Lemma~\ref{lemmatail}, we can obtain that 
   \begin{equation*}
       \pi_n\Big(\sqrt{n}\|\wt I^{-\frac{1}{2}}(\theta-\wh\theta)\|\leq 2\sqrt{c_1 d}\Big)\geq 1-\exp(-1).
   \end{equation*}
   Moreover, consider $K=\{\xi\,:\,\wt I^{-\frac{1}{2}}\xi\leq 2\sqrt{c_1d}\}$, then for any $\xi \in K$, we have 
   \begin{equation*}
       \|\xi\|\leq 2\mnorm{\wt I^{\frac{1}{2}}}_{\rm op}\sqrt{c_1d}\leq 2\sqrt{2c_1d}\mnorm{\m H_{\theta^*}^{-\frac{1}{2}}}_{\rm op} \leq c_2 d^{\frac{1+\gamma_0}{2}}.
   \end{equation*}
   Then by Lemma~\ref{thmala}, when  $d\leq c\frac{n^{\kappa_1}}{\log n}$ for a small enough $c$, for any $\xi\in K$, we have
   \begin{equation*}
     \Big|V_n(\xi)-\frac{\xi^T\m H_{\theta^*}\xi}{2}\Big| \leq \frac{1}{2}.
   \end{equation*}
   In addition, for any $\xi \in K$, we have 
\begin{equation*}
\begin{aligned}
    \underset{\xi\in K}{\sup}{\big|\xi^T(\wt I^{-1}-\m H_{\theta^*})\xi\big|}&= \underset{\|\xi\|\leq 2\sqrt{c_1d}}{\sup}{\big|\xi^T(I_d-\wt{I}^{\frac{1}{2}}\m H_{\theta^*}\wt{I}^{\frac{1}{2}})\xi\big|}\\
    &\leq 2c_1d\mnorm{I_d-\wt{I}^{\frac{1}{2}}\m H_{\theta^*}\wt{I}^{\frac{1}{2}}}_{\rm op}\\
    &\leq 2c_1C \,d^{\gamma_0+1}\big(d^{\frac{\gamma_3+1}{2}}\sqrt{\frac{\log n}{m}}\big)\vee \big(d^{\frac{\gamma_3+2}{2}}\frac{\log n}{m}\big)\vee \big(d^{\frac{1+\gamma_4}{2}+\gamma_2}\sqrt{\frac{\log n}{n}}\big)\\
    &\leq d^{\frac{1}{3}},
    \end{aligned}
\end{equation*}
 where the last inequality uses $d\leq c\frac{n^{\kappa_1}}{\log n}$ and $m\geq C_2\,d^{\gamma_3+2\gamma_0+\frac{7}{3}}$.  The desired result then follows from Lemma~\ref{lemma:Warm}.
   
 \subsection{Proof of Lemma~\ref{cor:pluginest}} 
 Since $\mb E[\widetilde{I}^{-1}]=\m H_{\wh\theta}$, we have
 \begin{equation*}
     \mnorm{\widetilde{I}^{-1}-\m H_{\theta^*}}_{\rm op}\leq \mnorm{\widetilde{I}^{-1}-\m H_{\wh\theta}}_{\rm op}+\mnorm{\m H_{\wh\theta}-\m H_{\theta^*}}_{\rm op}.
 \end{equation*}
 The second term can be bounded using Condition B.1.2 and equation~\eqref{diffhattheta} in the proof of Lemma~\ref{th1}, that is 
 \begin{equation*}
 \begin{aligned}
 \mnorm{\m H_{\wh\theta}-\m H_{\theta^*}}_{\rm op}\leq C\, d^{\gamma_2}\|\wh\theta-\theta^*\|\leq C\, d^{\frac{1+\gamma_4}{2}+\gamma_2}\sqrt{\frac{\log n}{n}}.
      \end{aligned}
 \end{equation*}
 The first term can be bounded using Bernstein's inequality. Let $m=|S|$, for $v,v'\in \mb S^{d-1}$ and $\theta,\theta'\in B_r(\theta^*)$,  we have
 \begin{equation}\label{pseudometric}
  \begin{aligned}
 &\sqrt{m^{-1}\sum_{i\in S} \big(v^T{\rm Hess}_{\theta}(\ell(X_i,
 \theta))v-{v'}^T{\rm Hess}_{\theta}(\ell(X_i,
 \theta'))v'\big)^2}\\
 &\leq C\, \sqrt{d} \|v-v'\|+ C\, d^{r_1} \|\theta-\theta'\|.
  \end{aligned}
 \end{equation}
 Then consider $\m N_{v}$ and $\m N_{\theta}$ to be the minimal $n^{-1}$ and $n^{-1}d^{-r_1}$ covering set of $\mb S^{d-1}$ and $B_r(\theta^*)$, then $\log |\m N_v|\leq C\, d\log n$ and $\log |\m N_\theta|\leq C\, d\log n$. Using the fact that 
 \begin{equation*}
     \underset{\theta\in B_r(\theta^*),\, X\in \m X}{\sup} \mnorm{{\rm Hess}_{\theta}(\ell(X,
 \theta))}_{\rm op}\leq C\, d^{\frac{\gamma_3}{2}};
 \end{equation*}
 \begin{equation*}
     \underset{\theta\in B_r(\theta^*)\, v\in \mb S^{d-1}}{\sup}\mb E \big[(v^T{\rm Hess}_{\theta}(\ell(X,\theta))^2\big]\leq \underset{\theta,\theta'\in B_r(\theta^*), \atop v\in \mb S^{d-1}}{\sup}\mb E \Big[\frac{(v^T\nabla\ell(X,\theta)-v^T\nabla\ell(X,\theta'))^2}{\|\theta-\theta'\|^2}\Big]
\leq C\, d^{\gamma_3},
 \end{equation*}
we can get by Bernstein's inequality and a simple union bound argument
that it holds with probability at least $1-n^{-c}$ that for any $v\in \m N_v$ and $\theta\in \m N_\theta$,
\begin{equation*}
    \underset{\theta\in B_r(\theta^*)}{\sup}\underset{v\in\mb S^{d-1}}{\sup}  \left(v^T(m^{-1}\sum_{i\in S} {\rm Hess}_{\theta}(\ell(X_i,
 \theta))-\m H_{\theta})v^T\right)  \leq C\, \big(d^{\frac{\gamma_3+1}{2}}\sqrt{\frac{\log n}{m}}\big)\vee \big(d^{\frac{\gamma_3+2}{2}}\frac{\log n}{m}\big).
 \end{equation*}

 \subsection{Proof of Corollary~\ref{co:quantile}}
 We will first check that Conditions B.1-B.3 hold for the quantile regression example under Condition D.1 and D.2. Consider the loss function 
\begin{equation*}
    \ell(X,\theta)=(Y-\widetilde X^T\theta)(\tau-\mathbf{1}(Y<\widetilde{X}^T\theta)),
\end{equation*}
 and its subgradient
 \begin{equation*}
  g(X,\theta)=(\mathbf{1}(Y<\widetilde{X}^T\theta)-\tau)\widetilde X.
\end{equation*}
Then we can write 
\begin{equation*}
    \m R(\theta)=\mb E[\ell(X,\theta)]=\mb E \big[\tau\,(Y-\widetilde X^T \theta)\big]-\mb E  \Big[\int_{-\infty}^{\widetilde X^T\theta-\widetilde X^T\theta^\ast}(\varepsilon+\widetilde{X}^T\theta^\ast-\widetilde X^T\theta) f_{e}(\varepsilon)d\varepsilon\Big].
\end{equation*}
Taking derivative of $\m R$ w.r.t $\theta$, we can obtain
\begin{equation*}
    \nabla \m R(\theta)= -\tau \cdot\mb E [\widetilde X]+\mb E [\mathbf{1}(Y<\widetilde{X}^T\theta)\widetilde X] =\mb E g(X,\theta).
\end{equation*}
Thus,
\begin{equation*}
    \m H_{\theta}=\mb E [f_{e}(\widetilde X^T\theta-\widetilde X^T\theta^\ast)\widetilde X\widetilde X^T].
\end{equation*}
Then for $\theta\in B_{c/\sqrt{d}}(\theta^\ast)$ with a small enough $c$, it holds that 
\begin{equation*}
    \frac{f_{e}(\widetilde X^T\theta-\widetilde X^T\theta^\ast)}{f_{e}(0)}\geq \frac{1}{2}.
\end{equation*}
Then by the fact that $\nabla \m R(\theta^\ast)=0$ and $\mb E[\widetilde X\widetilde X^T]\succeq C'd^{-\alpha_0}I_d$, we can obtain that  for any $\theta\in B_{\frac{c}{\sqrt{d}}}(\theta^\ast)$,
\begin{equation*}
    \m R(\theta)-\m R(\theta^\ast)\geq C_1\, d^{-\alpha_0}\|\theta-\theta^\ast\|^2;
\end{equation*}
on the other hand, for any $\theta\in  B_{\frac{c}{\sqrt{d}}}(\theta^\ast)^c$, 
\begin{equation*}
    \m R(\theta)-\m R(\theta^\ast)\geq \m R\Big(\theta^\ast+\frac{c(\theta-\theta^\ast)}{\sqrt{d}\|\theta-\theta^\ast\|}\Big)-\m R(\theta^\ast)\geq C_1\, d^{-\alpha_0-1},
\end{equation*}
hence for any $\theta\in \mb R^d$,
\begin{equation*}
    \m R(\theta)-\m R(\theta^\ast)\geq C_1\, d^{-\alpha_0}(d^{-1}\wedge \|\theta-\theta^\ast\|^2).
\end{equation*}
 Moreover, for any $\theta\in \Theta$ and $v\in \mb S^{d-1}$, 
 \begin{equation*}
     \begin{aligned}
      |v^T (\m H_{\theta}-\m H_{\theta^\ast})v|&\leq v^T\mb E\left[\left|f_{e}(\widetilde X^T\theta-\widetilde X^T\theta^\ast)-f_{e}(0)\right|\widetilde X\widetilde X^T\right]v\\
      &\leq C\, \mb E \left[|\widetilde{X}^T(\theta-\theta^\ast)|v^T\widetilde X \widetilde X^T v\right]\\
      &\leq C\, \|\theta-\theta^\ast\| \mb E\left(\left|\widetilde X  {(\theta-\theta^\ast)}/{\|\theta-\theta^\ast\|}\right|^3\right)^{\frac{1}{3}} (\mb E |v^T \widetilde X|^3)^{\frac{2}{3}}\\
      &\leq C\, d^{\alpha_1}\|\theta-\theta^\ast\|,
     \end{aligned}
 \end{equation*}
where the last inequality uses the assumption that $\sup_{\eta\in \mb S^{d-1}}\mb E[\eta^T \widetilde X]\leq C d^{\alpha_1}$. Thus we have Condition B.1 holds with $\gamma_0=\alpha_0$, $\gamma_1=1$, $\gamma_2=\alpha_1$. For Condition B.2, by $\m X={\rm supp}(\widetilde X)\subseteq [-C,C]^d$, we can obtain $\|g(X,\theta)\|\leq C\sqrt{d}$, thus for any $\theta,\theta'$, $|\ell(X,\theta)-\ell(X,\theta')|\leq C\sqrt{d}\|\theta-\theta'|$ and Condition B.2 and Condition B.3.1 hold with $\gamma=\frac{1}{2}$. For Condition B.3,  since for any $\theta,\theta'\in \Theta$, 
\begin{equation*}
\begin{aligned}
 \sqrt{\frac{1}{n}\sum_{i=1}^n \|g(X_i,\theta)-g(X_i,\theta')\|^2}&=\sqrt{\frac{1}{n}\sum_{i=1}^n \|\widetilde{X}_i\|^2(\mathbf{1}(Y<\widetilde X_i^T\theta)-\mathbf{1}(Y<\widetilde X_i^T\theta'))^2}\\
 &=\sqrt{d}\sqrt{\frac{1}{n}\sum_{i=1}^n (\mathbf{1}(Y<\widetilde X_i^T\theta)-\mathbf{1}(Y<\widetilde X_i^T\theta'))^2},
    \end{aligned}
\end{equation*}
 by Lemma 9.8 and Lemma 9.12 of~\cite{Kosorok2008},  the function class $\mathcal{F}=\{\mathbf{1}(Y\leq \theta^T\widetilde{X}), \theta\in \Theta\}$ is a VC-class with VC-dimension being bouned by $d+3$, then using Theorem 8.3.18 of~\cite{vershynin_2018} on the covering number's upper bound via VC dimension, we can verify Condition B.3.2. 

For Condition  B.3.3,   since for any $v\in \mb S^{d-1}$ and $\theta,\theta'\in \Theta$, 
\begin{equation*}
    \begin{aligned}
     \mb E(v^T g(X,\theta)-v^Tg(X,\theta'))^2&=\mb E[(\mathbf{1}(Y<\widetilde X_i^T\theta)-\mathbf{1}(Y<\widetilde X_i^T\theta'))^2(v^T\widetilde X)^2]\\
     &=\mb E \left[(v^T\widetilde X)^2 \int_{\widetilde X^T \theta\wedge \widetilde X^T\theta'}^{\widetilde X^T \theta\vee \widetilde X^T\theta'}  f(y-\widetilde{X}^T\theta^\ast|\widetilde X)\,\dd y\right]\\
     &\leq C\, \mb E \left[(v^T\widetilde X)^2 |\widetilde X^T\theta-\widetilde X^T\theta']\right]\\
     &\leq C\, \|\theta-\theta'\|\underset{v\in \mb S^{d-1}}{\sup}\mb E|v^T\widetilde X|^3\leq C\, d^{\alpha_1}\|\theta'-\theta\|;
    \end{aligned}
\end{equation*}
\begin{equation*}
    \begin{aligned}
     &\mb E \big[(\ell(X,\theta)-\ell(X,\theta')-g(X,\theta')(\theta-\theta'))^2\big]\\
     &=\mb E\big[\big(-(Y-\widetilde X^T\theta)\mathbf{1}(Y<\widetilde X_i^T\theta)+(Y-\widetilde X^T\theta')\mathbf{1}(Y<\widetilde X_i^T\theta')-\mathbf{1}(Y<\widetilde X_i^T\theta') \widetilde X^T(\theta-\theta')\big)^2\big]\\
     &=\mb E \left[\int_{\widetilde X^T \theta\wedge \widetilde X^T\theta'}^{\widetilde X^T \theta\vee \widetilde X^T\theta'} (y-\widetilde{X}^T\theta)^2 f(y-\widetilde{X}\theta^\ast|\widetilde X)\,\dd y\right]\\
     &\leq C\, \mb  E|\widetilde X^T\theta-\widetilde X^T\theta'|^3\leq C\, d^{\alpha_1}\|\theta'-\theta\|^3.
    \end{aligned}
\end{equation*}
Thus Condition B.3.3 holds with $\gamma_3=\alpha_1$ and $\beta_1=\frac{1}{2}$. For condition B.3.4, since
\begin{equation*}
    \mb E[g(X,\theta^\ast)g(X,\theta^\ast)^T]=\mb {E}\big[(\tau^2+\mathbf{1}(Y<\widetilde{X}^T\theta)-2\tau \mathbf{1}(Y<\widetilde{X}^T\theta))\widetilde X\widetilde X^T\big]=(\tau-\tau^2)\mb E[\widetilde X\widetilde X^T],
\end{equation*}
 and $J=\m H_{\theta^\ast}=f_{e}(0)\mb {E}[\widetilde X\widetilde X^T]$, we have 
 \begin{equation*}
     (\mb E[\widetilde X\widetilde X^T])^{\frac{1}{2}}J^{-1}(\mb E[\widetilde X\widetilde X^T])^{\frac{1}{2}}=f_{e}(0)^{-1} I_d,
 \end{equation*}
and thus $\gamma_4=\gamma_0$.  

\quad\\
Now we verify that the requirements of the $\wt I$ in Theorem~\ref{th:Gibbsmixing} are satisfied. Recall $\widetilde I^{-1}=\frac{1}{|S|}\sum_{i\in S} X_iX_i^T$, in order to show that
 $\mnorm{ \widetilde{I}^{-\frac{1}{2}}J^{-1} \widetilde{I}^{-\frac{1}{2}}}_{  \rm  op}\vee \mnorm{ \widetilde{I}^{\frac{1}{2}}J \widetilde{I}^{\frac{1}{2}}}_{  \rm  op}$  is bounded above by a constant, we will derive upper bound to the term of $\mnorm{  \widetilde{I}^{\frac{1}{2}}(\mb{E}[\wt X\wt X^T]) \widetilde{I}^{\frac{1}{2}}- I_d}_{  \rm  op}$. Let $m=|S|$, similar as the proof for Lemma~\ref{cor:pluginest},  we can obtain it holds with probability larger than $1-\frac{1}{n^2}$ that 
\begin{equation*}
\begin{aligned}
     \bbmnorm{ n^{-1}\sum_{i=1}^n \widetilde X_i\widetilde X_i^T-\mb E[\widetilde X\widetilde X^T]}_{\rm  op}&\leq C\,\underset{v\in \mb S^{d-1}}{\sup}\sqrt{\mb {E}|v^T\widetilde{X}|^4} d^{\frac{1}{2}}\sqrt{\frac{\log n}{m}}+d^2\frac{\log n}{m}\\
     \leq C\, d^{\frac{3}{4}+\frac{\alpha_1}{2}}\sqrt{\frac{\log n}{m}}+d^2\frac{\log n}{m},
\end{aligned}
\end{equation*}
where the last inequality is due to   $\underset{v\in \mb S^{d-1}}{\sup}\sqrt{\mb {E}|v^T\widetilde{X}|^4}\leq C\, {d}^{\frac{1}{4}}\underset{v\in \mb S^{d-1}}{\sup}\sqrt{\mb {E}|v^T\widetilde{X}|^3}\leq C\, d^{\frac{1+2\alpha_1}{4}}$. Then by $\mb E [\widetilde X\widetilde X^T]\succeq C'd^{-\alpha_0}I_d$, and $m\geq C_2\, d^{\alpha_1+2\alpha_0+3/2}\log n$, we can obtain 
\begin{equation*}
   \mnorm{\wt I}_{\rm op}\leq \frac{2}{C'}d^{\alpha_0}
\end{equation*}
Thus we have 
   \begin{equation*}
   \begin{aligned}
       \mnorm{\wt I^{\frac{1}{2}}(\mb{E}[\wt X\wt X^T]) \wt I^{\frac{1}{2}}-I_d}_{\rm op}\leq \mnorm{\wt I}_{\rm op}\mnorm{\widetilde{I}^{-1}-(\mb{E}[\wt X\wt X^T]) }_{\rm op}\leq C_1 \, d^{\alpha_0+\frac{3+2\alpha_1}{4}}\sqrt{\frac{\log n}{m}},
          \end{aligned}
   \end{equation*}
   which leads to 
   \begin{equation*}
      \frac{1}{2} I_d \preceq\wt I^{\frac{1}{2}} (\mb{E}[XX^T])\wt I^{\frac{1}{2}}\preceq 2I_d,
   \end{equation*}
   Thus 
   \begin{equation*}
       \frac{1}{2} f_{e}(0)I_d \preceq\wt I^{\frac{1}{2}}\m H_{\theta^*}\wt I^{\frac{1}{2}}\preceq 2f_{e}(0)I_d
   \end{equation*}
   Furthermore, by 
   \begin{equation*}
       \m H_{\theta^*}=f_{e}(0)\cdot \wt I^{-\frac{1}{2}} \big(\wt I^{\frac{1}{2}}(\mb{E}[XX^T])\wt I^{\frac{1}{2}}\big) \wt I^{-\frac{1}{2}},
   \end{equation*}
   we have 
   \begin{equation*}
       \begin{aligned}
          & \mnorm{\wt I}_{\rm op} \leq 2f_{e}(0)  \mnorm{\m H_{\theta^*}^{-1}}_{\rm op};\\
          &\mnorm{\wt I^{-1}}_{\rm op} \leq \frac{2}{f_{e}(0) } \mnorm{\m H_{\theta^*}}_{\rm op}.
       \end{aligned}
   \end{equation*}
We can then obtain that the requirements for the preconditioning matrix $\wt I$ in Theorem~\ref{th:Gibbsmixing} are satisfied with $\rho_2=2{f_{e}(0)}$ and $\rho_1=\frac{1}{2}f_{e}(0)$.

\end{document}